\documentclass[11pt,letterpaper]{article}
\usepackage[margin=1in]{geometry}
\usepackage{amsmath,amssymb,amsthm,booktabs,tabularx,enumitem}
\usepackage{mathtools,needspace,graphicx,float}
\usepackage[OT1]{fontenc}

\usepackage[expansion=true,protrusion=true]{microtype}
\usepackage{setspace,etoolbox}
\usepackage[authoryear,round,sort]{natbib}
\usepackage{xr-hyper}
\usepackage[dvipsnames]{xcolor}
\usepackage{tcolorbox}
\usepackage[colorlinks=true,linkcolor=MidnightBlue,citecolor=MidnightBlue,urlcolor=MidnightBlue]{hyperref}
\AtBeginEnvironment{figure}{\onehalfspacing}
\AtBeginEnvironment{table}{\onehalfspacing}
\def\spacingset#1{\renewcommand{\baselinestretch}%
{#1}\small\normalsize} \spacingset{1}
\theoremstyle{plain}
\newtheorem{theorem}{Theorem}
\newtheorem{lemma}{Lemma}
\newtheorem{proposition}{Proposition}
\newtheorem{corollary}{Corollary}
\theoremstyle{definition}
\newtheorem{definition}{Definition}
\newtheorem{assumption}{Assumption}
\theoremstyle{remark}
\newtheorem{remark}{Remark}
\newtheorem{example}{Example}
\newenvironment{keywords}{\par\noindent\textbf{Keywords.}\enspace\ignorespaces}{\par}

\newcommand{\RR}{\mathbb{R}}
\newcommand{\EE}{\mathbb{E}}
\newcommand{\PP}{\mathbb{P}}
\newcommand{\cM}{\mathcal{M}}
\newcommand{\cN}{\mathcal{N}}
\newcommand{\cS}{\mathcal{S}}
\newcommand{\cU}{\mathcal{U}}

\newcommand{\cI}{\mathcal{I}}
\newcommand{\cQ}{\mathcal{Q}}
\newcommand{\T}{\intercal}
\newcommand{\ev}{{\rm ev}}
\newcommand{\norm}[1]{||#1||}

\title{Identifiability of Latent Variable Models with Flexible Latent Distributions}
\author{Chengyu Cui and Gongjun Xu\\\vspace{0.5em} {Department of Statistics, University of Michigan}}
\date{}
\makeatletter
\renewcommand{\maketitle}{%
  \par\begingroup
  \centering\setstretch{1}\LARGE {\bf \@title}\par
  \vspace{1.5em}\large\@author\par
  \endgroup\vspace{2em}%
}
\makeatother
\hypersetup{pdftitle={Identifiability of Latent Variable Models with Flexible Latent Distributions},pdfauthor={Chengyu Cui and Gongjun Xu},pdfsubject={},pdfkeywords={},pdfcreator={LaTeX}}

\AtBeginDocument{\nocite{anandkumar2014tensor,deville1992calibration,kolda2009tensor,ostrovskii2021finite,qin1994empirical,rota1994classical,stein2003complex,withers1984asymptotic}}

\begin{document}
\maketitle

\spacingset{1.3}
\begin{abstract}

        Latent variable models are popularly used to characterize dependence among many observed variables through underlying latent variables.
    Recent advances in statistics and machine learning methods have encouraged growing interest in studying latent variable models with flexible latent distributions that capture complex latent structures. 
    This flexibility, however, raises a fundamental question of whether model parameters can be identified when the latent distribution is unknown and must itself be learned.
    In this paper, we show that exact identifiability fails almost everywhere, in the sense that distinct but arbitrarily close pairs of model parameters and latent distributions can yield exactly the same distribution for the observed variables. Nonetheless, the resulting identification ambiguity is small. In particular, we prove that the identification ambiguity, measured by the diameter of the set of parameters inducing the same marginal distribution, vanishes at a superpolynomial rate of $\exp\{-c(\log p)^2\}$ as the number of observed variables $p$ grows. Building on this result, we show that the failure of exact identifiability need not preclude consistent estimation or valid inference of the model parameters and latent distribution. These results suggest that, when many variables are observed, the identifiability issue is far less of an obstacle to estimation and inference in latent variable models than is commonly expected.
\end{abstract}

\section{Introduction}\label{sec_intro}


Advances in machine learning have greatly expanded our ability to learn
complex distributions~\citep{kingma2013auto,goodfellow2014generative,papamakarios2019flows,ho2020denoising}. These advances have, in turn, stimulated interest in flexible latent distributions within latent variable models, where the distribution of the observed variables is determined jointly by the conditional response model and the latent distribution. Moving beyond conventional Gaussian latent distributions allows such models to capture richer forms of population heterogeneity, and thereby to support inference on complex latent structures across the social, economic, and health sciences~\citep{lubke2005investigating,arellano2009robust,antonelli2016mitigating,xu2023identifiable}.
At the same time, these developments raise a fundamental question of whether model parameters remain identifiable when the latent distribution is left unrestricted, a property on which estimation, inference, and interpretation all depend.


This identifiability question is especially relevant under nonlinear models, where the model parameters and the latent distribution jointly affect the observed responses through a nonlinear relationship.
As a representative case, we consider a binary response vector $X\in\{0,1\}^{p}$, whose distribution depends on a $d$-dimensional latent random vector $Z\in\RR^d$. 
Conditional on $Z$, the components $\{X_j\}_{j\le p}$ are independent, with
\begin{equation}\label{eq_intro_model}
\mathbb P(X_j=1\mid Z=z) \; = \; f(\beta_j+\alpha_j^{\mathsf T}z), \quad\text{for}\quad j=1,\ldots,p,
\end{equation}
where $f:\mathbb R\to(0,1)$ is a specified response function. The model parameters are denoted by $\theta=(A,\beta)$, where $A = (\alpha_1,\dots,\alpha_p)^\T\in\RR^{p\times d}$ and $\beta=(\beta_1,\ldots,\beta_p)^{\T}\in\RR^p$. Since $Z$ is unobserved, the distribution of $X$ is obtained by integrating out \(Z\) with respect to its latent distribution $G$. Specifically, for $x\in\{0,1\}^p$, 
\begin{equation}
p_{\theta,G}(X = x) \;=\; \int_{\mathbb R^d} \prod_{j=1}^{p} \big\{f(\beta_j+\alpha_j^{\mathsf T}z)\big\}^{x_j} \big\{1-f(\beta_j+\alpha_j^{\mathsf T}z)\big\}^{1-x_j} \,dG(z).\label{eq_marginal_law}
\end{equation}
We denote the corresponding marginal distribution of the observed variable $X$ by \(P_{\theta,G}\), which depends on both the model parameter and latent distribution.
This model covers several important classes of latent models:

\begin{enumerate}[leftmargin=0.6cm]
    \item \textbf{Item response theory.}
 In item response theory and related measurement models for educational and psychological assessment, $Z$ represents unobserved abilities, attitudes, or psychological traits~\citep{reckase2009multidimensional}, and its distribution $G$ can be asymmetric, multimodal, and dependent to capture rich population heterogeneity~\citep{skrondal2004generalized}.
    Parameter \(\alpha_j\) describes how strongly and in which directions item \(j\) responds to variation in these latent traits.

    \item \textbf{Latent-class models.}
    When $G$ has finite support, \eqref{eq_intro_model} becomes a latent class model given by a finite mixture of product Bernoulli distributions~\citep{goodman1974exploratory,lindsay1991semiparametric,allman2009identifiability}. The support points of $G$ represent latent classes with distinct response profiles, while their probability masses give the corresponding population proportions~\citep{collins2013latent}. Here, $\alpha_j$ describes how the response probability for item $j$ varies across these latent classes.

\item \textbf{Generative models.}
Under the logistic link $f(\cdot)$, for instance, model~\eqref{eq_intro_model} describes the Bernoulli observation layer used in sigmoid belief networks and related generative models~\citep{hinton1995wake,kingma2013auto,tomczak2018vae}.
In this setting, $Z$ is the hidden representation, $G$ is its distribution, and $(A,\beta)$ parameterizes the observation layer.
\end{enumerate}


A fundamental issue shared by these applications is model identifiability.
Formally, identifiability requires that, for any $(\theta,G)$ and $(\tilde\theta,\tilde G)$ in the admissible parameter space, \begin{equation*}
P_{\theta,G}=P_{\tilde\theta,\tilde G}
    \quad\Longrightarrow\quad
    \theta=\tilde\theta.
\end{equation*}
Identifiability plays a fundamental role in both the theory and practice of latent variable modeling. Theoretically, it is a prerequisite for consistent parameter estimation and valid statistical inference. Practically, without identifiability, the data cannot distinguish between parameters that induce the same observed distribution, leaving the interpretation of the estimated model ambiguous.
We therefore ask
 \begin{center}\hypertarget{Q1}{(Q1)} \ {\it Are the model parameter $\theta$ and latent distribution $G$ identifiable in \eqref{eq_intro_model}?}\\ \hypertarget{Q2}{(Q2)} \ {\it If not, is reliable estimation still possible?}\end{center}

The primary difficulty in answering this is to uniquely disentangle these two components from the marginal distribution $P_{\theta,G}$. Since $G$ enters the marginal distribution via integration, the resulting ambiguity is of a different and less tractable nature than the familiar indeterminacies in $\theta$. To make the problem more tractable, most existing studies examine identifiability under additional restrictions, typically on $G$.
For instance, 
conventional analyses often take $G$ to be a $d$-dimensional Gaussian distribution, in which case the identifiability issue reduces to affine indeterminacies that can be resolved through appropriate constraints~\citep{goodman1974exploratory,andersonintroduction}. However, a Gaussian specification may fail to capture asymmetry, multimodality, and other features of heterogeneous latent populations~\citep{zhang2001linear,lubke2005investigating}.

Beyond the Gaussian case, several studies allow non-Gaussian latent distributions under other restrictions. For example, \citet{allman2009identifiability} consider latent distributions with a known finite number of support points and establish generic identifiability. 
Restrictions on the response structure are also considered, for instance by~\citet{xu2017} and~\citet{gu2020partial}, to establish identifiability of restricted latent class models.
In deep generative models, identifiability is often achieved using auxiliary variables, such as time segments or class labels~\citep{hyvarinen2016unsupervised,khemakhem2020variational}, or mixture priors on the latent variables~\citep{kivva2022identifiability}. 
Across these studies, identifiability is established under restrictions that can be difficult to verify in practice.  Moreover,  these restrictions narrow the model class, leaving open whether identifiability holds when the latent distribution $G$ is left essentially unrestricted.

In this work, we show that exact identifiability fails for almost every parameter in model~\eqref{eq_intro_model}, answering \hyperlink{Q1}{(Q1)} in the negative. This is perhaps unsurprising in light of existing non-identifiability results for the one-dimensional item response theory model~\citep{haberman2005identifiability,san2011bayesian}, yet our result is considerably stronger.
For almost every parameter $\theta$ and latent distribution $G$, we construct alternatives $(\tilde\theta,\tilde G)$ satisfying $P_{\theta,G}=P_{\tilde\theta,\tilde G}$, with $\tilde\theta\neq \theta$ and $\tilde G$ matching the moments of $G$ up to any prespecified order. This reveals a fundamental obstacle to identifiability under flexible latent distributions that finitely many moment conditions on $G$ cannot remove.


The more surprising finding is the superpolynomial rate at which this ambiguity vanishes as $p$ grows.
Specifically, under suitable regularity conditions, the discrepancy $D_p(\theta,\tilde\theta)$ between parameters $\theta,\tilde\theta$ (defined in \eqref{eq_def_Dp} below) satisfies
\begin{equation*}
    D_p(\theta,\tilde\theta) \; \le \; C\exp\{-c(\log p)^2\} \quad\text{whenever}\quad P_{\theta,G} \, = \, P_{\tilde\theta,\tilde G},
\end{equation*}
uniformly over admissible $(\theta,G)$, $(\tilde\theta,\tilde G)$, for constants $C,c>0$. 
Thus, the ambiguity decays faster than any polynomial rate of $p$ and becomes negligible on every polynomial scale. 
The key technical challenge is to bound $D_p(\theta,\tilde\theta)$ uniformly over all parameter pairs that yield the same marginal distribution. We address this challenge by constructing high-order approximations based on conditional moments of some carefully selected observed responses, given the rest of the observed responses. The nonlinear link leads to complicated interactions among conditional latent moments of different orders, requiring careful control of the expansion coefficients and remainder. Choosing an expansion order proportional to $\log p$ balances the decay of the remainder against the growth of the coefficients in the expansion, yielding the $\exp\{-c(\log p)^2\}$ bound.


This superpolynomial rate provides a basis for answering the second question \hyperlink{Q2}{(Q2)}. In particular, we show that, even without exact identifiability, estimation errors for both the model parameters and the latent distribution can converge to zero as the sample size and $p$ diverge. We also provide conditions of valid inference for functionals of the model parameters and population-level quantities defined by the latent distribution.
These results provide theoretical justification for the estimability of the model parameters and latent distribution when they are learned jointly~\citep[e.g.][]{woods2006item,antonelli2016mitigating,paganin2023computational,li2026estimating,cui2026multidimensional}, in the sense that non-identifiability need not pose an irreducible identifiability barrier to accurate estimation when $p$ is large.

\subsection{Related Work}
The vanishing identification ambiguity described above can be termed asymptotic identifiability. For unidimensional nonparametric item response models, \citet{douglas2001asymptotic} established this property for the response functions, and \citet{he2024extended} extended the analysis by relaxing derivative bounds. Both analyses rely on transforming the latent variable to a uniform distribution, which is not applicable in multidimensional settings, and neither quantifies the rate of the vanishing ambiguity.
Related questions are studied in panel data models, where $p$ is the panel length. For instance, \citet{chernozhukov2013average} studied the identified sets of average and quantile effects in panel models and showed that they shrink to a point as the panel length increases. In fixed-effects logit models, \citet{davezies2021identification} derived sharp bounds on average causal effects.

A related literature examines estimation error arising from misspecification of the latent distribution. For logistic mixed-effects models, \citet{neuhaus1992effects} showed that misspecifying the mixing distribution can lead to inconsistent parameter estimates. For longitudinal models,
\citet{rizopoulos2008shared} showed that the effects of such misspecification vanish for certain parameters as the number of measurements increases. For nonlinear panel models, \citet{arellano2009robust} analysed estimation bias under misspecification of the individual-effects. For generalized linear latent variable models, \citet{ma2010explicit} developed semiparametric estimating equations under misspecified latent distributions.
Our work studies a different identification problem when the latent distribution is unknown and must be estimated. The resulting identification ambiguity  is intrinsic to the model, independent of any specific estimation procedure.

Another related literature develops estimation methods under flexible latent distributions, including spline-based densities~\citep{woods2006item}, Dirichlet-process mixtures~\citep{paganin2023computational}, kernel density methods~\citep{li2026estimating}, and normalizing flows~\citep{cui2026multidimensional}, among others. These studies report accurate  estimation of the model parameters when the latent distribution is fitted jointly, but do not provide statistical guarantees for the resulting estimators.
This gap cannot be addressed by classical theory for latent variable models, which generally presupposes that the model parameters are identifiable~\citep{lindsay1991semiparametric,ma2010explicit}, a condition that fails for model~\eqref{eq_intro_model} when $G$ is unrestricted.


Identifiability is also a central concern in generative modelling and representation learning, where the goal is to uncover the underlying generative mechanisms or latent representations, under a flexible latent distribution.
Such questions arise in unsupervised disentanglement~\citep{locatello2019challenging}, identifiable generative models~\citep{khemakhem2020variational,xi2023indeterminacy}, and causal representation learning~\citep{scholkopf2021toward}. A prevailing view in this literature is that complex latent structures cannot be identified without additional information or structural constraints, such as auxiliary variables or inductive biases~\citep{locatello2019challenging,khemakhem2020variational}. Notably, \citet{xi2023indeterminacy} characterized the indeterminacies of generative models with flexible nonlinear generators. Their framework requires a continuous data-generating process and thus excludes discrete observations. Our work offers another route that accepts the failure of exact identifiability and instead quantifies the ambiguity that remains.


\subsection{Notation}
For a positive integer $n$, let $[n]=\{1,\ldots,n\}$. For a finite set $S$, $|S|$ denotes its cardinality, and $x_S$ is the subvector of $x$ indexed by $S$. For a vector $x$, $\|x\|$ and $\|x\|_\infty$ denote its $\ell_2$ and $\ell_{\infty}$ norms, respectively.
For a matrix $M=(M_{ij})_{n\times m}$, write $\|M\|_\infty=\max_{i\in[n],j\in[m]}|M_{ij}|$ and $\|M\| = \max_{\|x\|=1}\|Mx\|$, and write ${\rm col}(M)$ for its column space. Let  $I_d$ be the $d\times d$ identity matrix.
The class $C^\infty$ consists of functions with continuous derivatives of all orders.
For a multi-index $a=(a_1,\ldots,a_d)\in\{0,1,\ldots\}^d$, write $|a|=\sum_{k=1}^d a_k$ and $z^a=\prod_{k=1}^d z_k^{a_k}$.
Unless otherwise stated, $c,C>0$ denote constants that may change from line to line.

\section{Problem Setup}\label{sec_setup}
In this section, we formalize the model identifiability for~\eqref{eq_intro_model}. The model suffers from the  well-known affine indeterminacy as in latent factor models~\citep{andersonintroduction}. Specifically, for any $\mu\in\mathbb R^d$ and invertible $H\in\mathbb R^{d\times d}$, let $\tilde G$ denote the distribution of $\tilde Z=H^{-1}(Z-\mu)$ when $Z\sim G$, and set $\tilde\theta=(AH,\beta+A\mu)$. A change of variables then gives $P_{\theta,G} = P_{\tilde\theta,\tilde G}$.
We first address this affine indeterminacy, and then formulate identifiability in terms of the remaining ambiguity arising from the interaction between the unknown latent distribution and the model parameters, which is the main focus of our analysis.


First, we adopt a common normalization on the latent variables: 
\begin{equation}\label{eq_unv_normalize}
    \mathbb E_G (Z)=0,\qquad \operatorname{Cov}_G(Z)=I_d.
\end{equation}
 When $\mathbb E_G Z=m$ and $\operatorname{Cov}_G(Z)=\Sigma\succ0$, we can simply take $\mu=m$ and $H=\Sigma^{1/2}$ to yield a standardized latent distribution without changing the marginal distribution. This normalization can also be replaced by more interpretable procedures; see Remark~\ref{remark_inter_id} below.

Normalization \eqref{eq_unv_normalize} still leaves a rotational indeterminacy, as for any orthogonal matrix $H\in\RR^{d\times d}$, the preceding reparameterization preserves both normalization \eqref{eq_unv_normalize} and the marginal distribution. To account for this, we measure the distance metric between any $\theta = (A,\beta)$, $\tilde\theta = (\tilde A,\tilde\beta)$ as
\begin{equation}\label{eq_def_Dp}
    D_p(\theta,\tilde\theta) \; = \; \|\beta - \tilde\beta\|_{\infty} + \|AA^{\T} - \tilde A\tilde A^{\T}\|_{\infty} \; = \; \max_{j\in [p]}|\beta_j-\tilde\beta_j| +\max_{i,j\in [p]} |\alpha_i^{\mathsf T}\alpha_j -\tilde\alpha_i^{\mathsf T}\tilde\alpha_j|.
\end{equation}
This is a rotation-invariant distance. In particular, $D_p(\theta,\tilde\theta)=0$ if and only if $(A,\beta) = (\tilde A H,\tilde\beta)$ for some orthogonal matrix $H$, and $D_p(\theta,\tilde\theta)$ controls the maximum rowwise discrepancy between
 $A$ and $\tilde A$ after some orthogonal Procrustes alignment~\citep{ten1977orthogonal}; see Lemma~\ref{lem_design_alignment_sub}.

Our notion of model identifiability is then formulated as follows.
\begin{definition}\label{def_id}
    Consider model~\eqref{eq_intro_model}. For each $p$, let $\cM_p$ be a nonempty set of pairs $(\theta,G)$ satisfying~\eqref{eq_unv_normalize}. The model is said to be \textbf{exactly identifiable} over $\cM_p$ if, for any $(\theta,G),(\tilde\theta,\tilde G)\in\cM_p$,
    \begin{equation*}P_{\theta,G} \, = \, P_{\tilde\theta,\tilde G} \quad\Longrightarrow\quad D_p(\theta,\tilde\theta) = 0 \end{equation*} 
    The model is said to be \textbf{asymptotically identifiable} over $(\cM_p)_{p\ge1}$ if, for any two sequences satisfying $(\theta_p,G_p),(\tilde\theta_p,\tilde G_p)\in\cM_p$, \begin{equation*}P_{\theta_p,G_p} \, = \, P_{\tilde\theta_p,\tilde G_p},\;\forall p \quad\Longrightarrow\quad D_p(\theta_p,\tilde\theta_p)\to 0 \;\text{ as }\; p\to\infty \end{equation*}
\end{definition}
Exact identifiability is a standard assumption underlying consistency and inference results, whereas asymptotic identifiability is a weaker requirement that applies as the number of observed variables $p$ grows. When $P_{\theta,G}=P_{\tilde\theta,\tilde G}$, the two parameterizations induce the same distribution of the observed data and hence cannot be distinguished by any data-driven procedure, with the discrepancy $D_p(\theta,\tilde\theta)$ quantifying the resulting population-level ambiguity. Classical theory for consistent point estimation and inference at fixed $p$ therefore typically requires this discrepancy to be zero. Nevertheless, under suitable scaling of $p$ and the sample size, asymptotic identifiability may still permit consistent estimation and valid inference, particularly when the ambiguity vanishes rapidly.

\begin{remark}\label{remark_inter_id} \em
In Definition~\ref{def_id}, the affine indeterminacy is resolved by the normalization~\eqref{eq_unv_normalize} together with the rotationally invariant discrepancy $D_p(\cdot,\cdot)$. Besides these mathematically convenient conventions, rotation methods are often used in practice to select a simpler and more interpretable representation among equivalent ones~\citep{browne2001overview}. Two common settings are as follows:
\begin{enumerate}[label=(\arabic*),leftmargin=0.7cm]
    \item \textbf{Orthogonal rotation.} Suppose we retain the normalization~\eqref{eq_unv_normalize}, and the remaining indeterminacy regards orthogonal transformations, as recovering $A$ from $AA^{\T}$ requires choosing an orientation within $\{AH:H^{\T}H=I_d\}$. Rotation methods such as varimax~\citep{kaiser1958varimax,Rohe2020VintageFA} can then be used to select some $H^*$ to obtain $AH^*$ with a simpler structure.
    \item \textbf{Oblique rotation.} To allow correlated latent variables, condition ${\rm Cov}_G(Z) = I_d$ can be relaxed to ${\rm diag}\{{\rm Cov}_G(Z)\} = I_d$. Starting from a representation satisfying~\eqref{eq_unv_normalize}, any invertible $H$ satisfying ${\rm diag}(HH^{\T})=I_d$ yields an equivalent representation
    \[
        \tilde Z=HZ,\qquad\; \tilde A=AH^{-1},\qquad
        \operatorname{Cov}(\tilde Z)=HH^{\T},
    \]with ${\rm diag}\{{\rm Cov}_G(\tilde Z)\} = I_d$.
    The transformed latent coordinates have unit variances and may be correlated. Oblique rotation methods~\citep{Liu2023RotationTS,cui2026beyond} can select some $H^*$ to yield a simpler structure in $A(H^*)^{-1}$.
\end{enumerate}
\end{remark}


\section{Main Results}\label{sec_main_res}

\subsection{Impossibility of Exact Identifiability}\label{sec_non_id}

When $G$ is unspecified, it ranges over an infinite-dimensional space, and it is commonly believed that exact identifiability cannot hold in this case~\citep{lindsay1991semiparametric,haberman2005identifiability,san2011bayesian}.
The following result confirms this belief and sharpens it by showing that exact identifiability fails almost everywhere in the parameter space, even when constraining finitely many moments of the latent distribution. 
\begin{theorem}[Non-identifiability]\label{thm_nonidentifiability}  Assume that $f(x)$ is strictly increasing and analytic with limits $f(+\infty) = 1$ and $f(-\infty) = 0$, and that $G$ has a positive $C^\infty$ density on $\RR^d$. Then for every finite $p\ge1$, for every $\theta$ outside a Lebesgue zero set, there exists a neighbourhood $\cU$ of $\theta$ such that each $\tilde\theta\in\cU$ admits some $\tilde G$ satisfying \eqref{eq_unv_normalize} and
\[
P_{\theta,G}\; =\; P_{\tilde\theta,\tilde G}.
\]
Consequently, the model is not exactly identifiable in the sense of Definition~\ref{def_id}. Moreover, for any fixed integer $s\ge1$ with $\EE_G(\|Z\|^{2s})<\infty$, $\tilde G$ can be chosen so that
\[
\EE_{\tilde G}[Z^a]\; =\; \EE_G[Z^a]\, , \text{ for all }\; |a|\le 2s,
\]
where $a$ is a multi-index and $Z^a=\prod_{k=1}^d Z_k^{a_k}$. 
\end{theorem}
This result establishes that, when $G$ is unspecified, exact identifiability fails for almost every parameter $\theta$, in the sense that the exception is a Lebesgue zero set. This exceptional set is used to exclude degenerate parameter configurations, such as $\alpha_j=0$ or identical $(\alpha_j,\beta_j)$ across $j\in[p]$.
Moreover, 
this identifiability issue cannot be addressed by restricting to some local neighbourhood~\citep{rothenberg1971identification}, as every sufficiently small parameter perturbation can be offset by a corresponding change in the latent distribution while preserving the marginal distribution. 
 Finally, in the proof, $G$ and $\tilde G$ are constructed to share the same moments up to any prespecified order, beyond the normalization in \eqref{eq_unv_normalize}, provided these moments exist.  Consequently, no fixed finite collection of moment restrictions on the latent distribution can restore exact identifiability.

\begin{example}\label{example_nonid}
We illustrate the construction in Theorem~\ref{thm_nonidentifiability} using a one-dimensional logistic model with $f(t)=(1+e^{-t})^{-1}$ and latent distribution $G=N(0,1)$, whose density is denoted by $\phi(\cdot)$. Fix an admissible parameter $\theta=(A,\beta)$ as in that theorem. For a fixed $j\in[p]$, consider the perturbed alternative $\theta_\varepsilon=(A,\beta+\varepsilon e_j)$, where $e_j\in\RR^p$ is the $j$th standard basis vector. We construct a perturbed distribution $G_\varepsilon$ that offsets this parameter change, so that $P_{\theta,G}=P_{\theta_\varepsilon,G_\varepsilon}$, while matching the moments of $G$ through order $2s$ and differing from $G$ only within $[-1,1]$.

Write $q_S(\theta,z)=\prod_{k\in S}f(\beta_k+\alpha_kz)$ for $S\subseteq [p]$ and define
\begin{equation*}
 v_\varepsilon(z) \; = \; \Big[
 1,z^2,\ldots,z^{2s},
 \big\{{q_S(\theta_\varepsilon,z)+q_S(\theta_\varepsilon,-z)}
 \big\}_{\varnothing\ne S\subseteq[p]}
 \Big]^{\T} \in \RR^{s + 2^p}.
\end{equation*}
This feature vector encodes the normalization, the latent moments through order $2s$, and the response probabilities of all subsets of $[p]$. Assume it is linearly independent on $(-1,1)$. To localize the perturbation to the distribution, we use the smooth bump
\begin{equation*}
 w(z) \; = \; \begin{cases}\, \exp\{-1/(1-z^2)\},& |z|<1,\\ \, 0,& |z|\ge1. \end{cases}
\end{equation*}
Set a weighted Gram matrix  
$M_\varepsilon=\int_{\RR}w(z)v_\varepsilon(z)v_\varepsilon(z)^{\T}\phi(z)\,dz$ and the vector of feature discrepancies $b_\varepsilon=\int_{\RR}\{v_0(z)-v_\varepsilon(z)\}\phi(z)\,dz$. For sufficiently small $|\varepsilon|$, $M_\varepsilon$ is invertible. Then define the perturbed distribution by $G_\varepsilon$ with density $g_\varepsilon(z)=\phi(z)\{1+w(z)v_\varepsilon(z)^{\T}M_\varepsilon^{-1}b_\varepsilon\}$. It is easy to check that $g_\varepsilon$ is a positive and symmetric density for sufficiently small $|\varepsilon|$ with $g_\varepsilon(z)=\phi(z)$ for $|z|\ge1$. 
By construction, one has
\begin{equation}\label{eq_example_constraint_matching}
\begin{aligned} \int_{\RR}v_\varepsilon(z)g_\varepsilon(z)\,dz
 =&\,\int_{\RR}v_\varepsilon(z)\phi(z)\,dz + \Big\{\int_{\RR}w(z)v_\varepsilon(z)v_\varepsilon(z)^{\T}\phi(z)\,dz\Big\}M_{\varepsilon}^{-1}b_{\varepsilon}\\=&\,\int_{\RR}v_\varepsilon(z)\phi(z)\,dz+b_\varepsilon
 =\int_{\RR}v_0(z)\phi(z)\,dz.\end{aligned}
\end{equation}
By symmetry of $g_{\varepsilon}$ and $\phi$, this identity implies that $\EE_{G}[Z^a] = \EE_{G_\varepsilon}[Z^a]$ for all $a\le 2s$ and $\EE_G[q_S(\theta,Z)] = \EE_{G_\varepsilon}[q_S(\theta_\varepsilon,Z)]$ for all $S$. Therefore, $P_{\theta,G} = P_{\theta_\varepsilon,G_\varepsilon}$ with $D_p(\theta_\varepsilon,\theta)=|\varepsilon|$. Analogous constructions can offset any other sufficiently small perturbation through changes in the latent distribution while preserving its moments to any finite order.


\end{example}



Theorem~\ref{thm_nonidentifiability} establishes non-identifiability when the latent dimension $d$ is fixed. We next show that the latent dimension $d$ itself can also be ambiguous.
Borrowing terminology from factor analysis~\citep{onatski2012asymptotics}, we call a latent variable a weak factor if the corresponding column of $A$ is small. The proposition shows that appending a weak factor to $A$ can be offset by a perturbation of the latent distribution, leaving the marginal distribution unchanged. 
This ambiguity further motivates our treatment of the latent dimension $d$ as specified throughout.

\begin{proposition}[Ambiguity from a weak factor]\label{prop_weak_factor}
Under the setup of Theorem~\ref{thm_nonidentifiability}, suppose $p>d$, ${\rm rank}(A)=d$, and that $\theta=(A,\beta)$ is outside the Lebesgue zero set specified by Theorem~\ref{thm_nonidentifiability}. Let $H$ be any distribution on $\RR$ with mean zero and unit variance, and let $G\otimes H$ denote the joint distribution of independent $Z\sim G$ and $U\sim H$. Fix any $v\in\RR^p\setminus{\rm col}(A)$ and an integer $s\ge1$ with $\EE_G\|Z\|^{2s}<\infty$. Then, for sufficiently small $|\varepsilon|$, there exists a $d$-dimensional latent distribution $\tilde G_\varepsilon$ satisfying \eqref{eq_unv_normalize} such that
\begin{equation*}
 P_{([A,\varepsilon v],\beta),\,G\otimes H} \; 
 = \; P_{(A,\beta),\,\tilde G_\varepsilon}
 \quad\text{and}\quad
 \EE_G[Z^a] \;  =\;  \EE_{\tilde G_\varepsilon}[Z^a]\ \text{ for all } |a|\le2s.
\end{equation*}
The additional variable $U$ is the weak factor, with the corresponding column $\varepsilon v$ of order $\varepsilon$.
\end{proposition}
This result shows that a $(d+1)$-dimensional model with one weak factor can yield the same marginal distribution for the observed variables as a $d$-dimensional model with a suitably perturbed latent distribution. This reflects the difficulty of detecting weak factors in factor analysis~\citep{onatski2012asymptotics} from the perspective of model identifiability. The equivalence also has a more positive interpretation that a flexible latent distribution can absorb the effect of omitted weak signal factors without changing the model parameters associated with strong factors.

\subsection{Possibility of Asymptotic Identifiability}

While Theorem~\ref{thm_nonidentifiability} settles \hyperlink{Q1}{(Q1)} in the negative, reliable estimation may still be possible if the model is asymptotically identifiable, which would answer \hyperlink{Q2}{(Q2)} in the positive.
In this subsection, we establish asymptotic identifiability of the model by deriving a quantitative bound on the discrepancy $D_p(\theta,\tilde\theta)$ that holds uniformly over all $\{(\theta, G), (\tilde\theta,\tilde G)\}$ that yield the same marginal distribution. The subsequent Section~\ref{sec_impl} draws two implications of this bound to answer \hyperlink{Q2}{(Q2)}. 

We first introduce two technical assumptions.
\begin{assumption}\label{ass_link}\em
The link function $f:\mathbb R\to(0,1)$ satisfies the following conditions:
\begin{enumerate}[label=$(\roman*)$]
    \item\label{ass_link_1}  $f$ is analytic with $f'(x)>0$ for every $x\in\mathbb R$, $\lim_{x\to+\infty}f(x)=1$, and $\lim_{x\to-\infty}f(x)=0$. Moreover, $f'(0)\ge m$ for some constant $m>0$.
    \item\label{ass_link_2} There exists a constant $K\ge 1$ such that for any integer $k\ge2$, $\sup_{x\in\mathbb R}\frac{|f^{(k)}(x)|}{k!f'(x)}
    \le K^{k-1}$.
\end{enumerate}
\end{assumption}
Assumption~\ref{ass_link} imposes regularity conditions on the link function and is satisfied by, for example, the logistic link $f(x)=(1+e^{-x})^{-1}$ and the Cauchy link $f(x)=\tfrac12+\pi^{-1}\arctan(x)$.
The first condition ensures a strictly increasing analytic response function with the usual boundary limits, with the lower bound on
$f'(0)$ preventing it from becoming arbitrarily flat at $0$.
Condition~\ref{ass_link_2} is the key technical requirement for the superpolynomial rate. It implies that $f$ admits a Taylor expansion to arbitrary order around any point with coefficients bounded relative to $f'$. This uniform control is crucial for our technique to apply; see Remark~\ref{remark_technique} for more discussion. 


\begin{assumption}\label{ass_model}\em
For $(\theta,G)\in\cM_p$ with $\cM_p$ given in Definition~\ref{def_id}, we assume
\begin{enumerate}[label=$(\roman*)$]
    \item\label{ass_model_1} for some constant $M > 0$, $\max_{j\in [p]}\|\alpha_j\|\le M$; $\max_{j\in [p]}|\beta_j|\le M$;
    \item\label{ass_model_2} $\lambda_{\min}(p^{-1}{A}^\T A)\ge \kappa >0 $ for some constant $\kappa$;
    \item\label{ass_model_3} there exist some constants $\nu> 0$ and $L$ such that $\mathbb E_G [\exp(\nu\|Z\|^2)]\le L$.
\end{enumerate}
\end{assumption}
Assumption~\ref{ass_model} imposes
boundedness and nondegeneracy conditions standard in the latent-variable literature. 
Assumption~\ref{ass_model}\ref{ass_model_3} imposes a uniform sub-Gaussian tail bound on the latent distribution. 
It accommodates Gaussian latent distributions, widely used in classical item response models~\citep{bock1981marginal,reckase2009multidimensional}, finitely supported discrete latent distributions, as arise in latent-class models \citep{allman2009identifiability,xu2017}, and bounded-support latent distributions~\citep{he2024extended,cui2025identifiability}.
The condition places no restriction on the shape of the latent distribution, allowing asymmetry and multimodality, and does not require a density or any parametric form. 

To quantify the identification ambiguity over the admissible set $\cM_p$, we define
\begin{equation*}
    \omega_p(\mathcal M_p) \;  = \; 
    \sup\big\{D_p(\theta,\tilde\theta):  (\theta,G),(\tilde\theta,\tilde G)\in\mathcal M_p, \ P_{\theta,G}=P_{\tilde\theta,\tilde G}\ \big\}.
\end{equation*}
This gives the supremum of the $D_p$-diameters of sets of all admissible parameters that yield some common marginal distribution. In Example~\ref{example_nonid}, for instance, it upper bounds $|\varepsilon|$ for every perturbation that preserves the marginal distribution. 
We have the following superpolynomial bound for $\omega_p(\cM_p)$.

\begin{theorem}[Asymptotic identifiability]\label{thm_asymptotic}
    Under Assumptions~\ref{ass_link} and~\ref{ass_model}, the model \eqref{eq_intro_model} is asymptotically identifiable. In particular, there exist $c,C,p_0$ such that for every $p\ge p_0$, it holds that \begin{equation*}
    \omega_p(\mathcal M_p) \; \le\; C\exp\{-c(\log p)^2\}.
\end{equation*}
\end{theorem}

Theorem~\ref{thm_asymptotic} shows that the identification ambiguity decays as fast as $\exp\{-c(\log p)^2\}=p^{-c\log p}$, and hence faster than any polynomial rate in $p$. The bound holds uniformly over all pairs in $\cM_p$, with $D_p(\theta,\tilde\theta)$ controlling the entrywise discrepancy between $\theta$ and $\tilde\theta$. This essentially quantifies the size of the neighbourhood $\mathcal U$ in Theorem~\ref{thm_nonidentifiability}, bounding the discrepancy of admissible alternatives by $C\exp\{-c(\log p)^2\}$ in metric $D_p(\cdot,\cdot)$. It is therefore suggested that non-identifiability may pose less of an impediment to inference in latent variable models than is commonly expected.
The result is stronger than existing asymptotic identifiability results~\citep{douglas2001asymptotic,he2024extended}, which only establish a qualitative result under unidimensional models. Establishing the superpolynomial rate in Theorem~\ref{thm_asymptotic} requires a more delicate analysis, which we illustrate as follows.

\begin{remark}\label{remark_technique}\em
The primary challenge in establishing Theorem~\ref{thm_asymptotic} is to approximate the parameters via observation $Y$ to a sufficiently high order. We start with a conditioning argument. Select a small set $S\subseteq [p]$ containing the indices of the parameters to approximate together with a few reference indices, and we condition on the remaining responses $Y = X_{S^c}$. Then $Z$ can be approximated by its conditional mean $\EE(Z\mid Y)$, with residual $V=Z-\EE(Z\mid Y)$. Then for $i,j\in S$, we can write
\begin{equation*}
    \beta_i \, = \, \EE(T_i)+ \EE(\delta_i)\quad\text{ and }\quad \alpha_i^{\T}\alpha_j \, = \, \operatorname{Cov}(T_i+\delta_i,T_j+\delta_j)
+\EE\{(\alpha_i^{\T}V)(\alpha_j^{\T}V)\},
\end{equation*}
where $T_i=f^{-1}\{\EE(X_i\mid Y)\}$ is determined by the marginal distribution, 
whereas $\delta_i=\EE(\alpha_i^{\T}Z+\beta_i\mid Y)-T_i$ depends on unknown model parameters.  Here, $\delta_i$ is a correction term that we separate out to simplify the derivation. 

The key difficulty is to approximate $\delta_i$ with sufficiently small error. We work with the conditional latent moments $M_{k,i}=\EE\{(\alpha_i^{\T}V)^k\mid Y\}$, which satisfy $\|M_{k,i}\|\le e^{Ck^2}p^{-k/2}$ for $k\ge2$ (Lemma~\ref{lem_conditional_fluctuations}), so that expansions in these moments become more accurate as the order increases. For $\delta_i$, for instance, to fourth order, it can be approximated by
\begin{align*}
\delta_i\approx -\frac{f''(T_i)}{2f'(T_i)}M_{2,i}
-\frac{f'''(T_i)}{6f'(T_i)}M_{3,i}
-\frac{f^{(4)}(T_i)}{24f'(T_i)}M_{4,i}+\left\{\frac{f''(T_i)f'''(T_i)}{4f'(T_i)^2}
-\frac{f''(T_i)^3}{8f'(T_i)^3}\right\}M_{2,i}^{\,2}
+\cdots
\end{align*}
Here, $M_{4,i}$ and $M_{2,i}^2$ are both of order $p^{-2}$ but must be tracked separately. Higher orders generate more moment products of this type, with growing complexity in their coefficients, and it is thus difficult to ensure the precision of the approximation. Moreover, the latent moments are not observed and thus must themselves be approximated from observed quantities, which introduces further error terms that have to cancel successively to reach the desired order.

To handle these calculations, we introduce a symbolic calculus that organizes moment products as the expansion order $N$ increases. This separates algebraic bookkeeping from the probabilistic estimates needed for the generally dependent conditional moments. Then for pairs $(\theta,G)$ and $(\tilde\theta,\tilde G)$ that yield the same marginal distribution, we can construct the expansions with $\theta$ first, and then apply the same expansions to $\tilde\theta$. After careful control of the higher-order terms, we obtain
\[
D_p(\theta,\tilde\theta) \; \le \; CD_p(\theta,\tilde\theta)/p+e^{CN^2}p^{-(N+1)/2}.
\]
The first term controls the parameter-dependent differences between
the expansions, with the factor $p^{-1}$ reflecting the order at which the correction terms begin. The second term bounds the order-$N$ approximation error uniformly over $\theta$ and $\tilde\theta$, which aligns with the bound for latent moments $M_{k,i}$, with $e^{CN^2}$ accounting for the growth of the expansion coefficients.
Absorbing the first term and choosing $N=\lfloor\tau\log p\rfloor$
for sufficiently small $\tau>0$ gives the superpolynomial bound.
\end{remark}

\begin{remark}\em
The technique discussed in Remark~\ref{remark_technique} could potentially be extended to more general models, such as response types beyond binary data. Such extensions can be made via similar calculations for the expansion and control for latent moments. We focus on binary responses here to present the argument in a representative nonlinear model and, to our knowledge, to give the first quantitative bound on identification ambiguity in this setting. We leave extensions under more complex models to future work.
\end{remark}

We further specialize our results to three concrete examples of latent variable models. These specializations follow directly from our main results, but the results may be of independent interest in each setting.
\begin{example}[Bifactor Models] \label{example_bifactor}
The bifactor model is an important latent trait model in psychological and educational assessment~\citep{reise2012rediscovery}. It specifies one general variable $Z_0$ and $d$ group variables $Z_1,\ldots,Z_d$, giving $d+1$ latent variables in total. Each observed response $X_j$ is allowed to depend only on the general variable $Z_0$ and at most one group variable $Z_{g(j)}$, where $g(j)\in\{1,\ldots,d\}$. Its conditional response probability therefore takes the form
\begin{equation*}
\operatorname{pr}(X_j=1\mid Z) \; = \; f\big\{\ \beta_j+\lambda_{j0}Z_0+\lambda_{j1}Z_{g(j)}\ \big\}.
\end{equation*}
The vector $Z = (Z_0,Z_1,\dots,Z_d)^\T$ has mean zero and covariance $\Sigma={\rm diag}(1,R)$, where $R$ has unit diagonal. This covariance structure requires $Z_0$ to be uncorrelated with the group variables, a common restriction used to facilitate interpretation and identification. In particular, \citet{fang2021identifiability} established identifiability under this restriction with Gaussian latent variables.

Our theory establishes asymptotic identifiability for a broader class of latent distributions beyond the Gaussian family, while retaining the covariance structure. Let \(\lambda_j\in\mathbb R^{d+1}\) have first entry \(\lambda_{j0}\), entry \(g(j)+1\) equal to \(\lambda_{j1}\), and all remaining entries zero. Write \(\Lambda=(\lambda_1,\ldots,\lambda_p)^\T\). Under suitable regularity conditions (see details in Section~\ref{app_sec_bif} in the Supplementary Material), Theorem~\ref{thm_asymptotic} implies that, for $(\Lambda,\beta,G)$ and $(\tilde\Lambda,\tilde\beta,\tilde G)$ that yield the same marginal distribution,
\begin{equation*}\|\beta-\tilde\beta\|_\infty+\|\tilde\Lambda-\Lambda\|_{\infty}+\|\tilde\Sigma-\Sigma \|_{\infty} \; \le \; C\exp\{-c(\log p)^2\},\end{equation*}
for some constants $C,c>0$, after proper sign change of the latent variables, where $\Sigma$ and $\tilde\Sigma$ are the covariance matrices of $G$ and $\tilde G$, respectively.
\end{example}

\begin{example}[Independent Component Analysis] 
Independent component analysis (ICA) aims to recover independent latent variables, also referred to as latent sources, from observed mixtures~\citep{hyvarinen1999nonlinear,khemakhem2020variational}. Model~\eqref{eq_intro_model} can be viewed as a binary ICA model when $G$ is restricted to have independent components $Z_1,\ldots,Z_d$, with $A$ being the mixing matrix specifying how the latent sources are combined to determine the response probabilities. Theorem~\ref{thm_asymptotic} implies that $\min_{Q^{\T}Q=I_d}\max_{j\in [p]}\|\tilde\alpha_j-Q\alpha_j\|\le C\exp\{-c(\log p)^2\}$. Thus, as the number $p$ of observed mixtures grows, the mixing matrix is asymptotically identifiable up to an orthogonal transformation, an indeterminacy that is commonly resolved by independence and non-Gaussianity conditions on the sources~\citep{hyvarinen2000independent}. This result provides a quantitative perspective on the identifiability theory of ICA.
\end{example}

\begin{example}[Generative Models with Binary Decoders] 
Model~\eqref{eq_intro_model} can be viewed as the generative model underlying a variational autoencoder for binary data~\citep{kingma2013auto}, with prior $G$ and decoder $m_{\theta}(z)=\{f(\beta_j+\alpha_j^{\T}z)\}_{j\in[p]}$. Learning the prior, for example through normalizing flows~\citep{papamakarios2019flows}, can introduce additional indeterminacies. For instance, \citet{xi2023indeterminacy} characterized such indeterminacies, but their analysis is limited to continuous decoders. Our result shows that, under suitable regularity conditions and with $G$ left unspecified, any two observationally equivalent representations satisfy $\|\beta-\tilde\beta\|_\infty+\max_{j\in [p]}\|\tilde\alpha_j-Q\alpha_j\|\le C\exp\{-c(\log p)^2\}$ for some orthogonal matrix $Q$. Consequently, for each fixed $R>0$, the decoder discrepancy satisfies
\begin{equation*}
\sup_{\|z\|\le R}
\|m_{\tilde\theta}(Qz)-m_\theta(z)\|_\infty
\; \le \; C_R\exp\{-c(\log p)^2\},\end{equation*}where $C_R$ is some constant depending on the radius $R$.
Thus the decoder is asymptotically identifiable up to a rotation of the latent space, even when the prior is learned flexibly.
\end{example}

\subsection{Implications}\label{sec_impl}
This subsection presents two consequences of Theorem~\ref{thm_asymptotic}, concerning estimability of parameter $\theta$ and latent distribution $G$ to address \hyperlink{Q2}{(Q2)}. The first concerns parameter estimation.


\begin{corollary}[Estimation of model parameter]\label{cor_estimation}
Under Assumptions~\ref{ass_link} and~\ref{ass_model}, suppose $(\theta_{0},G_{0})\in\mathcal M_p$ and $(\theta^{\dagger}, G^{\dagger})\in\cM_p$ yield the same marginal distribution $P_{\theta_0,G_0}$ under model \eqref{eq_intro_model}. For any estimator $\hat\theta_n$ based on $n$ independent responses, we have
\begin{equation}\label{eq_estimation_bound}
    D_p(\hat\theta_{n},\theta_{0})
    \; = \; D_p(\hat\theta_{n},\theta^{\dagger}) + O_{\mathbb P}\big(\exp\{-c(\log p)^2\}\big).
\end{equation}
Thus when $n,p\to\infty$ and $D_p(\hat\theta_{n},\theta^{\dagger})\to0$, the estimator is consistent up to some orthogonal transformation. 
Moreover, suppose $\psi(\cdot)$ is a scalar functional satisfying $|\psi(\theta)-\psi(\tilde\theta)| \le L_\psi D_p(\theta,\tilde\theta)$ for a constant $L_\psi$ independent of $p$. If $\hat s_{n,p}^{\, -1}\{\psi(\hat\theta_n)-\psi(\theta^{\dagger})\}$ has a limiting distribution as $n,p\to\infty$, for some sequence $\hat s_{n,p}$ satisfying $\hat s_{n,p}^{\, -1}\exp\{-c(\log p)^2\}\rightarrow0$,
then $\hat s_{n,p}^{\, -1}\{\psi(\hat\theta_n)-\psi(\theta_0)\}$ has the same limiting distribution.
\end{corollary}
Corollary~\ref{cor_estimation} separates the estimation error into two parts. The identification ambiguity, captured by the second term in \eqref{eq_estimation_bound}, is controlled by Theorem~\ref{thm_asymptotic} regardless of the estimator.
The remaining $D_p(\hat\theta_n,\theta^\dagger)$ measures the sampling error relative to a target $\theta^{\dagger}$ satisfying $P_{\theta^{\dagger},G^{\dagger}} =P_{\theta_0,G_0} $ for some $G^{\dagger}$, which depends on the specific estimation procedure.
The corollary also allows valid inference for $\psi(\theta_0)$, provided that $\hat s_{n,p}^{-1}{\psi(\hat\theta_n)-\psi(\theta^\dagger)}$ converges to a normal distribution and the identification ambiguity is negligible relative to the standard error. For root-$n$ inference, for instance, the
additional requirement is $ n^{1/2}\exp\{-c(\log p)^2\}\to0$, which holds, for example, whenever $\log n=o\{(\log p)^2\}$.
The result thus provides a basis for consistent estimation and valid inference for procedures converging to population solutions that reproduce the true marginal distribution. Such population solutions arise naturally in marginal likelihood estimation. Identification ambiguity contributes only a superpolynomially small error and therefore need not obstruct consistent estimation.

There have been studies developing asymptotic theory for estimation of model \eqref{eq_intro_model}. 
In the semiparametric literature, for instance, \citet{lindsay1991semiparametric} linked conditional likelihood and mixture likelihood estimation in the Rasch model and related exponential response models.
For generalized linear latent variable models, \citet{ma2010explicit} developed semiparametric estimating equations and established root-$n$ consistency and asymptotic normality. Both studies take model identifiability as given.
In a growing-$p$ regime, \citet{cui2026marginal} recently established consistency and asymptotic normality of marginal maximum likelihood estimators under a working Gaussian latent distribution, requiring only that the true latent distribution be sub-Gaussian. 
In their analysis, the working density contributes a leading first-order term to the likelihood expansion, and the resulting asymptotic theory targets a canonical representation that need not coincide with the true parameter.

\medskip

We next consider the estimation and inference of the latent distribution $G$. We measure the discrepancy between two latent distributions by the $1$-Wasserstein distance
\begin{equation*}
 W_1(G,\tilde G) \; =\inf_{\pi\in\Pi(G,\tilde G)} \int_{\RR^d\times \RR^d}\|z-\tilde z\|\,d\pi(z,\tilde z),
\end{equation*}
where $\Pi(G,\tilde G)$ is the set of couplings of $G$ and $\tilde G$. For an orthogonal matrix $Q$, write $Q_{\#}G$ for the distribution of $QZ$ when $Z\sim G$. We have the following result.

\begin{corollary}[Estimation of the latent distribution]\label{prop_prediction}
Under Assumptions~\ref{ass_link} and~\ref{ass_model}, suppose $(\theta_0,G_0),(\theta^{\dagger},G^{\dagger})\in\mathcal M_p$ and $P_{\theta_0,G_0}=P_{\theta^{\dagger},G^{\dagger}}$.
There exists an orthogonal matrix $Q$ such that, for any estimator $\hat G_n$ based on $n$ independent responses and having a finite first moment,
\begin{equation}\label{eq_distribution_estimation_bound}
 W_1(Q_\#\hat G_n,G_0)
 \; = \; W_1(\hat G_n,G^{\dagger})+O_{\mathbb P}(p^{-1/2}).
\end{equation}
Consequently, when $n,p\to\infty$ and $W_1(\hat G_n,G^{\dagger})\to0$, $Q_{\#}\hat G_n$ is consistent for $G_0$ in $W_1$. Moreover, for every $L_h$-Lipschitz function $h:\RR^d\to\RR$, satisfying $|h(z)-h(z')|\le L_h\|z-z'\|$, we have, almost surely for all sufficiently large $p$,
\begin{equation}\label{eq_prediction_bound}
 \big|\ \EE_{\hat G_n}[h(QZ)]\ -\ \EE_{G_0}[h(Z)]\ \big| \; 
 \le \; L_h\big\{W_1(\hat G_n,G^{\dagger})+Cp^{-1/2}\big\}.
\end{equation}
\end{corollary}

Corollary~\ref{prop_prediction} shows in \eqref{eq_distribution_estimation_bound} that, as $p\to\infty$, $Q_{\#}\hat G_n$ can be consistent for $G_0$ provided that $\hat G_n$ converges in $W_1$ to some $G^{\dagger}$ that, together with some $\theta^{\dagger}$, induces the same marginal distribution as $(\theta_0,G_0)$.
The identification ambiguity for the latent distribution is of order $p^{-1/2}$. For comparison, the minimax $W_1$ rate for estimating a distribution on $\RR^d$ from $n$ observations is $n^{-1/d}$ for $d\ge3$ without involving latent variables~\citep{singh2018minimax}, so the ambiguity is negligible when $p\gg n^{2/d}$.

Equation~\eqref{eq_prediction_bound} is useful from a practical view. It shows that the failure of exact identifiability need not prevent accurate estimation of population features of the latent distribution, which are often the quantities of scientific interest.
For instance, in a latent-trait measurement model, $h(z)$ may be the success probability on a new testing item for an individual with latent trait $z$, so that $\EE_{G_0}[h(Z)]$ is the population success rate on that item. In a Bernoulli generative model for binary images, $h(z)$ may be the mean activation probability over a specified image region under a given decoder, so that $\EE_{G_0}[h(Z)]$ is the expected fraction of active pixels in that region. These guarantees provide a theoretical basis for downstream population-level inference. 

\section{Simulation Study}\label{sec_simulation}
We illustrate Theorems~\ref{thm_nonidentifiability} and~\ref{thm_asymptotic} via a numerical study. We use a logistic function for $f(\cdot)$. Consider $d=1$ and $d=3$, and
for each $d$, we draw ten sets of parameters $\theta_0=(A,\beta)$ with independent rows. For $d=1$, $\alpha_j\sim U(0.5,1.5)$. For $d=3$, $\alpha_j=r_ju_j$, where $r_j\sim U(0.5,1.5)$ and $u_j$ is uniform on the unit sphere of $\RR^3$. In both cases, $\beta_j\sim U(-1,1)$. A model with $p$ responses uses the first $p$ rows. We consider two choices of latent distributions $G_0$: $N(0,I_d)$, and the mixture $0.7\,N(0,0.25I_d)+0.3\,N(2d^{-1/2}\mathsf{1}_d,0.25I_d)$, with $\mathsf{1}_d$ denoting the $d$-dimensional all-$1$ vector. 


For a candidate parameter $\tilde\theta$, we measure how well $P_{\theta_0,G_0}$ can be reproduced with $\tilde\theta$ by
\begin{equation}\label{eq_sim_delta}
    \Delta_p(\tilde\theta)=\min_{\tilde G\in\mathcal G_{{\rm grid},d}}\frac{\|p_{\tilde\theta,\tilde G}-p_{\theta_0,G_0}\|}{\|p_{\theta_0,G_0}\|},
\end{equation}
where $p_{\theta,G}=\{p_{\theta,G}(x)\}_{x\in\{0,1\}^p}$ is the vector of response probabilities defined in~\eqref{eq_marginal_law}. The class $\mathcal G_{{\rm grid},d}$ consists of distributions supported on a fixed grid that satisfy~\eqref{eq_unv_normalize} and the tail bound $\EE_{\tilde G}\exp(\|Z\|^2/4)\le 2^{(d+1)/2}$ as required by Assumption~\ref{ass_model}\ref{ass_model_3}.
This grid contains $121$ equally spaced points on $[-12,12]$ for $d=1$ and $17^3$ points on $[-6,6]^3$ for $d=3$. We evaluate $P_{\theta_0,G_0}$ using a finer discretization of $G_0$, with $241$ and $41^3$ points on the respective domains, followed by exact standardization. For $d=1$ and $p=7,8,9$, we discretize $G_0$ and $\tilde G$ on shared standardized grids of $257$, $513$ and $1025$ points, respectively, and increase the arithmetic precision to $2500$, $4000$ and $6500$ bits to resolve the smaller perturbations. Thus, $\Delta_p(\tilde\theta)=0$ if and only if there exists an admissible $\tilde G \in \mathcal G_{{\rm grid},d}$ such that $P_{\tilde\theta,\tilde G}=P_{\theta_0,G_0}$. Numerical equivalence is assessed via this discretized form, with details described in Appendix~\ref{app_sim}. 

\begin{figure}[t]
    \centering
    \includegraphics[width=0.9\textwidth]{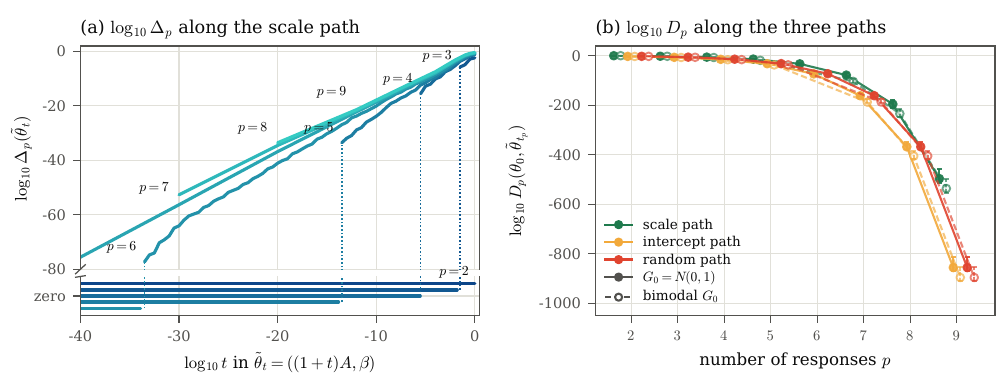}
    \caption{Exact equivalence at the population level for $d=1$.} 
    \label{fig_sim_population}
\end{figure}
\begin{figure}[t]
    \centering
    \includegraphics[width=0.9\textwidth]{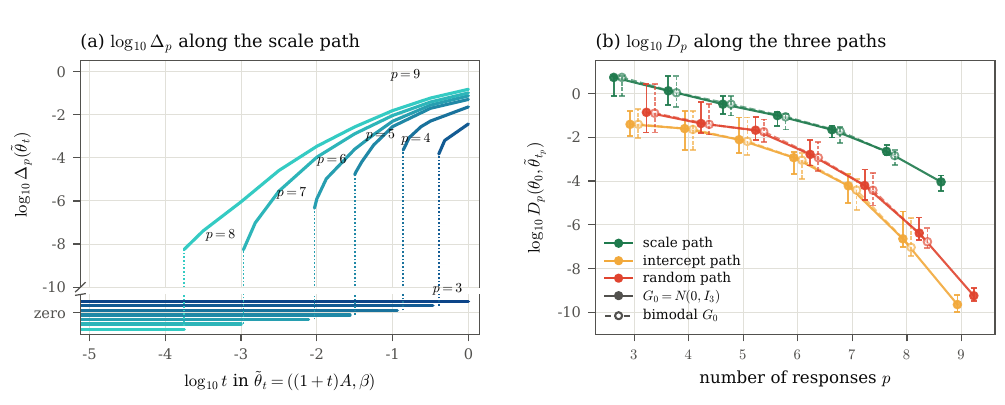}
    \caption{Exact equivalence at the population level for $d=3$. }
    \label{fig_sim_population_d3}
\end{figure}

We examine three perturbation paths: the scale path $\tilde\theta_t=((1+t)A,\beta)$, the intercept path $\tilde\theta_t=(A,\beta+te_1)$, and the random path $\tilde\theta_t=\theta_0+tU$, where $U$ is uniform on the unit sphere of $\RR^{p(d+1)}$. Along each path, we estimate $t_p$ as the transition at numerically zero value of $\Delta_p(\tilde\theta_t)$. The estimated transitions illustrate the ambiguity along corresponding directions, which serve as a proxy for $\omega_p$. We adopt this proxy because computing the full diameter $\omega_p$ would require optimizing over all admissible pairs and can be computationally infeasible. Figures~\ref{fig_sim_population} and~\ref{fig_sim_population_d3} summarize the results for $d=1$ and $d=3$, respectively. 


Figure~\ref{fig_sim_population} shows that the identification ambiguity decays rapidly as $p$ increases when $d=1$. Panel (a) illustrates the discrepancy $\Delta_p(\tilde\theta_t)$ along the scale path over a high-resolution grid; see Appendix~\ref{app_sim}for detail. The flat segments in the bottom band mark intervals of numerically zero discrepancy, illustrating the non-identifiability established in Theorem~\ref{thm_nonidentifiability}. These intervals contract sharply, with the endpoint approximately $t_6\approx10^{-34}$ at $p=6$, indicating that the ambiguity is extremely small when $d=1$. Panel (b) shows a similar decline in $D_p(\theta_0,\tilde\theta_{t_p})$ along all three paths, with comparable patterns under the Gaussian and bimodal latent distributions.
Figure~\ref{fig_sim_population_d3} shows a similar pattern for $d=3$, with the ambiguity larger over the displayed range. This can be explained by the greater flexibility resulting from a higher-dimensional latent distribution, which may require more observed responses to achieve a comparable scale of identification ambiguity.  
Overall, numerical results illustrate how the ambiguity depends on the latent dimension, validating the vanishing ambiguity  as $p$ grows, established by Theorem~\ref{thm_asymptotic}.


\section{Conclusion}\label{sec_conclude}

This paper studies identifiability in latent variable models with flexible latent distributions. We show that, when the latent distribution is unknown and learned jointly with the model parameters, exact identifiability fails for almost every parameter, even under moment constraints of any prespecified finite order. Despite this non-identifiability, the model parameters remain estimable. Specifically, we establish that the resulting identification ambiguity vanishes uniformly over admissible parameters at the superpolynomial rate $\exp\{-c(\log p)^2\}$, where $p$ is the number of observed variables. Consequently, any parameter that induces the same marginal distribution as the true parameter can serve as a population target for consistent estimation and valid inference. These results provide a theoretical basis for jointly estimating model parameters and latent distributions as both the sample size and $p$ grow.



Our work opens up several possible future research directions. First, an important next step is to develop joint estimators of the model parameters and latent distribution that converge to population solutions inducing the true marginal distribution, as required by Corollaries~\ref{cor_estimation} and~\ref{prop_prediction}. Flexible learned priors, including mixture priors~\citep{dunson2009nonparametric,tomczak2018vae} and normalizing flows~\citep{papamakarios2019flows}, offer promising constructions. Second, our proof techniques may provide a starting point for studying asymptotic identifiability under more general nonlinear decoders and more complex latent structures. Third, from a practical perspective, it would be valuable to develop computable diagnostics that quantify the identification ambiguity at a given finite $p$, which could guide applied researchers in deciding how many observed variables are needed for the estimated parameters to be reliably interpreted.

\bibliography{reference}

@article{kingma2013auto,
  title={Auto-encoding variational bayes},
  author={Kingma, Diederik P and Welling, Max},
  journal={arXiv preprint arXiv:1312.6114},
  year={2013}
}

@article{papamakarios2019flows,
  title={Normalizing flows for probabilistic modeling and inference},
  author={Papamakarios, George and Nalisnick, Eric and Rezende, Danilo Jimenez and Mohamed, Shakir and Lakshminarayanan, Balaji},
  journal={Journal of Machine Learning Research},
  volume={22},
  number={57},
  pages={1--64},
  year={2021}
}

@article{goodfellow2014generative,
  title={Generative adversarial nets},
  author={Goodfellow, Ian J and Pouget-Abadie, Jean and Mirza, Mehdi and Xu, Bing and Warde-Farley, David and Ozair, Sherjil and Courville, Aaron and Bengio, Yoshua},
  journal={Advances in Neural Information Processing Systems},
  volume={27},
  year={2014}
}

@article{ho2020denoising,
  title={Denoising diffusion probabilistic models},
  author={Ho, Jonathan and Jain, Ajay and Abbeel, Pieter},
  journal={Advances in Neural Information Processing Systems},
  volume={33},
  pages={6840--6851},
  year={2020}
}

@article{fang2021identifiability,
  title={Identifiability of bifactor models},
  author={Fang, Guanhua and Guo, Jinxin and Xu, Xin and Ying, Zhiliang and Zhang, Susu},
  journal={Statistica Sinica},
  volume={31},
  pages={2309--2330},
  year={2021},
  publisher={JSTOR}
}

@article{reise2012rediscovery,
  title={The rediscovery of bifactor measurement models},
  author={Reise, Steven P},
  journal={Multivariate Behavioral Research},
  volume={47},
  number={5},
  pages={667--696},
  year={2012},
  publisher={Taylor \& Francis}
}

@inproceedings{xi2023indeterminacy,
  title={Indeterminacy in generative models: Characterization and strong identifiability},
  author={Xi, Quanhan and Bloem-Reddy, Benjamin},
  booktitle={International Conference on Artificial Intelligence and Statistics},
  pages={6912--6939},
  year={2023},
  organization={PMLR}
}

@article{hyvarinen1999nonlinear,
  title={Nonlinear independent component analysis: Existence and uniqueness results},
  author={Hyv{\"a}rinen, Aapo and Pajunen, Petteri},
  journal={Neural Networks},
  volume={12},
  number={3},
  pages={429--439},
  year={1999},
  publisher={Elsevier}
}

@article{onatski2012asymptotics,
  title={Asymptotics of the principal components estimator of large factor models with weakly influential factors},
  author={Onatski, Alexei},
  journal={Journal of Econometrics},
  volume={168},
  number={2},
  pages={244--258},
  year={2012},
  publisher={Elsevier}
}

@inproceedings{tomczak2018vae,
  title={{VAE} with a {V}amp{P}rior},
  author={Tomczak, Jakub and Welling, Max},
  booktitle={International Conference on Artificial Intelligence and Statistics},
  pages={1214--1223},
  year={2018},
  organization={PMLR}
}

@article{xu2023identifiable,
  title={Identifiable and interpretable nonparametric factor analysis},
  author={Xu, Maoran and Winter, Steven and Herring, Amy H and Dunson, David B},
  journal={arXiv preprint arXiv:2311.08254},
  year={2023}
}

@article{dunson2009nonparametric,
  title={Nonparametric {B}ayes modeling of multivariate categorical data},
  author={Dunson, David B and Xing, Chuanhua},
  journal={Journal of the American Statistical Association},
  volume={104},
  number={487},
  pages={1042--1051},
  year={2009},
  publisher={Taylor \& Francis}
}

@article{withers1984asymptotic,
  title={Asymptotic expansions for distributions and quantiles with power series cumulants},
  author={Withers, Christopher S},
  journal={Journal of the Royal Statistical Society Series B: Statistical Methodology},
  volume={46},
  number={3},
  pages={389--396},
  year={1984},
  publisher={Oxford University Press}
}

@article{ostrovskii2021finite,
  title={Finite-sample analysis of {M}-estimators using self-concordance},
  author={Ostrovskii, Dmitrii M and Bach, Francis},
  journal={Electronic Journal of Statistics},
  volume={15},
  number={1},
  pages={326--391},
  year={2021}
}

@book{lawson1995solving,
  title={Solving least squares problems},
  author={Lawson, Charles L and Hanson, Richard J},
  year={1995},
  publisher={SIAM}
}

@article{davezies2021identification,
  title={Identification and estimation of average causal effects in fixed effects logit models},
  author={Davezies, Laurent and D'Haultf{\oe}uille, Xavier and Laage, Louise},
  journal={arXiv preprint arXiv:2105.00879},
  year={2021}
}

@article{chernozhukov2013average,
  title={Average and quantile effects in nonseparable panel models},
  author={Chernozhukov, Victor and Fern{\'a}ndez-Val, Iv{\'a}n and Hahn, Jinyong and Newey, Whitney},
  journal={Econometrica},
  volume={81},
  number={2},
  pages={535--580},
  year={2013},
  publisher={Wiley Online Library}
}

@article{kolda2009tensor,
  title={Tensor decompositions and applications},
  author={Kolda, Tamara G and Bader, Brett W},
  journal={SIAM Review},
  volume={51},
  number={3},
  pages={455--500},
  year={2009},
  publisher={SIAM}
}

@book{stein2003complex,
  title={Complex analysis},
  author={Stein, Elias M and Shakarchi, Rami},
  volume={2},
  year={2003},
  publisher={Princeton University Press, Princeton}
}

@article{hinton1995wake,
  title={The" wake-sleep" algorithm for unsupervised neural networks},
  author={Hinton, Geoffrey E and Dayan, Peter and Frey, Brendan J and Neal, Radford M},
  journal={Science},
  volume={268},
  number={5214},
  pages={1158--1161},
  year={1995},
  publisher={American Association for the Advancement of Science}
}

@article{hyvarinen2000independent,
  title={Independent component analysis: algorithms and applications},
  author={Hyv{\"a}rinen, Aapo and Oja, Erkki},
  journal={Neural Networks},
  volume={13},
  number={4-5},
  pages={411--430},
  year={2000},
  publisher={Elsevier}
}

@article{hyvarinen2016unsupervised,
  title={Unsupervised feature extraction by time-contrastive learning and nonlinear ica},
  author={Hyvarinen, Aapo and Morioka, Hiroshi},
  journal={Advances in Neural Information Processing Systems},
  volume={29},
  year={2016}
}

@article{kivva2022identifiability,
  title={Identifiability of deep generative models without auxiliary information},
  author={Kivva, Bohdan and Rajendran, Goutham and Ravikumar, Pradeep and Aragam, Bryon},
  journal={Advances in Neural Information Processing Systems},
  volume={35},
  pages={15687--15701},
  year={2022}
}

@inproceedings{locatello2019challenging,
  title={Challenging common assumptions in the unsupervised learning of disentangled representations},
  author={Locatello, Francesco and Bauer, Stefan and Lucic, Mario and Raetsch, Gunnar and Gelly, Sylvain and Sch{\"o}lkopf, Bernhard and Bachem, Olivier},
  booktitle={International Conference on Machine Learning},
  pages={4114--4124},
  year={2019},
  organization={PMLR}
}

@article{zhang2001linear,
  title={Linear mixed models with flexible distributions of random effects for longitudinal data},
  author={Zhang, Daowen and Davidian, Marie},
  journal={Biometrics},
  volume={57},
  number={3},
  pages={795--802},
  year={2001},
  publisher={Wiley Online Library}
}

@article{lubke2005investigating,
  title={Investigating population heterogeneity with factor mixture models.},
  author={Lubke, Gitta H and Muth{\'e}n, Bengt},
  journal={Psychological Methods},
  volume={10},
  number={1},
  pages={21--39},
  year={2005},
  publisher={American Psychological Association}
}

@book{reckase2009multidimensional,
  title={Multidimensional Item Response Theory},
  author={Reckase, Mark D.},
  year={2009},
  publisher={Springer New York, NY}
}

@article{goodman1974exploratory,
  title={Exploratory latent structure analysis using both identifiable and unidentifiable models},
  author={Goodman, Leo A},
  journal={Biometrika},
  volume={61},
  number={2},
  pages={215--231},
  year={1974},
  publisher={Oxford University Press}
}

@article{allman2009identifiability,
  title={Identifiability of parameters in latent structure models with many observed variables},
  author={Allman, Elizabeth S and Matias, Catherine and Rhodes, John A},
  journal={The Annals of Statistics},
  volume={37},
  number={6A},
  pages={3099--3132},
  year={2009}
}

@article{gu2020partial,
  title={Partial identifiability of restricted latent class models},
  author={Gu, Yuqi and Xu, Gongjun},
  journal={The Annals of Statistics},
  volume={48},
  number={4},
  pages={2082--2107},
  year={2020},
  publisher={JSTOR}
}

@article{rothenberg1971identification,
  title={Identification in parametric models},
  author={Rothenberg, Thomas J},
  journal={Econometrica: Journal of the Econometric Society},
  volume={39},
  number={3},
  pages={577--591},
  year={1971},
  publisher={JSTOR}
}

@article{kaiser1958varimax,
  title={The varimax criterion for analytic rotation in factor analysis},
  author={Kaiser, Henry F},
  journal={Psychometrika},
  volume={23},
  number={3},
  pages={187--200},
  year={1958},
  publisher={Springer}
}

@article{Rohe2020VintageFA,
  title={Vintage Factor Analysis with {V}arimax Performs Statistical Inference},
  author={Karl Rohe and Muzhe Zeng},
  journal={Journal of the Royal Statistical Society Series B: Statistical Methodology},
  volume={85},
  number={4},
  pages={1037-1060},
  year={2023}
}

@book{andersonintroduction,
  title={An introduction to multivariate statistical analysis},
  author={Anderson, Theodore Wilbur},
   year={2003},
  publisher={Wiley New York}
}

@article{rota1994classical,
  title={The classical umbral calculus},
  author={Rota, G-C and Taylor, Brian D},
  journal={SIAM Journal on Mathematical Analysis},
  volume={25},
  number={2},
  pages={694--711},
  year={1994},
  publisher={SIAM}
}

@article{li2026estimating,
  title={Estimating Latent Distribution of Item Response Theory Using Kernel Density Method},
  author={Li, Seewoo and Lee, Guemin},
  journal={Psychometrika},
  volume={91},
  number={2},
  pages={620--640},
  year={2026},
  publisher={Cambridge University Press}
}

@article{ma2010explicit,
  title={Explicit estimating equations for semiparametric generalized linear latent variable models},
  author={Ma, Yanyuan and Genton, Marc G},
  journal={Journal of the Royal Statistical Society Series B: Statistical Methodology},
  volume={72},
  number={4},
  pages={475--495},
  year={2010},
  publisher={Oxford University Press}
}

@article{san2011bayesian,
  title={On the {B}ayesian nonparametric generalization of {IRT}-type models},
  author={San Mart{\'\i}n, Ernesto and Jara, Alejandro and Rolin, Jean-Marie and Mouchart, Michel},
  journal={Psychometrika},
  volume={76},
  number={3},
  pages={385--409},
  year={2011},
  publisher={Cambridge University Press \& Assessment}
}

@inproceedings{khemakhem2020variational,
  title={Variational autoencoders and nonlinear {ICA}: A unifying framework},
  author={Khemakhem, Ilyes and Kingma, Diederik and Monti, Ricardo and Hyvarinen, Aapo},
  booktitle={International Conference on Artificial Intelligence and Statistics},
  pages={2207--2217},
  year={2020},
  organization={PMLR}
}

@article{scholkopf2021toward,
  title={Toward causal representation learning},
  author={Sch{\"o}lkopf, Bernhard and Locatello, Francesco and Bauer, Stefan and Ke, Nan Rosemary and Kalchbrenner, Nal and Goyal, Anirudh and Bengio, Yoshua},
  journal={Proceedings of the IEEE},
  volume={109},
  number={5},
  pages={612--634},
  year={2021},
  publisher={IEEE}
}

@article{lindsay1991semiparametric,
  title={Semiparametric estimation in the {R}asch model and related exponential response models, including a simple latent class model for item analysis},
  author={Lindsay, Bruce and Clogg, Clifford C and Grego, John},
  journal={Journal of the American Statistical Association},
  volume={86},
  number={413},
  pages={96--107},
  year={1991},
  publisher={Taylor \& Francis}
}

@article{haberman2005identifiability,
  title={Identifiability of parameters in item response models with unconstrained ability distributions},
  author={Haberman, Shelby J},
  journal={ETS Research Report Series},
  volume={2005},
  number={2},
  pages={1--22},
  year={2005},
  publisher={Wiley Online Library}
}

@article{cui2026beyond,
  title={Beyond vintage rotation: Bias-free sparse representation learning with oracle inference},
  author={Cui, Chengyu and Chen, Yunxiao and Ouyang, Jing and Xu, Gongjun},
  journal={arXiv preprint arXiv:2602.22590},
  year={2026}
}

@article{browne2001overview,
  title={An overview of analytic rotation in exploratory factor analysis},
  author={Browne, Michael W},
  journal={Multivariate Behavioral Research},
  volume={36},
  number={1},
  pages={111--150},
  year={2001},
  publisher={Taylor \& Francis}
}

@article{qin1994empirical,
  title={Empirical likelihood and general estimating equations},
  author={Qin, Jin and Lawless, Jerry},
  journal={The Annals of Statistics},
  volume={22},
  number={1},
  pages={300--325},
  year={1994},
  publisher={Institute of Mathematical Statistics}
}

@article{deville1992calibration,
  title={Calibration estimators in survey sampling},
  author={Deville, Jean-Claude and S{\"a}rndal, Carl-Erik},
  journal={Journal of the American Statistical Association},
  volume={87},
  number={418},
  pages={376--382},
  year={1992},
  publisher={Taylor \& Francis}
}

@article{singh2018minimax,
  title={Minimax distribution estimation in {W}asserstein distance},
  author={Singh, Shashank and P{\'o}czos, Barnab{\'a}s},
  journal={arXiv preprint arXiv:1802.08855},
  year={2018}
}

@article{Liu2023RotationTS,
  title={Rotation to Sparse Loadings Using ${L}^p$ Losses and Related Inference Problems},
  author={Xinyi Liu and Gabriel Wallin and Yunxiao Chen and Irini Moustaki},
  journal={Psychometrika},
  year={2023},
  volume={88},
  number={2},
  pages={527-553}
}

@article{ten1977orthogonal,
  title={Orthogonal {P}rocrustes rotation for two or more matrices},
  author={Ten Berge, Jos M.F.},
  journal={Psychometrika},
  volume={42},
  number={2},
  pages={267--276},
  year={1977},
  publisher={Springer-Verlag}
}

@article{arellano2009robust,
  title={Robust priors in nonlinear panel data models},
  author={Arellano, Manuel and Bonhomme, St{\'e}phane},
  journal={Econometrica},
  volume={77},
  number={2},
  pages={489--536},
  year={2009},
  publisher={Wiley Online Library}
}

@article{rizopoulos2008shared,
  title={Shared parameter models under random effects misspecification},
  author={Rizopoulos, Dimitris and Verbeke, Geert and Molenberghs, Geert},
  journal={Biometrika},
  volume={95},
  number={1},
  pages={63--74},
  year={2008},
  publisher={Oxford University Press}
}

@article{neuhaus1992effects,
  title={The effects of mixture distribution misspecification when fitting mixed-effects logistic models},
  author={Neuhaus, John M and Hauck, Walter W and Kalbfleisch, John D},
  journal={Biometrika},
  volume={79},
  number={4},
  pages={755--762},
  year={1992},
  publisher={Oxford University Press}
}

@article{anandkumar2014tensor,
  title={Tensor decompositions for learning latent variable models},
  author={Anandkumar, Animashree and Ge, Rong and Hsu, Daniel and Kakade, Sham M and Telgarsky, Matus},
  journal={The Journal of Machine Learning Research},
  volume={15},
  number={1},
  pages={2773--2832},
  year={2014},
  publisher={JMLR. org}
}

@article{antonelli2016mitigating,
  title={Mitigating bias in generalized linear mixed models: The case for {B}ayesian nonparametrics},
  author={Antonelli, Joseph and Trippa, Lorenzo and Haneuse, Sebastien},
  journal={Statistical Science},
  volume={31},
  number={1},
  pages={80--95},
  year={2016}
}

@article{cui2026multidimensional,
  title={Multidimensional Item Response Theory under General Latent Distributions},
  author={Cui, Chengyu and Chen, Taoyi and Wang, Chun and Xu, Gongjun},
  journal={arXiv preprint arXiv:2605.30471},
  year={2026}
}

@article{paganin2023computational,
  title={Computational strategies and estimation performance with {B}ayesian semiparametric item response theory models},
  author={Paganin, Sally and Paciorek, Christopher J and Wehrhahn, Claudia and Rodr{\'\i}guez, Abel and Rabe-Hesketh, Sophia and de Valpine, Perry},
  journal={Journal of Educational and Behavioral Statistics},
  volume={48},
  number={2},
  pages={147--188},
  year={2023},
  publisher={Sage Publications Sage CA: Los Angeles, CA}
}

@article{woods2006item,
  title={Item response theory with estimation of the latent population distribution using spline-based densities},
  author={Woods, Carol M and Thissen, David},
  journal={Psychometrika},
  volume={71},
  number={2},
  pages={281--301},
  year={2006},
  publisher={Cambridge University Press \& Assessment}
}

@article{cui2026marginal,
  title={Marginal maximum likelihood estimation and asymptotic theory for latent variable models in high dimensions},
  author={Cui, Chengyu and Xu, Gongjun},
  journal={arXiv preprint arXiv:2609.16653},
  year={2026}
}

@article{douglas2001asymptotic,
  title={Asymptotic identifiability of nonparametric item response models},
  author={Douglas, Jeffrey A},
  journal={Psychometrika},
  volume={66},
  number={4},
  pages={531--540},
  year={2001},
  publisher={Springer}
}

@article{he2024extended,
  title={Extended Asymptotic Identifiability of Nonparametric Item Response Models},
  author={He, Yinqiu},
  journal={Psychometrika},
  volume={89},
  number={3},
  pages={958--973},
  year={2024},
  publisher={Springer}
}

@article{xu2017,
author = "Xu, Gongjun",
journal = "The Annals of Statistics",
number = "2",
pages = "675--707",
publisher = "The Institute of Mathematical Statistics",
title = "Identifiability of restricted latent class models with binary responses",
volume = "45",
year = "2017"
}

@article{cui2025identifiability,
  title={Identifiability and inference for generalized latent factor models},
  author={Cui, Chengyu and Xu, Gongjun},
  journal={Journal of the American Statistical Association},
  number={just-accepted},
  pages={1--29},
  year={2026},
  publisher={Taylor \& Francis}
}

@book{collins2013latent,
  title={Latent class and latent transition analysis: With applications in the social, behavioral, and health sciences},
  author={Collins, Linda M and Lanza, Stephanie T},
  year={2013},
  publisher={John Wiley \& Sons}
}

@book{skrondal2004generalized,
  title={Generalized latent variable modeling: Multilevel, longitudinal, and structural equation models},
  author={Skrondal, Anders and Rabe-Hesketh, Sophia},
  year={2004},
  publisher={Chapman and Hall/CRC}
}

@article{bock1981marginal,
  title={Marginal maximum likelihood estimation of item parameters: Application of an {EM} algorithm},
  author={Bock, R Darrell and Aitkin, Murray},
  journal={Psychometrika},
  volume={46},
  number={4},
  pages={443--459},
  year={1981},
  publisher={Springer-Verlag}
}
\clearpage

\appendix
\renewcommand{\theequation}{S\arabic{equation}} 
\setcounter{equation}{0}  

\renewcommand{\theproposition}{\thesection\arabic{proposition}}
\renewcommand{\thelemma}{\thesection.\arabic{lemma}}
\renewcommand{\thecorollary}{\thesection\arabic{corollary}}


\section{Preliminaries}

We introduce the following additional notation. The unit sphere in $\RR^d$ is $\mathbb S^{d-1}=\{u\in\RR^d:\|u\|=1\}$, and $e_i^{(n)}$ denotes the $i$th standard basis vector in $n$-dimensional Euclidean space, with the superscript $^{(n)}$ omitted when the dimension is clear from the context. We write $\mathbb N_0=\{0,1,2,\ldots\}$.  For a matrix $M\in\RR^{n\times m}$, we write $\|M\|_{\mathrm F}=(\sum_{i,j}M_{ij}^2)^{1/2}$ for the Frobenius norm. For index sets $\mathcal I\subseteq[n]$, $\mathcal J\subseteq[m]$, write $(M)_{\mathcal I,\mathcal J}=(M_{ij})_{i\in\mathcal I,j\in\mathcal J}$. When $\cI = [n]$, we write $(M)_{,\mathcal{J}}$ for short and $(M)_{\cI,}$ when $\mathcal{J}=[m]$. 

For a function $h$, $\|h\|_\infty$ denotes its supremum norm. Recall that $C^\infty$ consists of functions with continuous derivatives of all orders. We further let $C_c^\infty$ denote those with compact support. A function on an open set is analytic if, near every point, it agrees with a convergent power series about that point. A Lebesgue zero set, or null set, is a set of Lebesgue measure zero. An integral with respect to a latent variable in $\RR^d$ is written as $\int \cdot d z$ unless a domain is specified. For a random vector or scalar $W$ and $1\le q<\infty$, write \begin{equation*}\|W\|_q=(\EE\|W\|^q)^{1/q},\end{equation*}where the norm inside the expectation is the Euclidean norm for vectors and absolute value for scalars. We distinguish this $L_q$ norm from the unsubscripted notation $\norm{\cdot}$, which denotes the Euclidean norm for vectors and the operator norm for matrices.

We also introduce basic polynomial notation. For a commutative coefficient ring $\mathcal A$, the polynomial ring $\mathcal A[x_1,\ldots,x_r]$ consists of finite sums $P(x)=\sum_a c_a x^a$, where $c_a\in\mathcal A$, $a\in\mathbb N_0^r$, and $x^a=\prod_{j=1}^r x_j^{a_j}$ is a monomial of total degree $|a|$. Equality of polynomials means equality of the coefficients of every monomial. An outer polynomial $R(x_1,\ldots,x_r)$ is a polynomial whose arguments may themselves be polynomials. For instance, if $F_1(u),\ldots,F_r(u)$ are polynomials, then $R(F_1(u),\ldots,F_r(u))$ denotes an outer polynomial of $u$ obtained by replacing $x_i$ by $F_i(u)$. For example, $R(x_1,x_2)=x_1^2+x_2$ gives $R(F_1(u),F_2(u))=F_1(u)^2+F_2(u)$.

%




\section{Proofs of Main Results}\label{app_proofs}
\renewcommand{\theequation}{\thesection.\arabic{equation}}
\setcounter{equation}{0}

\subsection{Proof of Theorem~\ref{thm_nonidentifiability}}
\label{app_nonidentifiability}

Given $(\theta,G)$, we construct perturbation of $G$ to offset any sufficiently small perturbation of $\theta$ and preserves the moments such that $\EE_{\tilde G}[Z^a]=\EE_G[Z^a]$ for all $|a|\le 2s$. The construction builds on the idea of moment-matching structure of calibration weighting and empirical likelihood
\citep{deville1992calibration,qin1994empirical}.

Let $q_S(\theta,z)=\prod_{j\in S}f(\beta_j+\alpha_j^{\T}z)$ for any $S\subseteq [p]$. When $S = \varnothing$, let $q_{\varnothing}(\theta,z) = 1$. By conditional independence, $\EE_{\theta,G}\prod_{j\in S}X_j=\int q_S(\theta,z)\,dG(z)$.
For $x\in\{0,1\}^p$ and $S_x=\{j:x_j=1\}$, inclusion--exclusion gives
\begin{equation}\label{eq_subset_law}
 P_{\theta,G}(X=x)
 =\sum_{T\subseteq[p]\setminus S_x}(-1)^{|T|}
   \int q_{S_x\cup T}(\theta,z)\,dG(z).
\end{equation}
Thus, it suffices to construct $(\theta,G)$ and $(\tilde\theta,\tilde G)$ such that $\int q_S(\theta,z)\,dG(z)=\int q_S(\tilde\theta,z)\,d\tilde G(z)$ for every $S\ne\varnothing$. The proof proceeds in two steps. First, we define a Lebesgue-null exceptional set $\cN_p$ with the required property. Second, for any $\theta\notin\cN_p$, any sufficiently close $\tilde\theta$, and any admissible distribution $G$ with density $g$, we construct a perturbed density of the form $g(z)(1 + \mathsf{ correction})$ so that the distribution $\tilde G$ satisfies $\int q_S(\theta,z)\,dG(z)=\int q_S(\tilde\theta,z)\,d\tilde G(z)$, as well as other requirements in Theorem~\ref{thm_nonidentifiability}.

\medskip

\paragraph{\underline{(1) Defining the exceptional set.}}
We first construct a Lebesgue null set $\cN_p$ such that, for any $\theta\notin\cN_p$, $q_S(\theta,z)$ are linearly independent as functions of $z$, namely, $\sum_{S\subseteq[p]}c_Sq_S(\theta,z)=0$ for every $z\in\mathbb R^d$ implies $c_S=0$ for every $S\subseteq[p]$. It then suffices to find one $\theta^*$ for which $Q_S(\theta^*,r)$, defined as 
\begin{equation*}
    Q_S(\theta,r)=\int_{\mathbb S^{d-1}}q_S(\theta,ru)\,d\sigma_{d-1}(u), \;S\subseteq[p]
\end{equation*}are linear independent as functions of $r>0$. Here, $\sigma_{d-1}$ is uniform probability measure on $\mathbb S^{d-1}$
($\sigma_0=(\delta_{-1}+\delta_1)/2$). To see this, suppose such $\theta^*$ exists. Then the $q_S(\theta^*,\cdot)$ are linearly independent, since any relation among them remains a relation after spherical averaging to $Q_{S}(\theta,r)$, so we can find $z_1,\dots,z_{2^p}\in\RR^d$ such that $\det\big(q_S(\theta^*,z_i)\big)_{1\le i\le2^p,\,S\subseteq[p]}\neq 0$. Next, define
\begin{equation}
    \cN_p = \big\{\theta:\det\big\{q_S(\theta,z_i)\big\}_{1\le i\le2^p,\,S\subseteq[p]}=0\big\}.\label{eq_precise_define_cNp}
\end{equation}
Since $f$ is analytic, $\det\big\{q_S(\theta,z_i)\big\}_{1\le i\le2^p,\,S\subseteq[p]}$ is also analytic in $\theta$. By the definition of $\theta^*$ and  $\{z_i\}_{i\in[2^p]}$, the determinant is nonzero at $\theta^*$, and therefore its zero set $\cN_p$ has Lebesgue measure zero. For $\theta\notin\cN_p$, matrix $\big(q_S(\theta,z_i)\big)_{i\in[2^p],S\subseteq[p]}$ is nonsingular, so the $q_S(\theta,\cdot)$ are linearly independent.

To show that $Q_S(\theta^*,r)$ are linear independent, we design $\theta^*$ inductively. Suppose $A^*$ in $\theta^*$ has only the first column nonzero, and let 
\begin{equation*}
    v_S(u)=\prod_{j\in S}f(\beta_j^*+\alpha_{j1}^*u)\quad\text{and}\quad
 \bar v_S(r)=\int_{\mathbb S^{d-1}}v_S(ru_1)\,d\sigma_{d-1}(u).
\end{equation*}Here, $u_1$ is the first component in $u$. We then choose $\beta_j^*\in\RR$ and $\alpha_{j1}^*\in\RR$ inductively so that $(\bar v_S)_{S\subseteq[k]}$ is linearly independent. 

When $k=0$, $\bar v_{\varnothing} = 1$, and there is nothing to prove. Suppose the assertion holds for $k\ge0$. Let $n=2^k$. Then $(v_S)_{S\subseteq[k]}$ is also
linearly independent due to the spherical averaging. 
 By the induction hypothesis there are radii $0<r_1<\cdots<r_n$ such that $(\bar v_S(r_i))_{i\le n,S\subseteq[k]}$ is nonsingular, and we fix $t_0>r_n$ as a threshold for new functions when moving from $k$ to $k+1$:
\begin{equation*}
\overline v_{S\cup\{k+1\}}(r)
=\int_{\mathbb S^{d-1}}
   v_S(ru_1)f_\ell(ru_1)\,d\sigma_{d-1}(u),
\qquad S\subseteq[k],
\end{equation*}
For $r\neq t_0$, let $\mathcal T_dv_S(r)
 :=\int_{ru_1>t_0}v_S(ru_1)\,d\sigma_{d-1}(u)$, which equals $0$ when $0<r<t_0$.
For $j=k+1$ and any $\ell>0$, define $\alpha_{j1}^*=\ell$, $\beta_{j}^*=-\ell t_0$, and $f_\ell(t)=f\{\ell(t-t_0)\}$. Then as $\ell\to\infty$, $f_{\ell}(t)\to 1_{\{t>t_0\}}$ for $t\neq t_0$, and thus since $\{u:ru_1=t_0\}$ is null under measure $\sigma_{d-1}$ for $r\neq t_0$, we have convergence $\int_{\mathbb S^{d-1}}v_S(ru_1)f_\ell(ru_1)\,d\sigma_{d-1}(u)\to \mathcal T_dv_S(r)$.

 We then show that $(\mathcal T_dv_S)_{S\subseteq[k]}$ are linearly independent on $(t_0,\infty)$. Let $v(u)=\sum_{S\subseteq[k]}c_Sv_S(u)$ with $c\neq0$. Since the $v_S(u)$ are linearly independent and analytic, $v(u)\not\equiv0$ has isolated zeros. Thus for some $\delta>0$, $v(u)$ has no zero within $(t_0,t_0+\delta)$, namely, it has a constant sign on that interval. For $r\in(t_0,t_0+\delta)$, the integrand of $\mathcal T_dv(r)=\int_{u_1>t_0/r}v(ru_1)\,d\sigma_{d-1}(u)$ has fixed nonzero sign throughout the integration domain $u>t_0/r$, within which we have $t_0< ru_1\le r <t_0+\delta$ as $u_1\le 1$ when $u\in \mathbb{S}^{d-1}$. This domain has positive measure under $\sigma_{d-1}$-measure, so $\mathcal T_dv(r)\neq0$.

 With $r_1<\cdots<r_n<t_0$ as above, we choose $t_0<r_{n+1}<\cdots<r_{2n}$ such that $(\mathcal T_dv_S(r_{n+i}))_{i\le n,S\subseteq[k]}$ is nonsingular. Note that $\bigl(\bar v_S(r_i)\bigr)_{1\le i\le 2^{k+1},\;S\subseteq[k+1]}$ converges to 
 \begin{equation*}
     \begin{pmatrix}(\bar v_S(r_i))_{i\le n,S\subseteq[k]} & 0\\ (\bar v_S(r_{n+i}))_{i\le n,S\subseteq[k]} & (\mathcal T_dv_S(r_{n+i}))_{i\le n,S\subseteq[k]} \end{pmatrix},
 \end{equation*}whose determinant is nonzero by the linear independence. Consequently, for sufficiently large $\ell$, $(\bar v_S)_{S\subseteq[k+1]}$ is  linearly independent on $(0,\infty)$, which completes the induction.

 Therefore, there is some $\theta^*$ such that $Q_S(\theta^*,r)$ are  linearly independent, and by the preceding discussion, $\cN_p$ constructed in \eqref{eq_precise_define_cNp} is the required Lebesgue null set.

\medskip

\paragraph{\underline{(2) Constructing the perturbed density.}}
Now fix $\theta\notin\cN_p$ and any $G$ satisfying \eqref{eq_unv_normalize}. Let $g\in C^\infty$ be the density of $G$. Fix $s\ge 1$ under which $\EE_G\|Z\|^{2s}<\infty$. The case $s=1$ proves the first assertion of the theorem.
 By definition, we know
 \begin{equation}
     \{z^a:|a|\le2s\}\ \cup\{q_S(\theta,z):\varnothing\ne S\subseteq[p]\}\label{eq_independent_constraints}
 \end{equation} is linear independent, as the first part concerns polynomials of $z$ and the second part concerns bounded functions of linear combination in $q_S(\theta,z)$. For simplicity, we write $2^{[p]} = \{S_0 = \varnothing, S_1,\dots,S_{n-1}\}$ where we slightly abuse the notation to denote $n = 2^p$. For any $w\in C_c^{\infty}(\RR^d)$, define the linear map 
 \begin{equation*}
     w(\cdot)\mapsto \Big[\big\{\int z^a w(z)g(z)\,dz\big\}_{|a|\le2s},
       \big\{\int q_{S_i}(\theta,z)w(z)g(z)\,dz\big\}_{i=1}^{n-1}\Big].
 \end{equation*}
This map is surjective. Otherwise, since the map is linear, the range of the map has dimension strictly less than $|\{a:|a|\le 2s\}| + n-1$. Then we can find some nonzero linear combination of the functions in \eqref{eq_independent_constraints} that integrates to zero for every $w(\cdot)g(\cdot)$. Since the density $g(\cdot)$ is positive, the combination would be zero, leading to the contradiction. The surjectivity then implies that there exists some $w_{l}\in { C_c^{\infty}(\RR^d)}$ mapped to $(0,e_l)$ for any $l\in[n-1]$, which thus satisfies
\begin{equation}
    \int z^a w_l(z)g(z)\,dz=0\quad\text{and}\quad \int q_{S_i}(\theta,z)w_l(z)g(z)\,dz=\delta_{il},\text{ for each }i\in[n-1]\label{eq_dual_bumps}
\end{equation}

The collection $\{w_l\}_{l\in[n-1]}$ serves the role of bump function $w(z)$ in Example~\ref{example_nonid}. For $\tilde\theta$, let 
\begin{equation*}
    M(\tilde\theta) = \big\{M_{il}(\tilde\theta)\big\}_{[n-1],[n-1]}\;\text{ and }\; b(\tilde\theta) = \big\{b_i(\tilde\theta)\big\}_{[n-1]}
\end{equation*}
where for each $i,l\in[n-1]$,
    $M_{il}(\tilde\theta)=\int q_{S_i}(\tilde\theta,z)w_l(z)g(z)\,dz$ and $b_i(\tilde\theta)=\int\{q_{S_i}(\theta,z)-q_{S_i}(\tilde\theta,z)\}g(z)\,dz$.
By definition, as $\tilde\theta\to \theta$, $M(\tilde\theta) \to M(\theta) = I_{n-1}$ and $b_i(\tilde\theta)\to b_i(\theta) = 0$. Then on a neighbourhood $\cU$ of $\theta$, we let 
\begin{equation}
    \tilde g(z)=g(z)\Big[1+\sum_{l=1}^{n-1}  \big[\{M(\tilde\theta)\}^{-1}b(\tilde\theta)\big]_l\,w_l(z)\Big].\label{eq_density_compensation}
\end{equation}
Because  $\{M(\tilde\theta)\}^{-1}b(\tilde\theta)\to 0$ as $\tilde\theta\to \theta$, we can shrink $\cU$ so that  $\|\sum_{l=1}^{n-1} [\{M(\tilde\theta)\}^{-1}b(\tilde\theta)]_l\,w_l\|_\infty<1/2$. This neighbourhood $\mathcal{U}$ depends on the chosen $s$, but we suppress this dependence in the notation for simplicity. Then $\tilde g>0$ and $\tilde g\in C^\infty(\RR^d)$. By \eqref{eq_dual_bumps}, we know that $\int z^a\tilde g(z)\,dz=\int z^a g(z)\,dz$ and similar to \eqref{eq_example_constraint_matching}, that
\begin{equation*}
    \int q_{S_i}(\tilde\theta,z)\tilde g(z)\,dz =\int q_{S_i}(\theta,z)g(z)\,dz.
\end{equation*}
Thus the pairs $(\theta,G)$ and $(\tilde\theta,\tilde G)$ satisfies the requirement in Theorem~\ref{thm_nonidentifiability}. Moreover, taking
$\tilde\theta=(A,\beta+u e_1)\in\cU$ with any $u\ne0$ gives $D_p(\theta,\tilde\theta)=|u|>0$, proving non-identifiability in the sense of Definition~\ref{def_id}.

\subsection{Proof of Proposition~\ref{prop_weak_factor}}
\label{app_weak_factor}
Following the definitions in Appendix~\ref{app_nonidentifiability}, let $n=2^p$ and use the dual functions $w_l$ from \eqref{eq_dual_bumps} at the fixed parameter $\theta$. Define $\theta_{\varepsilon} = ([A,\varepsilon v],\beta)$ and for $i\in[n-1]$
\begin{equation*}
    b_i^{\dagger}(\theta_{\varepsilon}) =\EE\Big[\prod_{j\in S_i} f(\beta_j+\alpha_j^{\T}Z+\varepsilon v_jU)\Big] -\EE[q_{S_i}(\theta,Z)],
\end{equation*}for $Z\sim G$ and $U\sim H$. Similar to \eqref{eq_density_compensation}, let
\begin{equation*}
    \tilde g_\varepsilon(z)=g(z)\Big[1+\sum_{l=1}^{n-1} b_l^{\dagger}(\theta_{\varepsilon})w_l(z)\Big].
\end{equation*}
 Similar to the proof of Theorem~\ref{thm_nonidentifiability}, one can check that, for sufficiently small $|\varepsilon|$, $\tilde g_\varepsilon(z)$ is a positive smooth density preserving all moments through order $2s$. For every $i$, it holds that 
\begin{equation*}
    \int q_{S_i}(\theta,z)\tilde g_\varepsilon(z)\,dz =\EE[q_{S_i}(\theta,Z)]+b_i^{\dagger}(\theta_{\varepsilon}) =\EE\Big[\prod_{j\in S_i} f(\beta_j+\alpha_j^{\T}Z+\varepsilon v_jU)\Big].
\end{equation*}
We therefore know $ P_{([A,\varepsilon v],\beta),\,G\otimes H} =P_{(A,\beta),\,\tilde G_\varepsilon}$ following a similar argument in the inclusion--exclusion argument for \eqref{eq_subset_law}. 
By the definition of $G$ and $H$, we know $G\otimes H$ has mean zero and covariance
$I_{d+1}$. The other requirements in Proposition~\ref{prop_weak_factor} can be similarly verified.

\subsection{Proof of Theorem~\ref{thm_asymptotic}}
\label{app_asymptotic}
\medskip\noindent{\underline{\it Step 0. Preliminaries.}}
To prove Theorem~\ref{thm_asymptotic}, we use a conditioning argument to compare $\theta$ and $\tilde\theta$ through high-order expansions of conditional moments. Establishing the superpolynomial rate requires careful control of the expansion terms and remainders as the expansion order increases. To organize these calculations, we introduce a symbolic calculus. We first collect the technical ingredients needed for the proof.

\smallskip
\paragraph{Notation.}
The entire proof concerns two \textbf{representations} $(\theta,G)$ and $(\tilde\theta,\tilde G)$ in $\mathcal M_p$ with $P_{\theta,G} = P_{\tilde\theta,\tilde G}$. For notational simplicity, we index them by $r=0,1$; refer them as model zero, model one; and write $\alpha_{0j}=\alpha_j$, $\alpha_{1j}=\widetilde\alpha_j$, $\beta_{0j}=\beta_j$, and $\beta_{1j}=\widetilde\beta_j$. 
The rows in $A_r$ are $\alpha_{rj}^{\T}$ and $Z_r\sim G_r$. All expectations of functions of the responses alone are under the common marginal distribution $P_{\theta,G}$, and we use $\EE_r$ to distinguish the two latent representations. 
With $H=2K$, we will frequently use the following relationships obtained from Assumption~\ref{ass_link}\ref{ass_link_2}
\begin{equation}\label{eq_linkslope}
 \left|\frac{f''(t)}{f'(t)}\right|\le H,\quad
 f'(t)\ge ce^{-H|t|},\quad\text{and}\quad
 e^{-H|u|}\le\frac{f'(t+u)}{f'(t)}\le e^{H|u|}.
\end{equation}
Write $D_p = D_p(\theta,\widetilde\theta)$ for short, when the argument is clear from the context.

\smallskip

\paragraph{Moment approximation.}
We now introduce the conditioning trick. Throughout, $N\ge4$ denotes an integer approximation order. We select a subset $\mathcal S$ of size $k=O(N)$ and condition on the retained responses $Y=X_{S^c}$. The set $\mathcal S$ contains a target set $\mathcal T$ and a reference set $\mathcal I$. The target set contains two elements $\{i,j\}$ as each time we compare $\beta_{0i},\beta_{1i}$ and products $\alpha_{0i}^{\T}\alpha_{0j}, \alpha_{1i}^{\T}\alpha_{1j}$. The reference set contains observed variables whose corresponding parameters are uniformly well-behaved in one representation, and are thereby used for robustly recovering conditional moments. At the same time, $\mathcal S$ must occupy only a small proportion of all observed variables, so that the rest responses remain informative about $Z$ and the conditional uncertainty is small. The following lemma guides the selection of these reference groups.

\begin{lemma}\label{lem_design_alignment}
Under Assumption~\ref{ass_model}, there are constants $c_0,c_1>0$ and integer $p_0$ such that for any $p\ge p_0$ and any $\mathcal{J}\subseteq[p]$, with $|\mathcal J|\le c_0p$, we have \begin{equation*}\sum_{j\in[p]\setminus\mathcal J}\alpha_{rj}\alpha_{rj}^{\T} \succeq(\kappa p/2)I_d, \end{equation*} for $r\in\{0,1\}$.
Moreover, given any integer $N$ and target set $\mathcal T\subseteq[p]$ with $|\mathcal T|\le2$ and $Nd+2\le c_0p$, one can choose disjoint reference sets $\mathcal I_1,\ldots,\mathcal I_N$, each of size $d$ and excluding $\mathcal T$, such that their row matrices $B_s: = (A_{0})_{\mathcal I_s,}$ under model zero satisfy $\sigma_{\min}(B_s)\ge c_1$. The constants are uniform over the two representations and the admissible class.

\end{lemma}
\begin{proof}
    See Appendix~\ref{app_prove_lem_design_alignment}.
\end{proof}

For a fixed target set $\mathcal T=\{i,j\}$, choose the reference groups from Lemma~\ref{lem_design_alignment} and set $\cI=\bigcup_{s=1}^N\cI_s$, $\mathcal S=\mathcal T\cup\cI$, and $Y=X_{\mathcal S^c}$. The same conditioning set is used in both representations. In the scalar illustration below, write $S=\mathcal S$.
\begin{example}
To better motivate the conditioning trick with moment approximation, we use the case $d=1$ as a toy example. For arbitrary $i,j\in[p]$, we seek to bound $\beta_j-\tilde\beta_j$ and $\alpha_i\alpha_{j}-\tilde\alpha_i\tilde\alpha_{j}$, and we have selected $i,j$ in $\cS$. Conditioning on $Y=X_{S^c}$, we decompose the latent predictor, for $i\in S$, as
\begin{align*}
    \eta_i=\alpha_i Z+\beta_i = T_i + \delta_i + \alpha_i V,
\end{align*} where
\begin{equation}\label{eq_Tdelta}
    T_i = f^{-1}{\EE(X_i\mid Y)},\qquad \delta_i = \alpha_i\EE(Z\mid Y)+\beta_i-T_i,\qquad V = Z - \EE(Z\mid Y).
\end{equation}
 This decomposition splits $\eta_i$ into two parts: $\alpha_iV$ contains the variation in $Z$ that cannot be predicted by $Y$; and the predictable part $\alpha_i\EE(Z\mid Y)+\beta_i$. When $f(\cdot)$ is linear, $\alpha_i\EE(Z\mid Y)+\beta_i = T_i$, and $\delta_i$ serves as correction arising from the nonlinearity of $f(\cdot)$.



As discussed in Remark~\ref{remark_technique}, our primary goal is to approximate $\delta_i$ and the conditional moments $M_k=\EE(V^k\mid Y)$ for $k\ge2$. Lemmas~\ref{lem_conditional_fluctuations}
and~\ref{lem_bias_exponential} below shows that $V=O_{L^q}(p^{-1/2})$ and $\delta_i=O_{L^q}(p^{-1})$ and $M_k=O_{L^q}(p^{-k/2})$  for $q\ge1$ and $k\ge 2$, where $O_{L^q}$ denotes an order bound in the $L^q$ norm. These rates therefore give only $O(p^{-1})$ bounds on the target parameter differences.
Achieving the superpolynomial bound, however, requires much more accurate common approximations to the correction terms. We first construct these under $\theta$.


To approximate $\delta_i$, rewrite the identity $\EE\{f(T_i+\alpha_iV+\delta_i)\mid Y\}=f(T_i)$ as
\begin{equation}\label{eq_conditional_mean_d1}\EE[g_{T_i}(\alpha_iV+\delta_i)\mid Y]=0,\end{equation} where $g_t(u)=\{f(t+u)-f(t)\}/f'(t)$. The denominator normalizes the coefficient of the linear term to one. Writing $c_{ki}=f^{(k)}(T_i)/{k!f'(T_i)}$, from Taylor expansion of \eqref{eq_conditional_mean_d1}, we can successively compute approximation to $\delta_i$ with $M_k$. For instance, one can expand through fourth order to get
\begin{equation}\label{eq_approximate_delta_4}
\delta_i = -c_{2i}\alpha_i^2M_2
            -c_{3i}\alpha_i^3M_3
            -c_{4i}\alpha_i^4M_4
            +(3c_{2i}c_{3i}-c_{2i}^3)\alpha_i^4M_2^2 + O_{L^q}( p^{-5/2}).\end{equation}
This shall be implemented in detail in Step 2 below. Next, we use conditional moments of distinct excluded responses following the moment method~\citep{anandkumar2014tensor} to approximate $M_k$. Choose four reference coordinates in $S\setminus\{i,j\}$ and relabel them $1,2,3,4$ without loss of generality. Define, for $r\in[4]$,
\begin{equation*}
\zeta_r=\frac{X_r-f(T_r)}{f'(T_r)},\quad\text{ such that }\quad\EE(\zeta_r\mid Z,Y)=g_{T_r}(\alpha_rV+\delta_r). ,\end{equation*}
Then each $O_k=(\prod_{r=1}^k\alpha_{r})^{-1}{\EE(\prod_{r=1}^k\zeta_{r}\mid Y)}
           $ has leading term is $M_k$ for $k=2,3,4$. Specifically, expanding through fourth order and letting $a_r=c_{2r}\alpha_r$ and $b_r=c_{3r}\alpha_r^2$, we obtain
\begin{equation*}
\begin{aligned}
O_2={}&M_2+(a_1+a_2)M_3
 +(b_1+b_2+a_1a_2)M_4 - (2a_1^2+2a_2^2+a_1a_2)M_2^2
 +O_{L^q}(p^{-5/2}),\\
O_3={}&M_3+(a_1+a_2+a_3)(M_4-M_2^2)
 +O_{L^q}(p^{-5/2}),\\
O_4={}&M_4+O_{L^q}(p^{-5/2}),\quad\text{and}\quad O_2^2=M_2^2+O_{L^q}(p^{-5/2}).
\end{aligned}
\end{equation*}
To eliminate the higher order terms contained in $O_2,O_3$, we successively define
\begin{equation}\label{eq_construct_example}
\begin{aligned}
\hat M_4= {}&O_4 \qquad \hat M_3 ={}O_3-(a_1+a_2+a_3)(O_4-O_2^2),\\
\hat M_2 ={}&O_2-(a_1+a_2)\hat M_3 -(b_1+b_2+a_1a_2)O_4 +(2a_1^2+2a_2^2+a_1a_2)O_2^2.
\end{aligned}
\end{equation}
Then $\hat M_k=M_k+O_{L^q}(p^{-5/2})$ for $k=2,3$. Substituting $\hat M_2$, $\hat M_3$, $\hat M_4$, and $O_2^2$ for $M_2$, $M_3$, $M_4$, and $M_2^2$ in \eqref{eq_approximate_delta_4} gives an approximation to $\delta_i$ with error $O_{L^q}(p^{-5/2})$. \end{example}

A general version of approximating the conditional latent moments to order $N$ shall be implemented in detail in Step 3 below.  This calculation illustrates how successive moment recovery sharpens the initial approximations at a fixed order. At higher orders, recursive substitution generates products such as $M_2M_3$ and $M_2^3$, making coefficient growth and accumulated errors increasingly difficult to track directly. To obtain bounds valid under growing $N, p$ and multidimensional $d$,  we introduce the symbolic calculus.

\smallskip
\paragraph{Symbolic calculus.}
We introduce a symbolic calculus to represent combinations of moment products as formal expressions and provide algebraic rules to carry out high-order expansions concisely. 
Similar to \eqref{eq_Tdelta}, let $V=Z-\EE(Z\mid Y)$ and $\delta_i=\EE(\alpha_i^{\T}Z+\beta_i\mid Y)-T_i$
from \eqref{eq_Tdelta}. We suppress $r$ when working within one representation. For a multi-index $\gamma=(\gamma_1,\ldots,\gamma_d)\in\mathbb N_0^d$, write $|\gamma|=\sum_{\ell=1}^d\gamma_\ell$ and $V^\gamma=\prod_{\ell=1}^d V_\ell^{\gamma_\ell}$.
Since $\EE(V\mid Y)=0$, the conditional moments of interest begin at $|\gamma|=2$. We symbolize $\delta_i$, $V$, and $\EE[V^\gamma\mid Y]$, and $\delta_i$, which are treated as indeterminates to form a commutative polynomial ring, while $T_i$ are treated as coefficient.

The following two lemmas bound these quantities, which hold under either representation.
\begin{lemma}\label{lem_conditional_fluctuations}
Under Assumptions~\ref{ass_link} and~\ref{ass_model}, there is a constant $C>0$ such that, for $p\ge p_0$, $k\ge2$ and $q\ge1$,
\begin{equation}\label{eq_moments}
 \EE\norm V^k\le e^{Ck^2}p^{-k/2}\quad\text{ and }\quad \big\|\EE(\norm V^k\mid Y)\big\|_q \le e^{Cqk^2}p^{-k/2}.
\end{equation}
The bounds hold uniformly over $\mathcal S$ with $|\mathcal S|\le c_0p$ and both representations, where $c_0,p_0$ are as in Lemma~\ref{lem_design_alignment}.
\end{lemma}
\begin{proof}
    See Appendix~\ref{app_prove_lem_conditional_fluctuations}.
\end{proof}

\begin{lemma}\label{lem_bias_exponential}
Under the setup of Lemma~\ref{lem_conditional_fluctuations}, for every $q\ge1$ and $t\ge0$, uniformly over $i\in\cS$,
\begin{align*}
    &\norm{\delta_i}_q\le e^{Cq}p^{-1},\\
    &\EE e^{t\norm V}+\EE e^{t|T_i|}+\EE e^{t|\delta_i|} \le e^{C(1+t^2)}.
    \end{align*} \end{lemma}
    \begin{proof}
        See Appendix~\ref{app_prove_lem_bias_exponential}.
    \end{proof}

We now define the symbolic calculus. First introduce a formal vector $v=(v_1,\ldots,v_d)$ to represent $V$, a formal scalar $s$ to represent $\delta_i$, and moment symbols $m_\gamma$ to represent $\EE(V^\gamma\mid Y)$ for $|\gamma|\ge2$, respectively. The commutative polynomial ring considered is then generated by the indeterminates \(v_1,\ldots,v_d\), \(s\), and \(m_\gamma\) over coefficients consisting of real-valued functions of the predictors \(\mathbf T :=(T_i)_{i\in\mathcal S}\), often denoted by  $c(T_i)$. The coefficients are also allowed to depend on model parameters.


 For any $u$ representing such as symbol, we denote ${\rm ev}(u)$ the evaluated value (at the observables) of this quantity when substituting $v=V$, $ s=\delta_i$, and $ m_\gamma=\EE(V^\gamma\mid Y)$ into a formal expression. ${\rm ev}(u)$ is thus a random quantity and $u$ is a deterministic symbol. The preceding lemmas give, for $q\ge1$,
$\ell\in[d]$, $i\in\cS$, and $|\gamma|\ge2$,
\begin{align}\label{eq_moment_symbol_ctrl}
 \big\|\ev(v_\ell)\big\|_q \le C_qp^{-1/2}, \qquad \big\|\ev(s)\big\|_q \le C_qp^{-1},\quad\text{ and }\quad
 \big\|\ev(m_\gamma)\big\|_q  \le e^{Cq|\gamma|^2}p^{-|\gamma|/2}.
\end{align}
 We use symbolic calculus to organize the expansions and track their coefficients and weights before replacing the corresponding random quantities with the formal symbols. Treating these symbols as algebraically independent provides the convenience of doing high-order expansion on the conditional moments without requiring the random quantities behind the symbols to be probabilistically independent.

According to their scaling with $p^{-1/2}$, we assign weight via ${\rm wt}(\cdot)$ as
\begin{equation*}
 \hspace*{2em}{\rm wt}(v_\ell)=1, \hspace*{5em} {\rm wt}(s)=2,\qquad\text{ and }\qquad {\rm wt}(m_\gamma)=|\gamma|.\hspace*{5em}
\end{equation*}
These weights tracks powers of $p^{-1/2}$, which are added under multiplication of these quantities. They are different from the degree of the polynomial.
For example, $m_{2e_1}^2$ ($\EE(V_1^2\mid Y)^2$) and $m_{4e_1}$ ($\EE(V_1^4\mid Y)$) both have weight $4$.

We regard $v,s,m_\gamma$ as algebraically independent variables and work in the polynomial ring they generate, with coefficients depending on the common $T_i$ and fixed $\alpha_{0i}$ or $\alpha_{1i}$. Addition and subtraction combine coefficients of identical monomials, while multiplication expands products and adds their weights. Under $d=1$ for example, the moment symbols can be written as $m_2,m_3,m_4,\ldots$ For a fixed
coefficient $\xi$, one can perform the following calculations
\begin{align*}
 (m_2+\xi m_3) + \xi m_3&=m_2+2\xi m_3;\\
 (m_2+\xi m_3)-\xi m_3&=m_2;\\
 (m_2+\xi m_3)m_2 + m_4&=m_2^2+\xi m_2m_3 + m_4. 
\end{align*}
The distinct monomials $m_2^2$ and $m_4$ remain separate.

Next, we define formal expectation $\EE_{\mathrm f}$ as the conditional expectation given $Y$ under the symbolic form. It treats $s$ and the moment symbols as coefficients in polynomials, satisfying
\begin{equation*}
 \EE_{\mathrm f}(1)=1,\qquad \EE_{\mathrm f}(v_\ell)=0,\qquad \EE_{\mathrm f}(v^\gamma)=m_\gamma,\quad\text{and}\quad \EE_{\mathrm f}(s^a m^\eta v^\gamma) =s^a m^\eta m_{\gamma}.
\end{equation*}
This defines a linear projection on the polynomial ring under symbolic calculation, as $\EE_{\rm f}[\EE_{\rm f}(\cdot)] = \EE_{\rm f}(\cdot)$. Here, $m^{\eta} = \prod_{\rho:|\rho|\ge 2} m _{\rho}^{\eta_{\rho}}$ is some monomial with the same multi-index convention as $V^\gamma$. For any polynomial $P$, formal expectation and evaluation satisfy
\begin{equation*}
 \ev\{\EE_{\mathrm f}(P)\}
 =\EE\{\ev(P)\mid Y\}.
\end{equation*}Finally, for a polynomial $P$, denote by $[P]_N$ the truncation that preserves only monomials of weight at most $N$.  For example, $[(v_1+s)^2]_3=v_1^2+2v_1s$ as $s^2$ has weight four.
Till now, we have defined a polynomial ring in conditional-moment symbols equipped with multiplication, truncation, and formal expectation.

\begin{remark}
    The use of formal symbols and an evaluation map is related to the algebraic viewpoint of classical umbral calculus
\citep{rota1994classical}. Organizing and truncating terms by their asymptotic order also parallels cumulant-based higher-order expansions \citep{withers1984asymptotic}.
Here the weights are used to reflect scaling with $p^{-1/2}$, allowing conditional latent moments, correction terms, and their products to be handled within one polynomial calculus.

\end{remark}

\smallskip
\paragraph{Outline.} To bound $D_p(\theta_0,\theta_1)$, we refine 
\begin{equation}\label{eq_init_approx}\begin{aligned}
 \beta_i&=\EE(T_i+\delta_i)\approx \EE[T_i],\\
 \alpha_i^{\T}\alpha_j &={\rm Cov}(T_i+\delta_i,T_j+\delta_j) +\EE\{(\alpha_i^{\T}V)(\alpha_j^{\T}V)\}\approx {\rm Cov}(T_i,T_j).
\end{aligned}\end{equation}
using high-order polynomial approximations to the bias and covariance corrections. The argument has three components: controlling approximation errors under finite products; constructing symbolic expansions of $\delta_i$ and its conditional-mean equation; and expressing the corrections through conditional response moments. We proceed one by one.

\medskip\noindent{\underline{\it Step 1. Weighted polynomial calculus.}}
For $\mathbf T=(T_i)_{i\in\mathcal S}$ and a polynomial $P(m)=\sum_\eta a_\eta(\mathbf T)m^\eta$, define
\begin{equation*}
\|P(m)\|_{\mathrm{coeff}}=\sum_\eta\sup_{\mathbf T}|a_\eta(\mathbf T)|.
\end{equation*}
Here, $a_{\eta}(\mathbf{T})$ denotes the coefficient in $P(m)$ associated with monomial $m^{\eta}$, which is itself a function of $T_i$s, such as $f''(T_i)/f'(T_i)$.
We use the same definition for polynomials involving $v,s$. This norm controls the accumulation of coefficients when the approximation order grows.
Under this norm, we have for any two polynomials $P$ and $Q$ that $\|PQ\|_{\mathrm{coeff}} \le\|P\|_{\mathrm{coeff}}\|Q\|_{\mathrm{coeff}}$, $\|\EE_{\mathrm f}(P)\|_{\mathrm{coeff}} \le\|P\|_{\mathrm{coeff}}$, and
\begin{equation*}
 \|\ev\{P(m)\}\|_q \le\sum_\eta\sup_{\mathbf T}|a_\eta(\mathbf T)|\,\|\ev\{m^\eta\}\|_q \le \|P\|_{\mathrm{coeff}}\max_{\eta:a_{\eta}\neq 0}\|\ev\{m^\eta\}\|_q.
\end{equation*}
For the polynomials constructed below, we will show that they satisfy coefficient bounds of the form $\|P\|_{\mathrm{coeff}}\le e^{CN^2}$, with $C$ independent of $N$.

As a concrete example of using the coefficient bound, 
consider a monomial $P =\prod_{j=1}^r u_j$ where $u_j$ represents one of the symbols \(v_\ell,s,m_\gamma\) and has weights $\eta_j$, and let $\eta=\sum_{j=1}^r\eta_j\ge1$. By \eqref{eq_moment_symbol_ctrl}, the $L^a$ norm of each evaluated $u_j$ is bounded by $e^{C(1+a)\eta_j^2}p^{-\eta_j/2}$. By H\"older inequality with norm orders $q\eta/\eta_j$ on each {evaluated} component of weight {$\eta_j$}, the $L^q$ norm of this evaluated monomial $P$ can be bounded by
\begin{equation}\label{eq_weighted}
\|\ev(P)\|_q \le \exp\Big(C\sum_{j=1}^r\eta_j^2+Cq\eta\sum_{j=1}^r\eta_j\Big)p^{-\eta/2}\le e^{C\eta^2+Cq\eta^2}p^{-\eta/2} =e^{C(1+q)\eta^2}p^{-\eta/2}.
\end{equation}
Moreover, by Jensen's inequality, we know $\|\EE[\ev(P)\mid Y]\|_q \le e^{C(1+q)\eta^2}p^{-\eta/2}$. When $\|P\|_{\mathrm{coeff}}\le e^{CN^2}$ for some $C$, and  all weights of $P$ lie between $N+1$ and $cN$ with fixed $c$, we have
\begin{equation}\label{eq_discarded}
    \|\ev(P)\|_q\le\|P\|_{\mathrm{coeff}} e^{C'(1+q)N^2}p^{-(N+1)/2}\le  e^{C''(1+q)N^2}p^{-(N+1)/2}.
\end{equation}

We next establish that for a sequence of random quantities, if they admits proper polynomial approximation whose coefficient can be controlled in terms of $\|P\|_{\mathrm{coeff}}$, then the product of these random quantities can be well-approximated by some polynomial as well. In particular, we have

\begin{lemma}
\label{lem_weighted_calculus}
Let $x_1,\ldots,x_r$ be random quantities with positive integer weights $w_1,\ldots,w_r$ satisfying $\sum_{j=1}^r w_j\le N$. Suppose that, for every $q\ge1$, the $L^q$ norm satisfies $\|x_j\|_q\le e^{C(1+q)w_j^2}p^{-w_j/2}$, and that each $x_j$ admits a polynomial approximation $F_j(v,s,m)$ with weights between $w_j$ and $N$, such that 
\begin{equation}
    \|x_j-\ev(F_j(v,s,m))\|_q\le e^{C(1+q)N^2}p^{-(N+1)/2}\label{eq_lemma_weighted_cal_cond1}
\end{equation}
for any $q\ge 1$ and the coefficient of each weight $w\in\{w_j,w_j+1,\dots,N\}$ can be controlled by
\begin{equation}
 \|[F_j]_w-[F_j]_{w-1}\|_{\mathrm{coeff}}\le C^w,\label{eq_lemma_weighted_cal_cond2}
\end{equation}
Then there is a polynomial $P_r(v,s,m)$ such that, for every $q\ge1$,
\begin{equation}\label{eq_remainder}
 \big\|\prod_{j=1}^r x_j-\ev(P_r(v,s,m))\big\|_q \le e^{C'(1+q)N^2}p^{-(N+1)/2},
\end{equation}
for each $q\ge 1$. In particular, $P_r(v,s,m)$ can be defined inductively by setting $P_0(v,s,m) = 1$ and $P_j(v,s,m)=[P_{j-1}(v,s,m)F_j(v,s,m)]_N$ for $j=1,\ldots,r$.
Moreover, for an outer polynomial $R(x_1,\dots,x_r)$ of weight at most $N$ with $\|R\|_{\mathrm{coeff}}\le e^{CN^2}$, we have
\begin{equation}
    \|R(x_1,\ldots,x_r) - \ev\{[R(F_1,\ldots,F_r)]_N\}\|_q\le e^{C''(1+q)N^2}p^{-(N+1)/2}.\label{eq_remainder_outer}
\end{equation}
\end{lemma}
\begin{proof}
See Appendix~\ref{app_prove_lem_weighted_calculus}.
\end{proof}

\medskip\noindent{\underline{\it Step 2. Symbolic expansion of $\delta_i$}}
We now construct polynomial approximations to $\delta_i$. Note that the definition of $\delta_i = \EE(\alpha_i^{\T}Z + \beta_i\mid Y) - T_i$ implies \begin{equation*}
\EE\{f(T_i+\alpha_i^{\T}V+\delta_i)\mid Y\}= \EE(X_i\mid Y)=f(T_i).   
\end{equation*} Define $g_t(u) = \{f(t+u) - f(t)\}/f'(t)$ and we have the conditional mean equation
\begin{equation}
    \EE\big\{g_{T_i}(\alpha_i^{\T}V + \delta_i)\mid Y\big\} =0 .\label{eq_conditional_mean}
\end{equation}The denominator $f'(T_i)$ makes the first-order coefficient equal to one. Under nonlinear $f$, directly expanding \eqref{eq_conditional_mean} to solve $\delta_i$ requires repeatedly substituting the expansion of \(\delta_i\) into itself, generating increasingly complicated products of conditional moments.

We instead replace $V,\delta_i$ and conditional moments of $V$ by their formal symbols $v,s$ and $m$, respectively. Specifically, we can symbolize \eqref{eq_conditional_mean} into
\begin{equation}
     \mathcal F_i(s,m):=\EE_{\mathrm f}g_{T_i}(\alpha_i^{\T}v+s)=0, \qquad {\rm wt}(s)=2.\label{eq_conditional_mean_symbol}
\end{equation}The weight constraint ${\rm wt}(s) = 2$ reflects the $\|\delta_i\|_q=O(p^{-1})$ by Lemma~\ref{lem_bias_exponential}, since each unit of weight corresponds to a factor $p^{-1/2}$. 
Lemma~\ref{lem_symbolic_bias} below gives a unique formal solution $s=H_i(m)$ starting at weight two, whose evaluated truncation approximates $\delta_i$ with $L^q$ error at most $e^{C_qN^2}p^{-(N+1)/2}$. This expresses $\delta_i$ in terms of conditional latent moments. 

We note that the expansion of $\mathcal F_i(s,m)$ consists of the term $s$, which has weight two, and remaining terms of weight at least two. We can therefore determine $H_i(m)$ implicitly by matching terms of equal weight. For example, let $d=1$ and write $c_k=c_k(T_i)$ for brevity. Expanding \eqref{eq_conditional_mean_symbol} through weight four gives\begin{equation}\label{Eq_H_example1}
\begin{aligned}
\mathcal F_i(s,m) =\, &s+c_2\alpha_i^2m_2+c_3\alpha_i^3m_3+c_4\alpha_i^4m_4\\ &\, +c_2s^2+3c_3\alpha_i^2s m_2 +\{\text{terms of weight at least five}\}. \end{aligned}\end{equation}
Writing $H_i(m)=H_{i,2}(m)+H_{i,3}(m)+\cdots$, where $H_{i,r}(m)$ has weight $r$, and matching terms of equal weight in $\mathcal F_i(H_i(m),m)=0$, we obtain $H_{i,2}(m) = -c_2\alpha_i^2m_2$, $H_{i,3}(m) = -c_3\alpha_i^3m_3$, and $H_{i,4}(m) = -c_4\alpha_i^4m_4 + (3c_2c_3 - c_2^3)\alpha_i^4m_2^2$. The same procedure determines all higher-order components, although their expressions become increasingly complicated. Our argument does not give a closed-form expression for the full series.
The coefficient and remainder bounds needed for these constructions are collected in the following lemma.
\begin{lemma}\label{lem_symbolic_bias}
Under the current setup, for either representation there is a unique formal series $H_i(m)$ starting at weight two such that $\mathcal F_i(H_i(m),m)=0$. Define
\begin{equation*}
G_{i,N}(v,m)=[g_{T_i}(\alpha_i^{\T}v+H_i(m))]_N.
\end{equation*}
The weight-$w$ coefficients of $H_i(m)$, for $w\ge2$, and of $G_{i,N}(v,m)$, for $1\le w\le N$, have coefficient norm at most $C^w$. The same bounds hold for their first derivatives with respect to each loading coordinate $\alpha_{i\ell}$. The weight-one term of $G_{i,N}(v,m)$ is $\alpha_i^{\T}v$. For every $q\ge1$,
\begin{align*}
\|\delta_i-\ev\{[H_i(m)]_N\}\|_q
&\le e^{C(1+q)N^2}p^{-(N+1)/2},\\
\|g_{T_i}(\alpha_i^{\T}V+\delta_i)-\ev\{G_{i,N}(v,m)\}\|_q
&\le e^{C(1+q)N^2}p^{-(N+1)/2}.
\end{align*}
Moreover, $\|\delta_i\|_q\le e^{C(1+q)}p^{-1}$ and $\|g_{T_i}(\alpha_i^{\T}V+\delta_i)\|_q\le e^{C(1+q)}p^{-1/2}$.
\end{lemma}
\begin{proof}
See Appendix~\ref{app_prove_lem_symbolic_bias}.
\end{proof}

Lemma~\ref{lem_symbolic_bias} provides the approximations and bounds needed to apply Lemma~\ref{lem_weighted_calculus} to $\delta_i$ and $g_{T_i}(\alpha_i^{\T}V+\delta_i)$. In the next step, we use these results to construct response-based approximations to ${\rm ev}(m_\gamma)$. Substituting these approximations into $[H_i(m)]_N$ and removing the resulting higher-weight terms yields an accurate approximation to $\delta_i$.

\medskip\noindent{\underline{\it Step 3. Corrections via conditional response moments.}}

Recall that our goal is to obtain sharper corrections than \eqref{eq_init_approx}. Both approximation $\delta_i\approx {\rm ev}\{[H_i(m)]_N\}$  established in Lemma~\ref{lem_symbolic_bias} and the conditional covariance $\EE\{(\alpha_j^{\T}V)(\alpha_i^{\T}V)\}$ involve conditional latent moments $\EE[V^{\gamma}\mid Y] = {\rm ev}(m_{\gamma})$. In this step, we recover these corrections from conditional moments of distinct reference responses from $\cI$ as given in Lemma~\ref{lem_design_alignment}. To better motivate the derivation, we first give a cubic calculation as an illustration. 

Inspired by the approximations $\hat M_2$ and $\hat M_3$ constructed in \eqref{eq_construct_example}, we aim to construct observable quantities to recover the conditional latent moments through combinations of 
\begin{equation*}
    \zeta_{h}^{(s)} = \frac{X_{j_h^{(s)}}-\EE(X_{j_h^{(s)}}\mid Y)}{f'(T_{j_h^{(s)}})},\;\text{ for }\; h\in[d], j_h\in\cI_s,\text{ and }s\in [N].
\end{equation*}Here, we denote the selected reference group as $\cI_s = \{j^{(s)}_{h}:h\in[d]\}\subseteq[p]$. Without loss of generality, assume $j^{(s)}_1<j^{(s)}_2<\dots<j^{(s)}_d$ for every $s\in[N]$. We aim to construct a linear system from these $\zeta_{h}^{(s)}$ to form $O_{\gamma}$, such that it matches $\EE[V^{\gamma}\mid Y]$ to the first order.

To achieve this, we work on the symbolic form. First note that under model zero, 
\begin{equation}\label{eq_xi_gt}\EE[\zeta_h^{(s)}\mid Z_0,Y] = g_{T_{j_h^{(s)}}}\big(\alpha_{0,j_h^{(s)}}^{\T}V_0 +\delta_{0,j_h^{(s)}}\big),\end{equation}as $f(T_i) = \EE[X_i\mid Y]$ and $\alpha_i^{\T}Z+\beta_i=T_i+\alpha_i^{\T}V+\delta_i$ for $i = j_h^{(s)}$.
 By Lemma~\ref{lem_symbolic_bias}, this conditional mean is approximated by $\ev_0\{G_{j_h^{(s)},N}(v,m)\}$, whose leading term is $\alpha_{0,j_h^{(s)}}^{\T}V_0$.
Recall that Lemma~\ref{lem_design_alignment} ensures that $B_s:=(A_0)_{\mathcal I_s,}$ satisfies $\sigma_{\min}(B_s)\ge c_1$. The rows of $B_s$, ordered as $j_1^{(s)},\ldots,j_d^{(s)}$, give the coefficients of the linear forms $\alpha_{0,j_h^{(s)}}^{\T}v$.

But we want to construct $v^{\gamma}$ from $\alpha_{0,j_h^{(s)}}^{\T}v$. We use a tensor-product change of basis~\citep{kolda2009tensor}.
For any multi-index $\gamma\in\mathbb N_0^d$ with $2\le k=|\gamma|\le N$, choose and fix an ordered list $\ell_1,\ldots,\ell_k\in[d]$ such that $v^\gamma=\prod_{s=1}^k v_{\ell_s}$. In what follows, we suppress the dependence of these indices on $\gamma$ when it is clear from the context. With this notation, define $b_h^{(s)}=(B_s^{-1})_{\ell_s h}$. Then $\sum_{h=1}^d b_h^{(s)}\alpha_{0,j_h^{(s)}}^{\T}v= e_{\ell_s}^{\T}v = v_{\ell_s}$ and therefore 
 \begin{equation}\label{eq_tensor_multi_index_equality}
      \prod_{s=1}^k\Big\{\sum_{h=1}^d b_h^{(s)}\alpha_{0,j_h^{(s)}}^{\T}v\Big\}=v^\gamma.
 \end{equation}
 Thus, we construct 
\begin{equation}
    O_\gamma= \sum_{h_1,\ldots,h_k=1}^d \Big\{\prod_{s=1}^k(B_s^{-1})_{\ell_s h_s}\Big\} \mathbb E\Big\{\prod_{s=1}^k\zeta_{h_s}^{(s)}\mid Y\Big\}\, ,\quad \text{for}\quad 2\le |\gamma|\le N\label{eq_reference_observable}
\end{equation}
Here, $\ell_s$ are implicitly determined by $\gamma$. Letting $U_s=\sum_{h=1}^d b_h^{(s)}\zeta_h^{(s)}$, $O_{\gamma} =  \EE_0(\prod_{s=1}^k U_s\mid Y)
 =\EE_0\big\{\prod_{s=1}^k
       \EE_0(U_s\mid Z_0,Y)\mid Y\big\} = \EE_0\big[\prod_{s=1}^k\{\sum_{h=1}^d b_h^{(s)}\EE_0[\zeta_h^{(s)}\mid Z_0,Y]\}\mid Y\big]$.
As shown above, $\EE_0[\zeta_h^{(s)}\mid Z_0,Y]$ has conditional mean approximated by some polynomial with leading term $\alpha_{0,j_h^{(s)}}{}^\T V_0$. Then \eqref{eq_tensor_multi_index_equality} and linearity in $O_{\gamma}$ shows that $O_{\gamma}$ has leading term $\EE_0(V_0^\gamma\mid Y)$.

By \eqref{eq_xi_gt} and the independence between responses conditional on $(Z_r,Y)$, we know that 
\begin{equation}\label{eq_observable}
   \frac{\EE\{\prod_{i\in S}[X_i-\EE(X_i\mid Y)]\mid Y\}} {\prod_{i\in S}f'(T_i)} =\EE_r\Big[\prod_{i\in S}g_{T_i}(\alpha_{ri}^{\T}V_r+\delta_{ri})\mathrel{\Big|}Y\Big],
\end{equation}for $S \subseteq \cS$ consisting of excluded indices for $\cS$ given in Lemma~\ref{lem_design_alignment}. We first show that this product can be approximated by some polynomial, and then we can obtain the approximation for $O_{\gamma}$ by taking $S = \{j_{h_s}^{(s)}:s = 1,\dots,k\}$.

We first consider model zero and let $\xi_i^{(r)}=g_{T_i}(\alpha_{ri}^{\T}V_r + \delta_{ri})$ for $r\in\{0,1\}$. When $|S|\le N$,  by Lemma~\ref{lem_symbolic_bias} and the fact that $G_{i,N}(v,m)$ has weights between one and $N$, with linear part $\alpha_{0i}^{\T}v$, we can apply Lemma~\ref{lem_weighted_calculus} treating $\xi_i^{(0)}$ as the random quantities with weights $1$ and $G_{i,N}(v,m)$ as the polynomial approximation to get
\begin{equation}\label{eq_reference_product_error}
 \Big\|\prod_{i\in S}\xi_i^{(0)}-\ev_0\{P_{S}(v,m)\}\Big\|_q \le e^{C(1+q)N^2}p^{-(N+1)/2},
\end{equation} where $S = \{i_1,\dots,i_{|S|}\}$ and $ P_{S}(v,m) := P_{S,|S|}(v,m)$ is inductively defined by $P_{S,0}(v,m)=1$ and $P_{S,t}(v,m)=[P_{S,t-1}(v,m)G_{i_t,N}(v,m)]_N$ for $t=1,2,\dots,|S|$. Thus $P_S(v,m) = [\prod_{i\in S}G_{i,N}(v,m)]_N$. It remains to pass the bound from $\prod_{i\in S}\xi_i^{(0)}$ to the conditional expectation as in \eqref{eq_observable}. Let $Q_S{(m)} = [\EE_{\rm f}\{P_S(v,m)\}]_N = \EE_{\rm f}\{P_S(v,m)\}$. Because $\EE_0[\ev_0\{P_{S}(v,m)\}\mid Y] = \ev_0([\EE_{\rm f}\{\prod_{i\in S}G_{i,N}(v,m)\}]_N)$ as formal expectation preserves the weights in each monomial, we have by Jensen inequality that 
\begin{align*}
     \Big\|\EE_0\big(\prod_{i\in S}\xi_i^{(0)}\mid Y\big) -\ev_0(Q_S(m))\Big\|_q = &\,\Big\|\EE_0\big\{\prod_{i\in S}\xi_i^{(0)} -\ev_0(P_S(v,m))\mid Y\big\}\Big\|_q\\
 \le& \, \Big\|\prod_{i\in S}\xi_i^{(0)}-\ev_0(P_S(v,m))\Big\|_q
 \le e^{C(1+q)N^2}p^{-(N+1)/2}.
\end{align*}

Therefore, for $O_{\gamma}$ given in \eqref{eq_reference_observable}, we can find a polynomial approximation given by 
\begin{equation}\label{eq_forward_polynomial}
 F_\gamma(m)= \Bigg[ \sum_{h_1,\ldots,h_k=1}^d \Big\{\prod_{s=1}^k(B_s^{-1})_{\ell_s h_s}\Big\} \EE_{\mathrm f}\Big\{\prod_{s=1}^kG_{j^{(s)}_{h_s},N}(v,m)\Big\} \Bigg]_N,
\end{equation}
such that $\varepsilon_\gamma=O_\gamma-\ev_0\{F_\gamma(m)\}$ satisfies
\begin{align}
     \|\varepsilon_\gamma\|_q &\le\sum_{h_1,\dots,h_k=1}^d\Big|\prod_{s=1}^k(B_s^{-1})_{\ell_s h_s}\Big|
 \left\|\EE_0\left(\prod_{s=1}^k\xi_{j_{h_s}^{(s)}}^{(0)} \mathrel{\Big|}Y\right) -\ev_0\Big(\Big[\EE_{\mathrm f}\Big\{\prod_{s=1}^kG_{j^{(s)}_{h_s},N}(v,m)\Big\}\Big]_N\Big)\right\|_q \nonumber\\
 &\le C_B^k e^{C(1+q)N^2}p^{-(N+1)/2} \le e^{C'(1+q)N^2}p^{-(N+1)/2}.\label{eq_forward}
\end{align} Here, we use $\prod_{s=1}^k\sum_{h=1}^d|(B_s^{-1})_{\ell_s h}| \le(d^{1/2}/c_1)^k\le C_B^k$ with $C_B := \max\{1,{d}^{1/2}/c_1\}$, by $\|B_s^{-1}\|\le c_1^{-1}$ from Lemma~\ref{lem_design_alignment}.  Moreover, by the previous argument to construct $O_{\gamma}$, it is not hard to verify that $[F_{\gamma}(m)]_k = m_{\gamma}$. Apply a similar argument in showing \eqref{eq_lemma_weighted_cal_cond2}, we can obtain 
\begin{equation}\label{eq_lemma_weighted_cal_cond_F_gamma}\|[F_\gamma]_w-[F_\gamma]_{w-1}\|_{\mathrm{coeff}}\le C^w,
\end{equation}
for $|\gamma|\le w\le N$, and $ \|\partial_{\alpha_{jl}}[F_\gamma]_w-\partial_{\alpha_{jl}}[F_\gamma]_{w-1}\|_{\mathrm{coeff}}\le C^w$. Also, for $\|O_{\gamma}\|_{q}$, apply 
\begin{equation*}
\prod_{s=1}^k\sum_{h=1}^d|(B_s^{-1})_{\ell_s h}|\le C_B^k,\qquad \|\xi_i^{(r)}\|_a\le e^{C(1+a)}p^{-1/2},
\end{equation*}
and Jensen inequality to get
\begin{equation}\begin{aligned}
    \|O_\gamma\|_q\le &\, C_B^k \sum_{h_1,\dots, h_k=1}^d\Big\|\EE_r\big(\prod_{s=1}^k\xi_{j_{h_s}^{(s)}}^{(r)}\mid Y\big)\Big\|_q\\ \le &\,C_B^k \sum_{h_1,\dots, h_k=1}^d\prod_{s=1}^k\|\xi_{j_{h_s}^{(s)}}^{(r)}\|_{qk}\\\le &\, e^{C'(1+q)k^2}p^{-k/2}.
\end{aligned}\label{eq_O_gamma_q}\end{equation}

Till now, we establish the polynomial approximations $F_{\gamma}(m)$ of order $N$ to all constructed observed quantities $O_{\gamma}$ for $2\le |\gamma|\le N$, which have their leading terms as $\ev(m_{\gamma})$. Next, we want to use these observed quantities to approximate $\delta_{0i}$ and $\EE_0\{(\alpha_{0i}{}^\T V)(\alpha_{0j}{}^\T V)\mid Y\} $. Recall that we have shown $\delta_i\approx \ev\{[H_i(m)]_N\}$ in Lemma~\ref{lem_symbolic_bias}. We aim to find some polynomial $R_i(O) \approx \delta_{0i}$. One may want to directly use $\ev_0\{[H_i(F(m))]_N\}$. While this can keep the leading terms aligned, interactions among the higher-order terms in $F(m)$ produce a more complicated expansion of $H_i(F(m))$. As in the successive recovery of $M_4,M_3,M_2$ in~\ref{eq_construct_example}, these terms must be removed before the desired accuracy is attained.

To cancel the higher order terms introduced in $F_{\gamma}$, we define an operator $\cU$. Let $\mathcal P_N$ be the space of polynomials in these moment symbols whose monomials have weight at most $N$, with coefficient functions of $\mathbf T = (T_i)_{i\in\cS}$ as in the coefficient norm defined above. For any $P(m) = \sum_{\eta}a_{\eta}(\mathbf T)m^{\eta}\in\mathcal{P}_N$, define 
\begin{equation}\label{eq_substitution_operator}
    \mathcal U P(m) :=\Big[\sum_\eta a_\eta(\mathbf T) \prod_{\gamma: 2\le |\gamma|\le N}F_\gamma(m)^{\eta_\gamma}\Big]_N.
\end{equation}
By definition, we have $(\mathcal{U} - \cI)^{N+1} = 0$ where $\cI$ is identity mapping here. Therefore we know the inverse of $\mathcal{U}$ can be expressed as $\cU^{-1}=\sum_{l=0}^N\{-(\cU-\cI)\}^l$, which will be used to cancel the high order terms introduced by polynomial approximation $F_{\gamma}$. We then introduce two technical lemmas.
\begin{lemma}\label{lem_design_alignment_sub}
Under Assumption~\ref{ass_model}, there is a constant $C>0$ such that,
for $p\ge p_0$,
\begin{equation}\label{eq_alignment}
 \min_{O^{\T}O=I_d}\max_{j\in[p]} \norm{\alpha_{1j}-O\alpha_{0j}} \le C\max_{i,j\in[p]} |\alpha_{1i}^{\T}\alpha_{1j} -\alpha_{0i}^{\T}\alpha_{0j}|.
\end{equation}
\end{lemma}
\begin{proof}
    See Appendix~\ref{app_prove_lem_design_alignment_sub}.
\end{proof}
The next lemma shows that replacing $m_{\gamma}$ with the corresponding $O_{\gamma}$ in $\mathcal{U}^{-1}[H_i(m)]_N$ can well approximate $\delta_{0i}$. Write $D_p = D_p(\theta_0,\theta_1)$ for short.

\begin{lemma}\label{lem_observable_corrections}
Suppose Assumptions~\ref{ass_link} and~\ref{ass_model} hold. Under model zero, define polynomials
\begin{equation}\label{eq_corrections}
 R_i(m)=\mathcal U^{-1}[H_i(m)]_N,\qquad
 R_{ij}(m)=\mathcal U^{-1}\EE_{\mathrm f} \{(\alpha_{0i}^{\T}v)(\alpha_{0j}^{\T}v)\}.
\end{equation}
They contain only terms of weights between two and $N$, and satisfy
\[
 \|R_i(m)\|_{\mathrm{coeff}}+\|R_{ij}(m)\|_{\mathrm{coeff}}\le e^{CN^2}.
\]
Let $R_i(O)$ and $R_{ij}(O)$ be evaluated value obtained by replacing each moment symbol $m_\gamma$ with $O_\gamma$. By Lemma~\ref{lem_design_alignment_sub}, we may assume without loss of
generality that the two representations have been orthogonally aligned so that $\max_{a\in[p]}\norm{\alpha_{1a}-\alpha_{0a}}\le CD_p$.
For $q\in\{1,2\}$ and $r\in\{0,1\}$, there is a common $\tau>0$ such that, whenever $4\le N\le\tau\log p$ and $Nd+2\le c_0p$, with $c_0$ as in Lemma~\ref{lem_design_alignment}, the following bound holds:
\begin{equation}\label{eq_common_corrections}
 \|R_i(O)-\delta_{ri}\|_q +\big\|R_{ij}(O)-\EE_r\{(\alpha_{ri}^{\T}V_r) (\alpha_{rj}^{\T}V_r)\mid Y\}\big\|_q\le rC_qD_p/p+e^{C(1+q)N^2}p^{-(N+1)/2}.
\end{equation}
Here, the first term on the right sides says that, for $r=0$, the $D_p/p$ term vanishes. Under this case, the inequality holds for any $q\ge 1$.
\end{lemma}
\begin{proof}
See  Appendix~\ref{app_prove_lem_observable_corrections}.

\end{proof}

\medskip\noindent{\underline{\it Step 4. Combining the bounds to control $D_p$.}}

Finally, we can compare
\begin{align*}
 \beta_{ri}&=\EE\{\EE_r(\alpha_{ri}^{\T}Z_r+\beta_{ri}\mid Y)\} = \EE[T_i+\delta_{ri}],\\
 \alpha_{ri}^{\T}\alpha_{rj} &={\rm Cov}\{\EE_r(\alpha_{ri}^{\T}Z_r+\beta_{ri}\mid Y), \EE_r(\alpha_{rj}^{\T}Z_r+\beta_{rj}\mid Y)\} +\EE_r\{(\alpha_{ri}^{\T}V_r)(\alpha_{rj}^{\T}V_r)\}.
\end{align*}
By Lemma~\ref{lem_observable_corrections}, we have 
\begin{align*}
    |\beta_{0i}-\beta_{1i}| \le &\, \big\|T_i+\delta_{0i} - (T_i+\delta_{1i})\big\|_1\\\le &\, \|\delta_{0i}-R_i(O)\|_1 + \|\delta_{1i}-R_i(O)\|_1\\\le&\, CD_p/p +e^{CN^2}p^{-(N+1)/2}.
\end{align*}
For $|\alpha_{0i}^{\T}\alpha_{0j}-\alpha_{1i}^{\T}\alpha_{1j}|$, we first note that by Jensen's inequality, $\|\EE_r(\alpha_{ri}^{\T}Z_r+\beta_{ri}\mid Y)\|_2 \le\|\alpha_{ri}^{\T}Z_r+\beta_{ri}\|_2 \le M(1+ d^{1/2})$ for $r\in\{0,1\}$. Here, $\|\cdot\|_2$ denotes the \(L^2\) norm for random variables. Then writing $\mu_{ri}:=\EE_r(\alpha_{ri}^{\T}Z_r+\beta_{ri}\mid Y) =T_i+\delta_{ri}$ for $r\in\{0,1\}$, we have $\|\mu_{ri}\|_2\le M(1+d^{1/2})$. With $\mu_{0i}-\mu_{1i}=\delta_{0i}-\delta_{1i}$ and Lemma~\ref{lem_observable_corrections}, we know 
\begin{equation*}
    \|\mu_{0i}-\mu_{1i}\|_2 \le\|\delta_{0i}-R_i(O)\|_2+\|\delta_{1i}-R_i(O)\|_2 \le CD_p/p+e^{CN^2}p^{-(N+1)/2}.
\end{equation*}
Finally, using  $|{\rm Cov}(U,V)|\le\|U\|_2\|V\|_2$,
one can compute
\begin{align*}
 |\alpha_{0i}^{\T}\alpha_{0j}-\alpha_{1i}^{\T}\alpha_{1j}|
 &\le|{\rm Cov}(\mu_{0i}-\mu_{1i},\mu_{0j})| +|{\rm Cov}(\mu_{1i},\mu_{0j}-\mu_{1j})|\\
 &\quad+\sum_{r=0}^1\Big\|\EE_r\{(\alpha_{ri}^{\T}V_r)
                (\alpha_{rj}^{\T}V_r)\mid Y\}-R_{ij}(O)\Big\|_1\\
 &\le M(1+ d^{1/2})
       \big\{\|\mu_{0i}-\mu_{1i}\|_2+\|\mu_{0j}-\mu_{1j}\|_2\big\}\\
 &\quad+\sum_{r=0}^1\Big\|\EE_r\{(\alpha_{ri}^{\T}V_r)
                (\alpha_{rj}^{\T}V_r)\mid Y\}-R_{ij}(O)\Big\|_1\\
 &\le CD_p/p+e^{CN^2}p^{-(N+1)/2},
\end{align*} where the last inequality uses \eqref{eq_common_corrections} at $q=1,2$.

Taking the two maxima defining $D_p$ and increasing $C$ gives
\begin{equation}\label{eq_final}
D_p\le CD_p/p+e^{CN^2}p^{-(N+1)/2},
\end{equation} with $4\le N\le\tau\log p$ and $Nd+2\le c_0p$ required by Lemma~\ref{lem_design_alignment}. When $p\ge 2C$, choose $N=\lfloor\tau\log p\rfloor$ and sufficiently small $\tau$, we know $D_p\le2e^{-(\tau/4)(\log p)^2}$. The proof is then completed by Lemma~\ref{lem_design_alignment_sub}.

\subsection{Proof of Corollary~\ref{cor_estimation}}
\label{app_estimation}

\begin{proof}
Both $(\theta_0,G_0)$ and $(\theta^{\dagger},G^{\dagger})$ belong to $\mathcal M_p$ and yield the same marginal distributions for the observations . Hence Theorem~\ref{thm_asymptotic} gives, for all sufficiently large $p$, that $D_p(\theta^{\dagger},\theta_0)
 \le C\exp\{-c(\log p)^2\}$.
The triangle inequality holds for $D_p(\cdot,\cdot)$ by definition, and therefore
\[
 \big|D_p(\hat\theta_n,\theta_0)  -D_p(\hat\theta_n,\theta^{\dagger})\big|  \le D_p(\theta^{\dagger},\theta_0),
\]
which proves \eqref{eq_estimation_bound}.

For the inference part, the Lipschitz condition and Theorem~\ref{thm_asymptotic} give $|\psi(\theta^{\dagger})-\psi(\theta_0)| \le CL_\psi\exp\{-c(\log p)^2\}$.
After division by $\hat s_{n,p}$, this difference is $o_{\mathbb P}(1)$, so the $\hat s_{n,p}^{\, -1}\{\psi(\hat\theta_n)-\psi(\theta^{\dagger})\}$ and $\hat s_{n,p}^{\, -1}\{\psi(\hat\theta_n)-\psi(\theta_0)\}$  have the same limiting distribution by Slutsky's theorem.
\end{proof}

\subsection{Proof of Corollary~\ref{prop_prediction}}\label{app_prediction}

For notational convenience we write $(\theta_1,G_1)$ as $(\theta^{\dagger},G^{\dagger})$.  By Lemma~\ref{lem_design_alignment_sub}, we can choose an orthogonal matrix $Q$ such that $\max_{j\in[p]}\|Q\alpha_{1j}-\alpha_{0j}\|\le CD_p(\theta_0,\theta_1)$. Applying $Z_1\mapsto QZ_1$ replaces the loading vectors by $Q\alpha_{1j}$ and the latent distribution by $Q_\#G^{\dagger}$. For simplicity, we set the aligned loading vectors and latent distribution again by $\alpha_{1j}$ and $G_1$, respectively, and thus $\max_{j\in[p]}\|\alpha_{1j}-\alpha_{0j}\|\le CD_p(\theta_0,\theta_1)$. In the proof of Lemma~\ref{lem_conditional_fluctuations}, we construct $\hat Z$ to approximate $Z$. Specifically, suppose $\hat Z_r(x)$ minimizes
\begin{equation*}
    \Phi_r(z;x)=\frac1p\sum_{j=1}^p \Big\{\int_0^{\alpha_{rj}^\T z + \beta_{rj}}f(s)ds -x_j(\alpha_{rj}^{\T}z+\beta_{rj})\Big\} +\frac{\|z\|^2}{2p}.
\end{equation*}
Each $p\Phi(z;x)$ is strongly convex, so its minimizer exists and is unique.
In the proof of Lemma~\ref{lem_conditional_fluctuations}, \eqref{eq_fluctuations_reconstruction_moments}  implies $\EE_{\theta_r,G_r}\|Z - \hat Z_r(X)\|^2\le Cp^{-1}$ for $r\in\{0,1\}$ here $X$ is the observed vector. To bound $\|\hat Z_0(X) - \hat Z_1(X)\|$, we note that because $\Phi_0(\cdot;x)$ is strongly convex and $\partial_{Z}\Phi_1(\hat Z_1(x);x) = 0$ for any $x\in\{0,1\}^p$, it holds that
\begin{equation*}
    \|\hat Z_0(x)-\hat Z_1(x)\| \le p\|\partial_{Z}\Phi_0(\hat Z_1(x);x)-\partial_{Z}\Phi_1(\hat Z_1(x);x)\|.
\end{equation*}
Since $\int_{-\infty}^{+\infty}f'(s)ds = 1$, by \eqref{eq_linkslope}, we know $1\ge f'(t)\int_{-\infty}^{+\infty}e^{-2K|s-t|}ds = f'(t)/K$, and thus $f$ is $K$-Lipschitz. Then by Assumption~\ref{ass_model}\ref{ass_model_1}, one can obtain 
\begin{align*}
    \|\nabla\Phi_0(z;x)-\nabla\Phi_1(z;x)\| \le &\, p^{-1}\Big\|\sum_{j=1}^p\big\{ \alpha_{0j}f(\alpha_{0j}^{\T} z + \beta_{0j}) - \alpha_{1j}f(\alpha_{1j}^{\T} z+\beta_{1j}) - (\alpha_{0j} - \alpha_{1j})x_j\big\}\Big\|\\\le  &\, CD_p(\theta_0,\theta_1)(1+  \norm z).
\end{align*}
For $\hat Z_0(X)$ and $\hat Z_1(X)$, we know by \eqref{Eq_bound_Z_norm} that $\|\hat Z_0(X)\|\le pM$ and $\|\hat Z_1(X)\|\le pM$. Combine these estimates and we know $\|\hat Z_0(X) - \hat Z_1(X)\|\le Cp^2D_p(\theta_0,\theta_1)$.

Now, we are ready to bound the Wasserstein distance between $G_0$ and $G_1$. Let $H_r$ be the distribution of $\hat Z_r(X)$. The joint distribution of $(Z,\hat z_r(X))$ couple $G_r$ and $H_r$, while $(\hat Z_0(X),\hat Z_1(X))$ couples $H_0$ and $H_1$. Then we have
\begin{align*}
 W_1(G_0,G_1)
 &\le W_1(G_0,H_0)+W_1(H_0,H_1)+W_1(H_1,G_1)\\ &\le\sum_{r=0}^1\EE_{\theta_r,G_r}\|Z-\hat Z_r(X)\| +\EE_P\|\hat Z_0(X)-\hat Z_1(X)\|\\
 &\le Cp^{-1/2}+Cp^2D_p(\theta_0,\theta_1) \le C'p^{-1/2}
\end{align*}
for all sufficiently large $p$, by Cauchy--Schwarz and Theorem~\ref{thm_asymptotic}. Then because orthogonal transformation preserve $W_1(\cdot,\cdot)$, we conclude that $|W_1(Q_\#\hat G_n,G_0)
             -W_1(\hat G_n,G^{\dagger})|\le W_1(G_0,Q_\#G^{\dagger})\le Cp^{-1/2}$.

Finally, for any $L_h$-Lipschitz function $h$ and distributions $\mu,\nu$, integrating $|h(z)-h(z')|\le L_h\|z-z'\|$ under any coupling and taking the infimum gives $|\EE_\mu h-\EE_\nu h|\le L_hW_1(\mu,\nu)$.
With $W_1(Q_\#\hat G_n,G_0) \le W_1(\hat G_n,G^{\dagger})+Cp^{-1/2}$, apply this inequality with $\mu=Q_\#\hat G_n$ and $\nu=G_0$ proves \eqref{eq_prediction_bound}.

\subsection{Details for Example~\ref{example_bifactor}}\label{app_sec_bif}

We give sufficient regularity conditions for all $p$:
\begin{enumerate}[label=$(\roman*)$]
    \item $d\ge2$ is fixed, and both $\theta,\tilde\theta$ have the same known group membership $g(j)\in\{1,\ldots,d\}$;
    \item $\Sigma^{-1/2}Z$ and $\tilde\Sigma^{-1/2}\tilde Z$ satisfy Assumptions~\ref{ass_link} and~\ref{ass_model} in dimension $d+1$;
    \item $\lambda_{\min}(\Sigma)$ and $\lambda_{\min}(\tilde\Sigma)$ are bounded away from zero;
    \item for each $g\in\{1,\ldots,d\}$, there are distinct $a_g,b_g$ with $g(a_g)=g(b_g)=g$ such that the $2\times 2$ matrices
\[
L_g=(\lambda_{jk})_{j\in\{a_g,b_g\},\,k\in\{0,1\}},
\quad\text{and}\quad
\tilde L_g=(\tilde\lambda_{jk})_{j\in\{a_g,b_g\},\,k\in\{0,1\}}
\]
satisfy $\sigma_{\min}(L_g)\ge c$ and $\sigma_{\min}(\tilde L_g)\ge c$ for some constant $c>0$.
\end{enumerate}

Let $a_j=\Sigma^{1/2}\lambda_j$ and $u_g=\Sigma^{1/2}e_{g+1}^{(d+1)}$, with analogous tilded notation, where $e_i^{(d+1)}$ denotes the $i$th standard basis vector in $\RR^{d+1}$. By the covariance structure, we have $a_j=\lambda_{j0}e_1^{(d+1)}+\lambda_{j1}u_{g(j)}$, where $\|u_g\|=1$ and $u_g^{\T}e_1^{(d+1)}=0$. Theorem~\ref{thm_asymptotic} and Lemma~\ref{lem_design_alignment_sub} yield an orthogonal $Q$ such that
\begin{equation}\label{eq_bifactor}
\|\beta-\tilde\beta\|_\infty
+\max_{j\in[p]}\|a_j-Q\tilde a_j\|\le C\exp\{-c(\log p)^2\}.
\end{equation}
Fix a group $g$ and let $(t_1,t_2)^\T=(\tilde L_g^\T)^{-1}(1,0)^\T$. Condition $(iv)$ gives $|t_1|+|t_2|\le C$. Since $\tilde a_j=\tilde\lambda_{j0}e_1^{(d+1)}+\tilde\lambda_{j1}\tilde u_g$ for items in group $g$, we have $e_1^{(d+1)}=t_1\tilde a_{a_g}+t_2\tilde a_{b_g}$. Therefore, we have
\begin{equation*}
\big\|Qe_1^{(d+1)}-(t_1a_{a_g}+t_2a_{b_g})\big\|
\le C\exp\{-c(\log p)^2\}.
\end{equation*}
So $Qe_1^{(d+1)}$ is within $C\exp\{-c(\log p)^2\}$ of $\operatorname{span}(e_1^{(d+1)},u_g)$ as $a_{a_g}$ and $a_{b_g}$ belongs to this space.

Condition $(iii)$ gives $1-|u_1^\T u_2|\ge c>0$, where $u_1=\Sigma^{1/2}e_2^{(d+1)}$ and $u_2=\Sigma^{1/2}e_3^{(d+1)}$. Since $u_1,u_2$ are unit vectors orthogonal to $e_1^{(d+1)}$ according to the structure of $\Sigma$, the largest eigenvalue of $u_1u_1^\T+u_2u_2^\T$ on this orthogonal complement is $1+|u_1^\T u_2|$. Using $\|Qe_1^{(d+1)}\|=1$ and the preceding bounds for the two planes, we obtain
\begin{align*}
\min_{s_0\in\{-1,1\}}\|Qe_1^{(d+1)}-s_0e_1^{(d+1)}\|^2
&=2\big[1-\big|(e_1^{(d+1)})^\T Qe_1^{(d+1)}\big|\big]\\
&\le2\big[1-\{(e_1^{(d+1)})^\T Qe_1^{(d+1)}\}^2\big]\\
&\le\frac{2}{1-|u_1^\T u_2|}\sum_{g=1}^2 \big[1-\{(e_1^{(d+1)})^\T Qe_1^{(d+1)}\}^2 -\{u_g^\T Qe_1^{(d+1)}\}^2\big]\\
&\le C\exp\{-2c(\log p)^2\}.
\end{align*}
Here each summand is the squared distance from $Qe_1^{(d+1)}$ to $\operatorname{span}(e_1^{(d+1)},u_g)$. Taking square roots gives $\|Qe_1^{(d+1)}-s_0e_1^{(d+1)}\|\le C\exp\{-c(\log p)^2\}$ for some $s_0\in\{-1,1\}$. With $\|a_j - Q\tilde a_j\|\le C\exp\{-c(\log p)^2\}$, projecting the aligned loadings onto the general direction and its orthogonal complement therefore gives, uniformly in $j$, we have
\begin{equation*}
|\lambda_{j0}-s_0\tilde\lambda_{j0}|
+\|\lambda_{j1}u_{g(j)}-\tilde\lambda_{j1}Q\tilde u_{g(j)}\| \le C\exp\{-c(\log p)^2\}.
\end{equation*}
By condition $(iv)$, each group has an anchor whose group loading is bounded away from zero. Choosing a sign $s_g$ using this anchor gives $\|u_g-s_gQ\tilde u_g\|\le C\exp\{-c(\log p)^2\}$ and $|\lambda_{j1}-s_g\tilde\lambda_{j1}|\le C\exp\{-c(\log p)^2\}$ for every item in that group. With $D=\operatorname{diag}(s_0,s_1,\ldots,s_d)$, these direction bounds give $\|Q\tilde\Sigma^{1/2}-\Sigma^{1/2}D\|\le C\exp\{-c(\log p)^2\}$, which yields $\|\tilde\Sigma-D\Sigma D\|\le C\exp\{-c(\log p)^2\}$. By condition $(iii)$, taking positive definite square roots preserves this order, so $\|\tilde\Sigma^{1/2}-D\Sigma^{1/2}D\|\le C\exp\{-c(\log p)^2\}$. Hence $\|(Q-D)\tilde\Sigma^{1/2}\|\le C\exp\{-c(\log p)^2\}$, and the uniform bound on $\|\tilde\Sigma^{-1/2}\|$ gives
\[
\|Q-D\|\le C\exp\{-c(\log p)^2\}.
\]
Combining these bounds with \eqref{eq_bifactor} proves the displayed bound in Example~\ref{example_bifactor}.



\section{Proofs of Technical Lemmas}\label{app_technical_lemmas}
\renewcommand{\theequation}{\thesection.\arabic{equation}}
\setcounter{equation}{0}

Throughout this appendix, we use the notation and conventions of
Appendix~\ref{app_asymptotic} .

\subsection{Proof of Lemma~\ref{lem_design_alignment}}\label{app_prove_lem_design_alignment}

By Assumption~\ref{ass_model}\ref{ass_model_2}, there is some $p_0$ such that for any $p\ge p_0$, $\sum_{j=1}^p\alpha_{rj}\alpha_{rj}^{\T} \succeq\kappa pI_d$ for $ r\in\{0,1\}$. Choose $c_0>0$ sufficiently small that $c_0M^2\le\kappa/2$. So for any $\mathcal J\subseteq[p]$ with $|\mathcal J|\le c_0p$, Assumption~\ref{ass_model}\ref{ass_model_1} gives
\begin{equation*}
 0\preceq\sum_{j\in\mathcal J}\alpha_{rj}\alpha_{rj}^{\T} \preceq |\mathcal J|M^2 I_d\preceq(\kappa p/2)I_d.
\end{equation*}
This proves the first part. 

For the reference sets, we show that they can be selected successively. Suppose that $\mathcal I_1,\ldots,\mathcal I_{s-1}$ have already been chosen and let $ \mathcal R_s=[p]\setminus ( \cup_{t=1}^{s-1}\mathcal I_t\cup\mathcal T)$.  Since $|\mathcal T|+(s-1)d\le Nd+2\le c_0p$, Cauchy--Binet formula imply
\begin{equation*}
 \sum_{\substack{I\subseteq\mathcal R_s\\|I|=d}} \det\{(A_0)_{I,}\}^2
 =\det\Big(\sum_{j\in\mathcal R_s} \alpha_{0j}\alpha_{0j}^{\T}\Big) \ge (\kappa p/2)^d.
\end{equation*}
There are at most $\binom pd\le p^d/d!$ terms in this sum. Thus there exists some $\mathcal I_s\subseteq\mathcal R_s$ with $|\mathcal I_s|=d$ satisfies $\big|\det\{(A_0)_{\cI_s}\}\big|\ge {(d!)}^{1/2}(\kappa/2)^{d/2}$. Since $\|B_s\|\le \|B_s\|_{\rm F}\le {d}^{1/2}{M}$ for $B_s = (A_0)_{\cI_s}$, this gives $\sigma_{\min}(B_s)\ge |\det(B_s)|/\|B_s\|^{d-1}\ge (d!)^{1/2}(\kappa/2)^{d/2}({d}^{1/2}M)^{1-d} := c_1$. In this way, we successively select the reference set with the required properties.

\subsection{Proof of Lemma~\ref{lem_conditional_fluctuations}}\label{app_prove_lem_conditional_fluctuations}

Fix either representation and suppress the index. We follow Proposition B.4(a) in~\citet{ostrovskii2021finite} to construct a regularized convex M-estimator to approximate $Z$ from $Y = X_{\cS^c}$. Specifically, let
\begin{equation*}
    \hat Z = \mathop{\arg\min}_{z\in \RR^d} \frac1p\sum_{j\in [p]\setminus\cS} \Big\{\int_0^{\alpha_j^{\T}z+\beta_j}f(s)ds -X_j(\alpha_j^{\T}z+\beta_j)\Big\} +\frac{\norm z^2}{2p}.
\end{equation*}
Then we establish the required bound by controlling $Z - \hat Z$ and $\EE(Z - \hat Z\mid Y)$.
First, it is not hard to check that this yields a unique minimizer $\hat Z$. The scores of the estimating equation at $Z$ and $\hat Z$ are given by
\begin{equation}\label{eq_define_aux_score}
    S = \frac1p\sum_{j\in [p]\setminus\cS} \alpha_j\{X_j-f(\alpha_j^{\T}Z+\beta_j)\}-\frac Zp\end{equation} with $\hat S = p^{-1}\sum_{j\in [p]\setminus\cS} \alpha_j\{X_j-f(\alpha_j^{\T}\hat Z+\beta_j)\}-p^{-1}{ \hat Z} = 0$.
Let $h = \hat Z - Z$, and we have
\begin{equation}\label{eq_fluctuations_score_identity}
    h^{\T} S = h^{\T} (S - \hat S) = h^{\T}\frac1p\sum_{j\in[p]\setminus\cS}\alpha_j\big\{f(\alpha_j^{\T}(Z+h)+\beta_j) -f(\alpha_j^{\T}Z+\beta_j)\big\} +\frac{\norm h^2}{p} 
\end{equation}When $h\neq 0$, by integral mean value theorem, the first term on the right side can be written as $(p\|h\|)^{-1}\sum_{j\in[p]\setminus\cS}(\alpha_j^\T h)^2\int_0^{\|h\|} f'\{\alpha_j^{\T}(Z+sh/\|h\|)+\beta_j\}\,ds$. On $0\le s\le 1\wedge\|h\|$,  by Assumption~\ref{ass_model}\ref{ass_model_1} and \eqref{eq_linkslope}, we know $f'(\cdot)\ge ce^{-C\|Z\|}$. With $\sum_{j\in[p]\setminus \cS}\alpha_j\alpha_j^{\T} \succeq (\kappa/2)pI_d$ by Lemma~\ref{lem_design_alignment}, we know this part of integral is lower bounded by 
\begin{equation*}
    ce^{-C\|Z\|}(\kappa/2)\|h\|\big(1\wedge\|h\|\big).
\end{equation*}
Thus on the event $\mathcal{E} := \{\|S\|\le ce^{-C\|Z\|}/2\}$, we know from \eqref{eq_fluctuations_score_identity} that $\|S\|\ge ce^{-C\|Z\|}(1\wedge\|h\|)$, and therefore $\norm h\le1/2$ and
\begin{equation}\label{eq_fluctuations_good_score}
 \norm h\le C e^{C\norm Z}\norm S.
\end{equation}
Apparently, \eqref{eq_fluctuations_good_score} holds when $h=0$. On the other side, on the complement event $\mathcal{E}^c = \{\|S\|> ce^{-C\|Z\|}/2\}$, we know by $\hat S =0 $ that 
\begin{equation}\norm{\hat Z} = \norm{\sum_{j\in [p]\setminus\cS} \alpha_j\{X_j-f(\alpha_j^{\T}\hat Z+\beta_j)\}}\le pM\label{Eq_bound_Z_norm}\end{equation} by Assumption~\ref{ass_model}\ref{ass_model_1}, which yields
\begin{equation}\label{eq_fluctuations_global_bound}
 \norm h\le pM+\norm Z.
\end{equation}

Next, because each summand in $\sum_{j\in [p]\setminus\cS} \alpha_j\{X_j-f(\alpha_j^{\T}Z+\beta_j)\}$ is conditionally independent of each other with mean zero, for each component of $\bar S = S +Z /p$, denoted by $\bar S_k$ for $k\in[d]$, one has by Hoeffding's inequality that
\begin{equation*}
    \EE\{e^{t\bar S_k}\mid Z=z\} \le\exp\Big\{\frac{t^2}{8p^2} \sum_{j\in[p]\setminus\cS}\alpha_{jk}^2\Big\}
 \le e^{M^2t^2/(8p)}.
\end{equation*}
Thus for \eqref{eq_define_aux_score}, Markov inequality yields $\EE(\norm S^k\mid Z=z) \le C^kk^{k/2}\{p^{-k/2}+p^{-k}\norm z^k\}$ for any $k\ge 2$. Then \eqref{eq_fluctuations_good_score} gives 
\begin{equation*}
    \EE(\norm h^k\mathbf 1_{\mathcal E}\mid Z=z) \le  e^{Ck\norm z} C^kk^{k/2}\{p^{-k/2}+p^{-k}\norm z^k\}.
\end{equation*}
For the complementary event, Markov's inequality at order $3k$ gives $\PP(\mathcal E^c\mid Z=z) \le C^{3k}({3k})^{3k/2}e^{Ck\norm z} \{p^{-3k/2}+p^{-3k}\norm z^{3k}\}$
Multiplying this bound by $(pM+\norm z)^k$ from \eqref{eq_fluctuations_global_bound} bounds its contribution to the $k$th moment. In particular, since $p\ge1$ and $pM+\norm z\le Cp(1+\norm z)$, the two event bounds together imply
\begin{equation*}
 \EE(\norm h^k\mid Z=z) \le \{C{(3k)}^{1/2}\}^{3k}p^{-k/2} e^{Ck\norm z}(1+\norm z)^{4k}.
\end{equation*}
Here, the powers of $p$ from the complementary event are $p^{-k/2}$ and $p^{-2k}$, both bounded by $p^{-k/2}$. Note that $\EE\|Z\|^k$ and $\EE e^{Ck\|Z\|}$ can be bounded by $c^{C'k^2}$ by the sub-Gaussianity of $Z$. Therefore, with $\{C({3k})^{1/2}\}^{3k}\le e^{C'k^2}$, and $(1+x)^{4k}\le e^{4kx}$, $e^{Ck\|z\|}\le e^{\nu \|z\|^2}e^{C'k^2}$, we arrive at
\begin{equation}\label{eq_fluctuations_reconstruction_moments}
 \EE\norm{Z-\hat Z}^k\le e^{Ck^2}p^{-k/2}.
\end{equation}

Finally, we note that $V=(Z-\hat Z)-\EE(Z-\hat Z\mid Y)$ as $\hat Z$ is measurable with respect to $Y$.
The triangle inequality and conditional Jensen's inequality give
$\norm V_k\le2\norm{Z-\hat Z}_k$, proving the first part in \eqref{eq_moments}. Applying conditional Jensen's inequality once more, with moment order $kq\ge2$, yields
\begin{equation*}
 \big\|\EE(\norm V^k\mid Y)\big\|_q \le\big(\EE\norm V^{kq}\big)^{1/q} \le e^{Cqk^2}p^{-k/2}.
\end{equation*}
which completes the proof.

\subsection{Proof of Lemma~\ref{lem_bias_exponential}}\label{app_prove_lem_bias_exponential}
Fix either representation and suppress the index. For any $i\in\cS$, by conditional independence, $f(T_i) = \EE\{f(\alpha_i^\T Z + \beta_i)\mid Y\}$. We first bound $T_i$. By \eqref{eq_linkslope}, the function
$\varphi(x)=f(H^{-1}\log x)$ is concave in $x>0$. By Jensen's inequality, we know $f(T_i)=\EE\{\varphi(e^{H(\alpha_i^\T Z + \beta_i)})\mid Y\} \le f[H^{-1}\log\EE\{e^{H(\alpha_i^\T Z + \beta_i)}\mid Y\}\big]$. Because $f$ is strictly increasing, we know 
\begin{equation*}
    T_i\le H^{-1}\log\EE(e^{H(\alpha_i^\T Z + \beta_i)}\mid Y).
\end{equation*}
A lower bound can be similarly obtained by considering function $1-f(-t)$:
\begin{equation*}
    T_i\ge -H^{-1}\log\EE(e^{-H(\alpha_i^\T Z + \beta_i)}\mid Y).
\end{equation*}
With $\delta_ i = \EE[\alpha_i^\T Z + \beta_i\mid Y] - T_i = H^{-1}\log \{e^{H \EE[\alpha_i^\T Z + \beta_i\mid Y] }\} - T_i$ and $V = Z- \EE(Z\mid Y)$, we know 
\begin{equation*}
     -H^{-1}\log\EE(e^{H(\alpha_i^\T V)}\mid Y) \le\delta_i\le H^{-1}\log\EE(e^{-H(\alpha_i^\T V)}\mid Y).
\end{equation*}
By $e^x\le1+x+x^2e^{|x|}/2$, $\log x\le x-1$ and $\EE[V\mid Y] = 0$, we have 
\begin{equation}\label{eq_bias_conditional_quadratic}
 |\delta_i| \le\frac H2\EE\{(\alpha_i^\T V)^2e^{H|\alpha_i^\T V|}\mid Y\}
 \le\frac{HM^2}{2}\EE\{\norm V^2e^{HM\norm V}\mid Y\}.
\end{equation}
For $e^{HM\norm V}$, by Cauchy's inequality and Jensen's inequality, one has 
\begin{align*}
 \EE e^{t\norm V}
 \le\EE\left[e^{t\norm Z}e^{t\norm{\EE(Z\mid Y)}}\right] \le(\EE e^{2t\norm Z})^{1/2} (\EE e^{2t\norm{\EE(Z\mid Y)}})^{1/2} \le\EE e^{2t\norm Z} \le L e^{t^2/\nu},
\end{align*}where the last inequality follows from $2t\norm Z\le\nu\norm Z^2+t^2/\nu$ and the sub-Gaussianity of $Z$. Then with $(\EE\norm{V}^{4q})^{1/(2q)}\le e^{Cq}p^{-1}$ by Lemma~\ref{lem_conditional_fluctuations}, we conclude with $\|\EE(\cdot\mid Y)\|_{q}\le \|\cdot\|_q$ and H\"older's inequality that
\begin{equation*}
    \norm{\delta_i}_q\le \frac{HM^2}{2} (\EE\norm V^{4q})^{1/(2q)} (\EE e^{2qHM\norm V})^{1/(2q)} \le e^{Cq}p^{-1}.
\end{equation*}
This proves the first part of Lemma~\ref{lem_bias_exponential}.

For the second part, it suffices to show that both $\EE( e^{t|\alpha_i^\T Z + \beta_i|})$ and $\EE[e^{t|T_i|}]$ can be bounded by $e^{C(1+t^2)}$, as Cauchy's inequality and Jensen's inequality  yield
\begin{equation*}
    \EE e^{t|\delta_i|}\le (\EE e^{2t|\EE[\alpha_i^\T Z + \beta_i\mid Y]|})^{1/2} (\EE e^{2t|T_i|})^{1/2} \le(\EE e^{2t|\alpha_i^\T Z + \beta_i|})^{1/2} (\EE e^{2t|T_i|})^{1/2} .
\end{equation*}
For $\EE( e^{t|\alpha_i^\T Z + \beta_i|})$, by Assumptions~\ref{ass_model}, it is easy to check that it can be bounded by $e^{C(1+t^2)}$. For $\EE[e^{t|T_i|}]$, when $t\ge H$, we note by Jensen's inequality that, for the convex map $x\mapsto x^{H^{-1}t}$, there is 
\begin{align*}\EE e^{t|T_i|}\le &\, \EE e^{tT_i}+\EE e^{-tT_i}\le\EE\{\EE(e^{H(\alpha_i^{\T}Z + \beta_i)}\mid Y)\}^{t/H}+\EE\{\EE(e^{-H(\alpha_i^{\T}Z + \beta_i)}\mid Y)\}^{t/H}\\
\le&\, \EE e^{t(\alpha_i^{\T}Z + \beta_i)}+\EE e^{-t(\alpha_i^{\T}Z + \beta_i)} \le e^{C(1+t^2)}.\end{align*} When $t<H$, $\EE[e^{t|T_i|}]\le \EE[e^{H|T_i|}\le c^{C(1+H^2)}$.
Therefore we complete the second part of Lemma~\ref{lem_bias_exponential}.

\subsection{Proof of Lemma~\ref{lem_weighted_calculus}}\label{app_prove_lem_weighted_calculus}

\begin{proof}
    Since the polynomials $F_j$ and $P_j$ share the formal arguments $(v,s,m)$, we suppress in this proof for simplicity. It suffices to prove \eqref{eq_remainder} for $P_r$ defined by setting $P_0 = 1$ and $P_j=[P_{j-1}F_j]_N$ for $j=1,\ldots,r$. By this definition, we have 
    \begin{align} \prod_{j=1}^r x_j-\ev(P_r) =&\, \sum_{j=1}^r\{\ev(P_{j-1})x_j-\ev(P_j)\}\prod_{l>j}x_l\nonumber\\ = &\,  \sum_{j=1}^r\Big[\underbrace{\ev(P_{j-1})\{x_j-\ev(F_j)\}\prod_{l>j}x_l}_{\gamma_{j1}} + \underbrace{\ev\big\{P_{j-1}F_j-[P_{j-1}F_j]_N\big\}\prod_{l>j}x_l}_{\gamma_{j2}}\Big].\label{eq_product}
\end{align}
To bound $\gamma_{j1}$, we first control $\ev(P_{j-1})$ by showing $\|P_{j-1}\|_{\mathrm{coeff}}\le e^{CN}$. Note that $P_j$ is defined inductively by multiplication and truncation at $N$, and thus for any $j\le w\le N$, we have
\begin{equation*}
    \big\|[P_j]_w-[P_j]_{w-1}\big\|_{\mathrm{coeff}}\le \sum_{\substack{a_1+\cdots+a_j=w\\a_1,\ldots,a_j\ge1}}\prod_{i=1}^j\big\|[F_i]_{a_i}-[F_i]_{a_i-1}\big\|_{\mathrm{coeff}} \le C^w\binom{w-1}{j-1}\le(2C)^w,
\end{equation*}by \eqref{eq_lemma_weighted_cal_cond2}. Then $\|P_{j-1}\|_{\mathrm{coeff}} \le \sum_{w=j-1}^N(2C)^w\le e^{CN}$ and $\|F_j\|_{\mathrm{coeff}} \le \sum_{w=w_j}^NC^w\le e^{CN}$. With \eqref{eq_weighted}, we know $\|\ev(P_{j-1})\|_a \le\|P_{j-1}\|_{\mathrm{coeff}}e^{C(1+a)N^2} \le e^{C(1+a)N^2}$. For $\prod_{l>j}x_l$, we know by $\|x_j\|_q\le e^{C(1+q)w_j^2}p^{-w_j/2}$ and H\"older inequality that for any $q\ge 1$,
\begin{align*}
    \Big\|\prod_{l>j}x_l\Big\|_q \le  \prod_{l>j}\|x_l\|_{q(\sum_{k>j}w_k)/w_l} 
 \le   e^{C(1+q)(\sum_{l>j}w_l)^2} p^{-\frac12\sum_{l>j}w_l} \le  e^{C(1+q)N^2}.
\end{align*}
Therefore, for $\gamma_{j1}$, we have by \eqref{eq_lemma_weighted_cal_cond1} that
\begin{align*}
 \|\gamma_{j1}\|_q \le  &\, \|\ev(P_{j-1})\|_{3q}\, \|x_j-\ev(F_j)\|_{3q}\, \Big\|\prod_{l>j}x_l\Big\|_{3q}\\
 \le &\,  e^{C(1+3q)N^2} \big\{e^{C(1+3q)N^2}p^{-(N+1)/2}\big\} e^{C(1+3q)N^2}\\ \le &\,  e^{C'(1+q)N^2}p^{-(N+1)/2},
\end{align*}

Next, we bound  $\gamma_{j2}$. Because $\|P_{j-1}F_j-[P_{j-1}F_j]_N\|_{\mathrm{coeff}}\le \|P_{j-1}F_j\|_{\mathrm{coeff}}\le \|P_{j-1}\|_{\mathrm{coeff}}\|F_j\|_{\mathrm{coeff}}\le e^{CN}$, we can apply \eqref{eq_discarded} to get $\|\ev\{P_{j-1}F_j-[P_{j-1}F_j]_N\}\|_{2q}\le e^{C(1+2q)N^2}p^{-(N+1)/2}$. Together with $\big\|\prod_{l>j}x_l\big\|_{2q} \le e^{C(1+2q)N^2}$ and the Cauchy inequality, we arrive at
\begin{equation*}
    \|\gamma_{j2}\|_q \le \|\ev\{P_{j-1}F_j-[P_{j-1}F_j]_N\}\|_{2q}\big\|\prod_{l>j}x_l\big\|_{2q }\le e^{C(1+q)N^2}p^{-(N+1)/2}.
\end{equation*}
Combine the upper bounds for $\gamma_{j1}$ and $\gamma_{j2}$ in \eqref{eq_product}, and using $r\le N$, yields \eqref{eq_remainder}. Similar calculation follows for outer polynomial $R(x_1,\dots,x_r)$ to show \eqref{eq_remainder_outer}.

\end{proof}

\subsection{Proof of Lemma~\ref{lem_symbolic_bias}}\label{app_prove_lem_symbolic_bias}
\begin{proof}
Fix either representation and suppress its index. Table~\ref{tab_step2_expansions} summarizes the formal expansions used throughout this proof.
Write $c_k(t)=f^{(k)}(t)/\{k!f'(t)\}$. The expressions $\mathcal F_i(s,m)$ and $H_i(m)$ are formal power series; their finite truncations are polynomials.

\begin{table}[htbp]
\centering
\small
\setlength{\tabcolsep}{4pt}
\renewcommand{\arraystretch}{1.2}
\caption{ Formal expansions used in the proof of Lemma~\ref{lem_symbolic_bias}.}
\label{tab_step2_expansions}
\begin{tabular}{@{}>{\raggedright\arraybackslash}p{0.18\textwidth}>{\raggedright\arraybackslash}p{0.80\textwidth}@{}}
\toprule
Expression & Definition \\
\midrule
$\mathcal F_i(s,m)$ & $\EE_{\mathrm f}\{g_{T_i}(\alpha_i^{\T}v+s)\}$ \\[6pt]
$P_{i,N}(v,s)$ & $\displaystyle\sum_{k=1}^N c_k(T_i)(\alpha_i^{\T}v+s)^k$ as Taylor expansion of $g_{T_i}(\alpha_i^{\T}v+s)$ to order $N$.\\[6pt]
$H_i(m)$ & Solution to $\mathcal F_i(H_i(m),m)=0$. \\[6pt]
$G_{i,N}(v,m)$ & $\{P_{i,N}(v,H_i(m))\}_N$\\
\bottomrule
\end{tabular}
\end{table}

The series $\mathcal F_i(s,m)$ provides approximation to the conditional-mean equation used to determine $H_i(m)$ so that  $G_{i,N}(v,m)$ replaces $s$ with $H(m)$ in $P_{i,N}(v,s)$. Both $P_{i,N}(v,s)$ and $G_{i,N}(v,m)$ approximate $g_{T_i}(\alpha_i^{\T}V+\delta_i)$ in the evaluated value. 

 We first bound the coefficients obtained by recursively solving $H_i(m)$ from \eqref{eq_conditional_mean_symbol}. Write $b_w = \|[H_i(m)]_w-[H_i(m)]_{w-1}\|_{\mathrm{ coeff}}$. For $g_t(u)$, by Assumption~\ref{ass_link}\ref{ass_link_2}, we know $|c_k(t)|\le K^{k-1}$ and we have by the integral form of remainder that
\begin{equation}\label{eq_taylor}
\begin{aligned}
    \Big|g_t(u) -\sum_{k=1}^N c_k(t)u^k\Big|\le &\,\frac{|u|^{N+1}}{N!f'(t)} \int_0^1(1-s)^N { |f^{(N+1)}(t+s u)|}\,ds\\ 
    \overset{(1)}{\le} &\, \frac{|u|^{N+1}}{N!} \int_0^1|1-s|^N \frac{(N+1)!K^Nf'(t+{s} u)}{f'(t)}\,ds\\
    \overset{(2)}{\le} &\,  |u|^{N+1}K^N e^{H|u|}
    \end{aligned}
\end{equation}
Here, $(1)$ follows from Assumption~\ref{ass_link}\ref{ass_link_2}; $(2)$ follows from  \eqref{eq_linkslope}. Since $c_1 = 1$, we can write \eqref{eq_conditional_mean_symbol} with $s = H_i(m)$ as 
\begin{equation}\label{eq_him_scalar_est_eq}
    H_i(m) = - \sum_{k\ge2}c_k(T_i)\, \EE_{\mathrm f}\big[\big\{\alpha_i^{\T} v + H_i(m)\big\}^k\big].
\end{equation}

To bound the coefficients of $H_i(m)$, we use a majorant argument that records coefficient norms by weight. Introduce a scalar variable $t$, and let $[t^w]$ extract the coefficient of $t^w$ in a scalar series. A moment symbol $m_\gamma$ is recorded by $t^{|\gamma|}$, and the coefficient at $t^w$ is denoted by $[t^w]t^{|\gamma|}$, which equals one when $w=|\gamma|$ and zero otherwise. For $H_i(m)$, the corresponding scalar series is $\sum_{\ell\ge2} \|[H_i]_\ell-[H_i]_{\ell-1}\|_{\mathrm{coeff}}\,t^\ell =\sum_{\ell\ge2}b_\ell t^\ell$, so $b_w=[t^w]\sum_{\ell\ge2}b_\ell t^\ell$.
Applying the coefficient norm inequalities to \eqref{eq_him_scalar_est_eq} yields recursive bounds for $b_w$.
The majorant argument then constructs a scalar series $h(t)=\sum_{w\ge2}h_wt^w$ with nonnegative coefficients satisfying $b_w\le h_w$. Showing that $h$ is analytic near the origin allows us to apply Cauchy's coefficient estimate~\citep[Chapter~2, Corollary~4.3]{stein2003complex}, giving $b_w\le h_w\le C^w$.

In particular, equating the weight-$w$ components on both sides of \eqref{eq_him_scalar_est_eq}, we obtain
\begin{equation*}
    b_w\le[t^w] \sum_{k= 2}^w\sup_{T_i}|c_k(T_i)| \Big\{\|\alpha_i^{\T}v\|_{\mathrm{coeff}} t + \sum_{\ell=2}^{w-1}b_{\ell}t^{\ell}\Big\}^k
\end{equation*}
By Assumption~\ref{ass_model}, we have $\|\alpha_i^{\T}v\|_{\mathrm{coeff}} =\sum_{\ell=1}^d|\alpha_{i\ell}|\le C_0$.
Together with $\sup_{T_i}|c_k(T_i)|\le K^{k-1}$, we have for a suitable constant $C\ge K(C_0\vee 1)$ that
\begin{align*}
    b_w 
    \le [t^w]\sum_{k=2}^w C^k \Big(t+\sum_{\ell=2}^{w-1}b_\ell t^\ell\Big)^k.
\end{align*}
To derive an explicit bound for $b_w$, we introduce a series $h(t) = \sum_{w\ge 2}h_wt^w $ satisfying 
\begin{equation*}
    h(t) = \sum_{k\ge 2}C^k\{t + h(t)\}^k = \frac{C^2\{t+h(t)\}^2}{1-C\{t+h(t)\}} := H(t,h).
\end{equation*}Since $H(0,0) = 0$ and $\partial_hH(t,h)|_{(t,h) = (0,0)} = 0$,  implicit function theorem implies that this equation admits some analytic solution with $h(0) = h'(0) = 0$. 
As in the recursive construction, we can recursively determine $h_w>0$ from the above equation. Moreover, if $b_{l}\le h_l$ for $2\le l<w$, we know \begin{equation*}
    b_w\le[t^w]\sum_{k=2}^w C^k
       \left(t+\sum_{\ell=2}^{w-1}h_\ell t^\ell\right)^k =h_w.
\end{equation*}
It is easy to check that $b_2\le h_2$. Thus inductively, we know $b_w\le h_w$ for every $w\ge2$. By Cauchy's estimate for analytic functions, we have
\begin{equation*}
 b_w = \|[H_i]_w-[H_i]_{w-1}\|_{\mathrm{coeff}} \le h_w \le r^{-w}\max_{|t|=r}|h(t)|\le C_1^w,
\end{equation*}

Next, we bound the coefficient $\|\partial_{\alpha_{i\ell}}H_{i,w}(m)\|_{\rm coeff}$, where $H_{i,w}$ is treated as a polynomial in parameter $\alpha_{i\ell}$. To bound the first-order derivatives, we allow these coordinates to be complex with $\max_\ell|\alpha_{i\ell}|\le M+1$ for $M$ given in Assumption~\ref{ass_model}\ref{ass_model_1}.
The preceding calculation still applies after increasing $C_0$ as $\sum_\ell|\alpha_{i\ell}|\le d(M+1)$. Then Cauchy's formula implies that, at any real $\alpha_i$, $\partial_{\alpha_{i\ell}}H_{i,w}(m;\alpha_i) = (2\pi {\mathrm i})^{-1}\oint_{|z|=1}z^{-2}{H_{i,w}(m;\alpha_i+ze_\ell)}dz$, and therefore
\begin{equation}\label{eq_cauchy_der_estimate}
    \|\partial_{\alpha_{i\ell}}H_{i,w}\|_{\mathrm{coeff}}\le\sup_{|z|=1} \|H_{i,w}(m;\alpha_i+ze_\ell)\|_{\mathrm{coeff}} \le C_2^w.
\end{equation}

Now we move forward to bound $\delta_i - {\ev }\big\{[H_i(m)]_N\big\}$. Factorize $\mathcal{F}_i(s,m)$ to get $\mathcal{F}_i(s,m) = (s-H_i(m))\mathcal Q_i(s,m)$ for  a formal series $\mathcal Q_i(s,m)$ given by 
\begin{equation}\label{Eq_define_series_cQ}
    \mathcal{Q}_i(s,m) = 1+\sum_{k\ge2}c_k(T_i)\EE_{\mathrm f}\Big\{ \sum_{\ell=0}^{k-1}(\alpha_i^{\T}v+s)^{k-1-\ell}(\alpha_i^{\T}v+H_i(m))^\ell\Big\}
\end{equation}
Writing $\mathcal{Q}_{i,N}^{-1}(s,m)=[\mathcal Q_i(s,m)^{-1}]_N$, we have
\[
[s-H_i(m)]_N=\big[[\mathcal F_i(s,m)]_N \mathcal{Q}_{i,N}^{-1}(s,m)\big]_N.
\]
In what follows, we omit the common arguments of polynomials when they are clear from the context. It suffices to control the evaluated right side.  
To control this term, following a similar strategy in the proof of Lemma~\ref{lem_weighted_calculus}, we first control $\|[\mathcal Q_i^{-1}]_w-[\mathcal Q_i^{-1}]_{w-1}\|_{\mathrm{coeff}}$. 
Apply $[t^w]$ to both side of the equation, and note that the contribution of $c_k(T_i)$, $\alpha_i^{\T}v$, $s$, and $H_i(m)$ are bounded by $K^{k-1}$, $C_0t$, $t^2$, and $h(t)$, respectively, which gives
\begin{equation*}
    \|[\mathcal Q_i-1]_w-[\mathcal Q_i-1]_{w-1}\|_{\mathrm{coeff}} \le[t^w]\sum_{k\ge2}kC^k\{t + t^2 + h(t)\}^{k-1},
\end{equation*}Let 
$q(t) = \sum_{k\ge 2}kC^k\{t + t^2 + h(t)\}^{k-1}$, which has nonnegative coefficients. Then 
$\|[(1-\mathcal Q_i)^j]_w- [(1-\mathcal Q_i)^j]_{w-1}\|_{\mathrm{coeff}}\le [t^w] q(t)^j$.
With $q(0) = 0$, we can select some radii $\rho$ such that $q(\rho)<1/2$ and the inverse $\mathcal Q_i^{-1}=\sum_{j\ge0}(1-\mathcal Q_i)^j$ satisfies
\begin{equation*}
    \big\|[\mathcal Q_i^{-1}]_w-[\mathcal Q_i^{-1}]_{w-1}\big\|_{\mathrm{coeff}} \le[t^w]\sum_{j=0}^{\infty}q(t)^j = [t^w]\frac{1}{1-q(t)} 
\end{equation*}
Suppose $(1-q(t))^{-1} = \sum_{n=0}^{\infty}d_nt^n$. Since $d_n\ge 0$ as $q(t)$ has nonnegative coefficients, we know at $\rho$, there is $d_w\rho^w\le \sum_{n=0}^{\infty}d_n\rho^n = (1-q(\rho))^{-1}$. Thus $[t^w](1-q(t))^{-1} = d_w$ and 
\begin{equation*}
    \big\|[\mathcal Q_i^{-1}]_w-[\mathcal Q_i^{-1}]_{w-1}\big\|_{\mathrm{coeff}}\le\frac{\rho^{-w}}{1-q(\rho)} \le2\rho^{-w}\le C^{w+1}.
\end{equation*}

We now bound $\delta_i-\ev\{[H_i(m)]_N\}$ using the preceding truncated identity. Note that 
\begin{equation}
\|\delta_i-\ev\{[H_i]_N\}\|_q \le\|\ev\{[\mathcal F_i]_N\}\|_{2q} \|\ev(\mathcal{Q}_{i,N}^{-1})\|_{2q}+\big\|\ev\big\{[\mathcal F_i]_N \mathcal{Q}_{i,N}^{-1} -\big[[\mathcal F_i]_N \mathcal{Q}_{i,N}^{-1}\big]_N\big\}\big\|_q.\label{Eq_bound_delta_poly}
\end{equation}

For $\ev\{[\mathcal F_i]_N\}$, let $P_{i,N}(v,s)=\sum_{k=1}^Nc_k(T_i)(\alpha_i^{\T}v+s)^k$. Then $[\mathcal F_i]_N=\EE_{\mathrm f}\{[P_{i,N}]_N\}$. It is easy to check by previous bounds that $\|P_{i,N}\|_{\mathrm{coeff}}\le\sum_{k=1}^NK^{k-1}(C_0+1)^k\le e^{CN}$ and $\|[\mathcal F_i]_N\|_{\mathrm{coeff}}\le\|P_{i,N}\|_{\mathrm{coeff}}\le e^{CN}$.
Let $u = \alpha_i^{\T}V + \delta_i$. With $\|u\|_r\le e^{C(1+r)}p^{-1/2}$ and $\EE e^{a|u|}\le e^{C(1+a^2)}$, we have by \eqref{eq_taylor} and H\"older's inequality, that
\begin{equation}
\begin{aligned}
\|\ev(P_{i,N})-g_{T_i}(u)\|_q
&\le K^N
   (\EE e^{2qH|u|})^{1/(2q)}
   \|u\|_{2q(N+1)}^{N+1}\\
&\le K^N e^{C(1+q)}
   e^{C(N+1)\{1+2q(N+1)\}}
   p^{-(N+1)/2}\\
&\le e^{C'(1+q)N^2}p^{-(N+1)/2}.
\end{aligned}\label{eq_taylor_remiander}
\end{equation}
We then invoke $\EE[g_{T_i}(u)\mid Y] = 0$ as in \eqref{eq_conditional_mean} to get
\begin{equation}\label{eq_discard_coro}\begin{aligned}
    \big\|\ev\{[\mathcal F_i]_N\}\big\|_q
 =&\,\Big\|\EE\big[\ev\{[P_{i,N}]_N- P_{i,N}\} + \ev(P_{i,N})-g_{T_i}(u) \mid Y\big]\Big\|_q\\\le&\,\|\ev\{P_{i,N}-[P_{i,N}]_N\}\|_q + \|\ev(P_{i,N})-g_{T_i}(u)\|_q \\
 \le&\, e^{C(1+q)N^2}p^{-(N+1)/2},
\end{aligned}\end{equation}where the last inequality is due to \eqref{eq_discarded} and \eqref{eq_taylor_remiander}.

For ev$(\mathcal{Q}_{i,N}^{-1})$, we note that $\|\mathcal{Q}_{i,N}^{-1}\|_{\mathrm{coeff}}\le1+\sum_{w=1}^NC^{w+1}\le e^{CN}$. By \eqref{eq_weighted}, we have that
\begin{equation*}
 \|\ev(\mathcal{Q}_{i,N}^{-1})\|_{2q} \le\|\mathcal{Q}_{i,N}^{-1}\|_{\mathrm{coeff}}e^{C(1+2q)N^2} \le e^{C(1+q)N^2}.
\end{equation*}

Finally, for the last term in \eqref{Eq_bound_delta_poly}, since the terms discarded from $[\mathcal F_i]_N \mathcal{Q}_{i,N}^{-1}$ have weights between $N+1$ and $2N$, we arrive at 
\begin{equation*}
    \big\|[\mathcal F_i]_N \mathcal{Q}_{i,N}^{-1}-\big[[\mathcal F_i]_N \mathcal{Q}_{i,N}^{-1}\big]_N\big\|_{\mathrm{coeff}}\le\|[\mathcal F_i]_N\|_{\mathrm{coeff}}
       \|\mathcal{Q}_{i,N}^{-1}\|_{\mathrm{coeff}}\le e^{CN}.
\end{equation*}
Combining these estimates into \eqref{Eq_bound_delta_poly} to yield 
\begin{equation}
    \|\delta_i-\ev\{[H_i]_N\}\|_q\le e^{C(1+q)N^2}p^{-(N+1)/2}.\label{eq_Herror}
\end{equation}

We next construct a polynomial approximation and establish an $L^q$ bound for $g_{T_i}(\alpha_i^{\T}V+\delta_i)$ similarly as for $\delta_i$.
Specifically, let  $G_{i,N}(v,m)=[g_{T_i}(\alpha_i^{\T}v+H_i(m))]_N$. 
Its weight-one term is $\alpha_i^{\T}v$. The other coefficient can be bounded in the Taylor series \eqref{eq_taylor} by 
\begin{equation}\label{eq_factorerror_coediff}
    \|[G_{i,N}(v,m)]_w-[G_{i,N}(v,m)]_{w-1}\|_{\mathrm{coeff}} \le [t^w]\sum_{k=1}^{\infty}K^{k-1}\{C_0t+h(t)\}^k=[t^w]\frac{C_0t+h(t)}{1-K\{C_0t+h(t)\}} \le C^w
\end{equation}
The same Cauchy derivative estimate as in \eqref{eq_cauchy_der_estimate} gives
\begin{equation*}
\max_{\ell\in[d]}\|\partial_{\alpha_{i\ell}}([G_{i,N}(v,m)]_w-[G_{i,N}(v,m)]_{w-1})\|_{\mathrm{coeff}}\le C^w.
\end{equation*}

To transfer the approximation of $\delta_i$ to $G_{i,N}(v,m)$, let
\[
 K_N(v,s,m)= \Big[ \frac{g_{T_i}(\alpha_i^{\T}v+s) - g_{T_i}(\alpha_i^{\T}v+H_i(m))} {s-H_i(m)} \Big]_N,
\]Invoke $\cQ_i(s,m)$ in \eqref{Eq_define_series_cQ} to get
\begin{equation*}
    K_N(v,s,m) = \Big[\sum_{k\ge1}c_k(T_i)\sum_{\ell=0}^{k-1} (\alpha_i^{\T}v+s)^{k-1-\ell} (\alpha_i^{\T}v+H_i(m))^\ell\Big]_N
\end{equation*}
Thus we have $\|K_N\|_{\mathrm{coeff}} \le\sum_{w=0}^N[t^w]\{1+q(t)\}\le e^{CN}$ and $\|\ev\{K_N\}\|_{2q}\le e^{C(1+q)N^2}$.

Suppose $K(v,s,m)$ is the untruncated version of $K_N$. Then $g_{T_i}(\alpha_i^\T v+s) -g_{T_i}(\alpha_i^\T v+H_i(m)) =(s-H_i(m))K(v,s,m)$. Recall that $P_{i,N}(v,s)=\sum_{k=1}^Nc_k(T_i)(\alpha_i^{\T}v+s)^k$ and $G_{i,N}(v,m)=[g_{T_i}(\alpha_i^{\T}v+H_i(m))]_N$. Truncation to $N$ yields
\begin{equation*}
    [P_{i,N}(v,s)]_N-G_{i,N}(v,m)=[(s-[H_i(m)]_N)K_N]_N
\end{equation*}
H\"older's inequality gives
\begin{align*}
    &\,\|\ev\{[P_{i,N}(v,s)]_N-G_{i,N}(v,m)\}\|_q\\ \le &\, \big\|\delta_i-\ev\{[H_i(m)]_N\}\big\|_{2q} \big\|\ev(K_N)\big\|_{2q}\\ &\,+\|\ev \{(s-[H_i(m)]_N)K_N-[(s-[H_i(m)]_N)K_N]_N\}\|_q\\
 \overset{(1)}{\le} &\, e^{C(1+q)N^2}p^{-(N+1)/2}+ \big\{1 + \|[H_i(m)]_N\|_{\rm coeff}\big\}\|K_N\|_{\rm coeff}e^{C'(1+q)N^2}p^{-(N+1)/2}\\\overset{(2)}{\le }&\, e^{C''(1+q)N^2}p^{-(N+1)/2}
\end{align*}
Here, in $(1)$, the first term is bounded by \eqref{eq_Herror} and \eqref{eq_discarded}, and the second part is due to truncated terms have weights between $N+1$ and $2N$ and coefficient norm at most $e^{CN}$; $(2)$ is due to $\|{[H_i(m)]_N}\|_{\mathrm{coeff}} \le e^{CN}$ and $\|K_N\|_{\mathrm{coeff}} \le e^{CN}$.
Therefore, we arrive at
\begin{equation}\label{eq_factorerror}\begin{aligned}
    \|g_{T_i}(\alpha_i^{\T}V+\delta_i)-\ev\{G_{i,N}(v,m)\}\|_q \le&\,\|g_{T_i}(u)-\ev\{P_{i,N}(v,s)\}\|_q\\ &\,+\|\ev\{P_{i,N}(v,s)-[P_{i,N}(v,s)]_N\}\|_q \\&\,+\|\ev\{[P_{i,N}(v,s)]_N-G_{i,N}(v,m)\}\|_q\\ \le &\, e^{C(1+q)N^2}p^{-(N+1)/2}.
\end{aligned}\end{equation}where the last terms on the right side  are bounded by \eqref{eq_discard_coro}.

Finally, similar to \eqref{eq_taylor_remiander}, we have by $|g_t(u)|\le e^{H|u|}|u|$ that 
\begin{align*}
\|g_{T_i}(\alpha_i^{\T}V+\delta_i)\|_q \le  \|u\|_{2q}(\EE e^{2qH|u|})^{1/(2q)} \le  e^{C(1+2q)}p^{-1/2} e^{C\{1+(2qH)^2\}/(2q)} \le e^{C(1+q)}p^{-1/2}.
\end{align*}
\end{proof}

\subsection{Proof of Lemma~\ref{lem_design_alignment_sub}}\label{app_prove_lem_design_alignment_sub}

Write $\Delta=\max_{i,j\in[p]}|\alpha_{1i}^{\T}\alpha_{1j}
-\alpha_{0i}^{\T}\alpha_{0j}|$. Recall that Lemma~\ref{lem_design_alignment} finds a reference set $\mathcal I$ of size $d$ for which $\sigma_{\min}(B_0)\ge c_1$, where $B_r=(A_r)_{\mathcal I,}$.
Set $Q=B_1^{\T}(B_0^{-1})^{\T}$, and we  have
\[
 \|Q^{\T}Q-I_d\|
 =\|B_0^{-1}(B_1B_1^{\T}-B_0B_0^{\T})(B_0^{-1})^{\T}\|
 \le dc_1^{-2}\Delta.
\]
When $dc_1^{-2}\Delta\le1/2$, we know $O=Q(Q^{\T}Q)^{-1/2}$ is orthogonal and satisfies $\|O-Q\|\le3 dc_1^{-2}\Delta$ by the spectral theorem.
Since $\max_{j\in[p]}\|B_1\alpha_{1j}-B_0\alpha_{0j}\|\le d^{1/2}\,\Delta$, we know that for this fixed $O$, the following holds
\begin{equation*}
 \norm{\alpha_{1j}-O\alpha_{0j}} = \norm{O^{\T}\alpha_{1j}-\alpha_{0j}}\le\|O-Q\|\norm{\alpha_{1j}}
       +\|B_0^{-1}\|\norm{B_1\alpha_{1j}-B_0\alpha_{0j}}
 \le C\Delta,
\end{equation*}
uniformly in $j\in[p]$. When $dc_1^{-2}\Delta>1/2$, take $O=I_d$ and use
$\max_{j\in[p]}\norm{\alpha_{1j}-\alpha_{0j}}\le2M
\le4dMc_1^{-2}\Delta$. This proves \eqref{eq_alignment}.

\subsection{Proof of Lemma~\ref{lem_observable_corrections}}\label{app_prove_lem_observable_corrections}
 Throughout this proof, the series $H_i(m)$ and the forward polynomials $F_\gamma(m)$ are those of model zero.  We omit the arguments in the polynomials when they are clear from the context. Write $\Gamma_N=\{\gamma\in\mathbb N_0^d:2\le|\gamma|\le N\}$, and let $P^{[w]}=[P]_w-[P]_{w-1}$ denote the component of weight $w$ in a polynomial $P$.

By the construction of $F_{\gamma}$, we know $F_\gamma=m_\gamma+$ terms of weight greater than $|\gamma|$.
Consequently, $\mathcal U-\cI$ strictly increases the minimum weight of every nonconstant polynomial and maps any constant to zero. Thus
\begin{equation}\label{eq_correction_finite_inverse}
    \mathcal U^{-1}=\sum_{l=0}^N(\cI - \cU)^l\quad\text{ and }\quad \cU\sum_{l=0}^N(\cI - \cU)^l=I.
\end{equation}
Next, we calculate the coefficient bound regarding operator $\cU$, which is defined by
\begin{equation*}
    \|\cU\|_{\mathrm{coeff}} =\sup_{P\ne0}\frac{\|\cU P\|_{\mathrm{coeff}}}{\|P\|_{\mathrm{coeff}}},
\end{equation*}where the supremum is taken with respect to all symbolic polynomial $P(m) = \sum_{\eta}a_{\eta}(\mathbf{T})m^{\eta}$ where $\mathbf{T} - (T_i)_{i\in\cS}$.
For any nonconstant monomial $m^{\eta} = \prod_{a=1}^b m_{\gamma_a}$ of weights $u\le N$. Recall that $\cU(m)$ is defined by substituting $F_{\gamma}$ for each moment symbol $M_{\gamma}$ and truncating to $N$ order:
\begin{equation*}
    \cU(m^\eta)=\Big[\prod_{a=1}^bF_{\gamma_a}(m)\Big]_N.
\end{equation*}
Since $2\le |\gamma_a|\le N$, we know $m$ is of weight $u\ge 2b$. For any $u\le w\le N$, the weight-$w$ component in $\cU(m^\eta)$ is given by choosing some weight $t_a$ component from $F_{\gamma_a}$ such that $\sum_{a=1}^bt_a = w$. Since $\|F_{\gamma_a}^{[t]}\|_{\mathrm{coeff}}\le C^t$ by \eqref{eq_lemma_weighted_cal_cond_F_gamma}, the coefficient norm at weight $w$ is bounded by
\begin{equation*}
 \Big\|\Big[\prod_{a=1}^bF_{\gamma_a}\Big]_{w} - \Big[\prod_{a=1}^bF_{\gamma_a}\Big]_{w-1}\Big\|_{\mathrm{coeff}} \le \sum_{\substack{t_1+\cdots+t_b=w\\t_a\ge|\gamma_a|}}    \prod_{a=1}^b C^{t_a} \le C^w\binom{w-u+b-1}{b-1}\le C_1^w.
\end{equation*}
For constant monomial, this bound holds trivially. Therefore, it is easy to see by summing over all $w\le N$ that $\|\cU\|_{\rm coeff} \le e^{CN}$ and $\|\cU - \cI\|_{\rm coeff} \le e^{CN}$, which gives 
\begin{equation}\label{eq_correction_inverse_bound}
 \|\mathcal U^{-1}\|_{\mathrm{coeff}} \le\sum_{l=0}^N\|\mathcal U-\mathcal I\|_{\mathrm{coeff}}^l \le e^{CN^2}.
\end{equation}

By the definition of $R_i(m)$ and $R_{ij}(m)$ in \eqref{eq_corrections} and $\cU$, we know $[R_i\{F(m)\}]_N=[H_i(m)]_N$ and $ [R_{ij}\{F(m)\}]_N=\EE_{\mathrm f} \{(\alpha_{0i}^{\T}v)(\alpha_{0j}^{\T}v)\}$. Both of these polynomials start at weight 2, and thus by Lemma~\ref{lem_symbolic_bias} and \eqref{eq_correction_inverse_bound}, we know $ \|R_i\|_{\mathrm{coeff}}\le e^{CN^2}$ and $\|R_{ij}\|_{\mathrm{coeff}} \le e^{CN^2} $. Follow the same calculation above, one can also show that for any $2\le w\le N$,
\begin{equation*}
    \big\|\{R_i(m)\}_{[w]} - \{R_i(m)\}_{[w-1]}\big\|_{\mathrm{coeff}}\le e^{Cw^2}\quad\text{ and }\quad\big\|\{R_{ij}(m)\}_{[w]} - \{R_{ij}(m)\}_{[w-1]}\big\|_{\mathrm{coeff}} \le e^{Cw^2}.
\end{equation*}

We first prove \eqref{eq_common_corrections} under model zero, i.e., $r=0$. Recall that we have shown in \eqref{eq_forward} and \eqref{eq_O_gamma_q} that, for any $q\ge 1$, $\|O_\gamma-\ev_0(F_\gamma(m))\|_q \le e^{C(1+q)N^2}p^{-(N+1)/2}$ and $\|O_\gamma\|_q\le e^{C(1+q)|\gamma|^2}p^{-|\gamma|/2}$. Therefore by Lemma~\ref{lem_weighted_calculus}, we know
\begin{equation*}
 \Big\|\prod_{\gamma\in\Gamma_N}O_\gamma^{\eta_\gamma} -\ev_0\Big\{\big[\prod_{\gamma\in\Gamma_N}   F_\gamma(m)^{\eta_\gamma}\big]_N\Big\}\Big\|_q
 \le e^{C(1+q)N^2}p^{-(N+1)/2}.
\end{equation*}
The coefficient bound for $R_{i}(m)$ and $R_{ij}(m)$ therefore gives 
\begin{align}\label{eq_correction_zero_substitution_1}
 \|R_i(O)-\ev_0\{[R_i(F(m))]_N\}\|_q\le e^{C'(1+q)N^2}p^{-(N+1)/2},\\
 \|R_{ij}(O)-\ev_0\{[R_{ij}(F(m))]_N\}\|_q\le e^{C'(1+q)N^2}p^{-(N+1)/2}.\label{eq_correction_zero_substitution_2}\end{align}
Then for $R_i(O)-\delta_{ri}$, combining \eqref{eq_correction_zero_substitution_1}, Lemma~\ref{lem_symbolic_bias}, and $[R_i\{F(m)\}]_N=[H_i(m)]_N$, we obtain the bound as required in \eqref{eq_common_corrections} under $r=0$.
For $R_{ij}(O)-\EE_r\{(\alpha_{ri}^{\T}V_r) (\alpha_{rj}^{\T}V_r)\mid Y\}$, combining \eqref{eq_correction_zero_substitution_2} and $ [R_{ij}\{F(m)\}]_N=\EE_{\mathrm f} \{(\alpha_{0i}^{\T}v)(\alpha_{0j}^{\T}v)\}$ proves the required bound in \eqref{eq_common_corrections} under $r=0$.

Finally, we show \eqref{eq_common_corrections} under $r=1$. Following the definition of $F_{\gamma}(m)$ in \eqref{eq_forward_polynomial}, let 
\begin{equation}\label{eq_correction_forward_models}
 F_\gamma^{(r)}(m)= \Bigg[\sum_{h_1,\ldots,h_k=1}^d \Big\{\prod_{s=1}^k(B_s^{-1})_{\ell_s h_s}\Big\} \EE_{\mathrm f}\Big\{\prod_{s=1}^k  G_{j_{h_s}^{(s)},N}^{(r)}(v,m)\Big\}\Bigg]_N.
\end{equation} Here, $k= |\gamma|$; $B_s^{-1}$ are fixed at model zero; $G_{i,N}^{(r)}(v,m)$ are the polynomials  from Lemma~\ref{lem_symbolic_bias} using $\alpha_{ri}$.
Then $F_{\gamma}^{(0)}(m) = F_{\gamma}(m)$. Under model one, $F_\gamma^{(r)}(m)$ has leading weight-k component $\EE_{\mathrm f}\big[ \prod_{s=1}^k\{ \sum_{h=1}^d(B_s^{-1})_{\ell_s h}\alpha_{1,j_h^{(s)}}^{\T}v\}\big]$. Similar to \eqref{eq_lemma_weighted_cal_cond_F_gamma}, we know $\|[F_\gamma^{(1)}]_w-[F_\gamma^{(1)}]_{w-1}\|_{\mathrm{coeff}}\le C^w$ for $k\le w\le N$. In addition, to bound the coefficient in $F_\gamma^{(1)}(m)-F_\gamma^{(0)}(m)$, interpolate $\alpha_{i}$  by $\alpha_a(\lambda)=(1-\lambda)\alpha_{0a}+\lambda\alpha_{1a}$ and let $G_{a,N}(\lambda)$ be the polynomial from Lemma~\ref{lem_symbolic_bias} with $\alpha_a(\lambda)$. Define $F_\gamma(\lambda)$ by \eqref{eq_correction_forward_models} with $G_{a,N}^{(r)}$ replaced by $G_{a,N}(\lambda)$. Then
\begin{align*}
    &\,\big\|\{F_\gamma^{(1)}(m)-F_\gamma^{(0)}(m)\}_{[w]} - \{F_\gamma^{(1)}(m)-F_\gamma^{(0)}(m)\}_{[w-1]}\big\|_{\mathrm{coeff}}
 \\\le &\,\int_0^1 \Big\|\big\{\frac{d}{d\lambda}F_\gamma(\lambda)\big\}_{[w]} - \big\{\frac{d}{d\lambda}F_\gamma(\lambda)\big\}_{[w-1]} \Big\|_{\mathrm{coeff}}\,d\lambda\\
 \le &\,C_B^k\max_{h_1,\dots,h_k\in[d]}\sum_{s=1}^k \sum_{\substack{t_1+\cdots+t_k=w\\t_1,\ldots,t_k\ge1}}{\Big\|\frac{d}{d\lambda}
\big[\{G_{a_s,N}(\lambda)\}_{[t_s]} - \{G_{a_s,N}(\lambda)\}_{[t_s-1]}\big]\Big\|_{\mathrm{coeff}}}
\\&\,\qquad\times \prod_{b\ne s} {\big\|\{G_{a_b,N}(\lambda)\}_{[t_b]} - \{G_{a_b,N}(\lambda)\}_{[t_b-1]}\}\big\|_{\mathrm{coeff}}}.
\end{align*}
Here, $a_s = j_{h_s}^{(s)}$ and the second inequality uses $\prod_{s=1}^k\sum_{h=1}^d|(B_s^{-1})_{\ell_s h}|\le C_B^k$. Since $\norm{\alpha_a(\lambda)}\le M$, with the coefficient and derivative bound in Lemma~\ref{lem_symbolic_bias}, one can show $\|G_{a,N}(\lambda)^{[t]}\|_{\mathrm{coeff}}\le C^t$. The chain rules gives 
\begin{align*}
 &\,\Big\|\frac{d}{d\lambda}
\big[\{G_{a_s,N}(\lambda)\}_{[t_s]} - \{G_{a_s,N}(\lambda)\}_{[t_s-1]}\big]\Big\|_{\mathrm{coeff}}
\\ \le&\, \sum_{\ell=1}^d|\alpha_{1a_s,\ell}-\alpha_{0a_s,\ell}| \left\|\partial_{\alpha_{a_s\ell}}\big[\{G_{a_s,N}(\lambda)\}_{[t_s]} - \{G_{a_s,N}(\lambda)\}_{[t_s-1]}\big]\right\|_{\mathrm{coeff}}
 \\ \le&\, D_p C^{t_s},
\end{align*}
where we denote $D_p = D_p(\theta_0,\theta_1)$. Therefore, 
\begin{equation}\label{eq_correction_forward_difference}
    \big\|\{F_\gamma^{(1)}(m)-F_\gamma^{(0)}(m)\}_{[w]} - \{F_\gamma^{(1)}(m)-F_\gamma^{(0)}(m)\}_{[w-1]}\big\|_{\mathrm{coeff}}\le D_pC_1^w.
\end{equation}

Write $P_{ij}^{(r)}(m)=\EE_{\mathrm f} \{(\alpha_{ri}^{\T}v)(\alpha_{rj}^{\T}v)\}$ for short. Then by definition, 
\begin{align*}
[\{R_i(F^{(1)}(m)\}]_N-[H_i^{(1)}(m)]_N  = &\, [\{R_i(F^{(1)}(m)\}-\{R_i(F^{(0)}(m)\}]_N -[H_i^{(1)}(m)-H_i^{(0)}(m)]_N,\\
[R_{ij}\{F^{(1)}(m)\}]_N-P_{ij}^{(1)}(m) = &\, [R_{ij}\{F^{(1)}(m)\}-R_{ij}\{F^{(0)}(m)\}]_N -\{P_{ij}^{(1)}(m)-P_{ij}^{(0)}(m)\}.
\end{align*}
We first show that the evaluated value of the two polynomials are bounded by $C_qD_p/p$, and then establish that replacing $m_{\gamma}$ by $O_{\gamma}$ introduces errors bounded by $e^{C(1+q)N^2}p^{-(N+1)/2}$.

First, note that $\prod_{a=1}^bF_{\gamma_a}^{(1)}- \prod_{a=1}^bF_{\gamma_a}^{(0)} =\sum_{a=1}^b(F_{\gamma_a}^{(1)}-F_{\gamma_a}^{(0)}) \prod_{c<a}F_{\gamma_c}^{(1)} \prod_{c>a}F_{\gamma_c}^{(0)}$.  Combining $\|\{R_i(m)\}_{[w]} - \{R_i(m)\}_{[w-1]} \|_{\mathrm{coeff}}\le e^{Cw^2}$, $\|\{R_{ij}(m)\}_{[w]} - \{R_{ij}(m)\}_{[w-1]}\|_{\mathrm{coeff}} \le e^{Cw^2}$, and the coefficient bounds for $F_{\gamma}^{(r)}$, $F_{\gamma}^{(1)} - F_{\gamma}^{(0)}$, we know the coefficients of weight $w$ in $\prod_{a=1}^bF_{\gamma_a}^{(1)}- \prod_{a=1}^bF_{\gamma_a}^{(0)} $ can be bounded by $D_pe^{Cw^2}$, as only components of weight at most $w$ are involved. With Lemma~\ref{lem_symbolic_bias} and $P_{ij}^{(r)}$ having leading term at least weight $2$, we can also show that the weight-$w$ components' coefficients in $\{R_i(F^{(1)}(m)\}]_N-[H_i^{(1)}(m)]_N$ and $[R_{ij}\{F^{(1)}(m)\}]_N-P_{ij}^{(1)}(m)$ can be bounded by $D_p e^{Cw^2}$ uniformly for $2\le w\le N$.
With \eqref{eq_weighted}, we conclude that 
\begin{equation}\label{eq_correction_stability_sum}\begin{aligned}
    &\,\Big\|{\rm ev}_1\Big(\big[R_i\{F^{(1)}(m)\}\big]_N - \big[H_i^{(1)}\big]_N\Big)\Big\|_q + \Big\|{\rm ev}_1\Big(\big[R_{ij}\{F^{(1)}(m)\}\big]_N - P_{ij}^{(1)}(m)\Big)\Big\|_q\\ \le &\, \sum_{w=2}^ND_p e^{Cw^2}e^{Cqw^2}p^{-w/2} \le C' \frac{D_p}{p}.\end{aligned}
\end{equation}

Finally, to replace $m_{\gamma}$ with $O_{\gamma}$, we note that since the conditional independent holds also under model one, \eqref{eq_observable} still hold with $r=1$, and therefore with Lemmas~\ref{lem_symbolic_bias} and~\ref{lem_weighted_calculus}, we know 
\begin{equation*}
 \|O_\gamma-\ev_1(F_\gamma^{(1)})\|_q
 \le e^{C(1+q)N^2}p^{-(N+1)/2}.
\end{equation*}
Each $F_\gamma^{(1)}$ starts at weight $|\gamma|$ and has the coefficient bounds established above, and \eqref{eq_forward}, \eqref{eq_O_gamma_q} analogously hold for $F_\gamma^{(1)}$.Repeating the monomial substitution argument leading to \eqref{eq_correction_zero_substitution_1}, \eqref{eq_correction_zero_substitution_2} yields 
\begin{align*}
 \|R_i(O)-\ev_1\{[R_i\{F^{(1)}(m)\}]_N\}\|_q\le e^{C(1+q)N^2}p^{-(N+1)/2};\\
 \|R_{ij}(O)-\ev_1\{[R_{ij}\{F^{(1)}(m)\}]_N\}\|_q\le e^{C(1+q)N^2}p^{-(N+1)/2};
\end{align*}
For $R_i$, combine this with \eqref{eq_correction_stability_sum} and the approximation of $\delta_{1i}$ by $\ev_1\{[H_i^{(1)}]_N\}$ from Lemma~\ref{lem_symbolic_bias} yields the required bound under $r=1$.
For $R_{ij}$, recall that $\ev_1(P_{ij}^{(1)}) =\EE_1\{(\alpha_{1i}^{\T}V_1)(\alpha_{1j}^{\T}V_1)\mid Y\}$. Use \eqref{eq_correction_stability_sum} and we prove \eqref{eq_common_corrections} under $r=1$.




\section{Additional Simulation Details}\label{app_sim}

This section describes how the quantities reported in Section~\ref{sec_simulation} are computed. The computation has two steps. First, for each candidate $\tilde\theta$, we approximate the discrepancy $\Delta_p(\tilde\theta)$ in~\eqref{eq_sim_delta} by a penalized nonnegative least-squares fit over the masses that $\tilde G$ places on the grid, with the constraints of $\mathcal G_{{\rm grid},d}$ imposed through a penalty. To distinguish exact equivalence from very small positive discrepancies under the cases with $d=1$, $p=7,8,9$ and $d=3$, $p=9$, we check directly by verifying whether admissible grid masses can reproduce $P_{\theta_0,G_0}$, while accounting for rounding error. Second, along each perturbation path $\tilde\theta_t$, we classify each evaluated value of \(t\) as numerically equivalent or nonequivalent using these computations, and locate the transition point $t_p$ by a coarse search followed by bisection in $\log_{10}t$. We describe the two steps in turn.

Let $z_1,\ldots,z_J$ be the points of the grid supporting $\mathcal G_{{\rm grid},d}$, let $p_0=p_{\theta_0,G_0}$, and for a candidate $\tilde\theta$ define the $2^p\times J$ matrix
$B_{\tilde\theta}=\{\operatorname{pr}_{\tilde\theta}(X=x\mid Z=z_\ell)\}_{x\in\{0,1\}^p,\ell\in[J]}$. A latent distribution with nonnegative masses $w=(w_1,\ldots,w_J)^\T$ on the grid then has response probability vector $B_{\tilde\theta}w$. Membership in $\mathcal G_{{\rm grid},d}$ imposes the linear constraints
\[
\sum_{\ell=1}^Jw_\ell=1,\quad
\sum_{\ell=1}^Jw_\ell z_\ell=0,\quad
\sum_{\ell=1}^Jw_\ell z_\ell z_\ell^{\T}=I_d,\quad
\sum_{\ell=1}^Jw_\ell\exp(\|z_\ell\|^2/4)+s=2^{(d+1)/2},\]
where the first three constraints encode~\eqref{eq_unv_normalize}, and the slack variable $s\ge0$ turns the tail bound into an equality. We write these $R=2+d+d(d+1)/2$ constraints compactly as $Cw+s e_R=b$, where \(e_R=(0,\ldots,0,1)^\T\in\mathbb R^R\) is the last standard basis vector. To put the constraints on a common scale, let $D=\operatorname{diag}(a_1^{-1},\ldots,a_R^{-1})$, where $a_r$ is the largest absolute coefficient in row $r$ of $[C,e_R]$. We then replace the constrained minimization in~\eqref{eq_sim_delta} by the penalized problem
\begin{equation}\label{eq_sim_penalty}
\min_{w\in[0,\infty)^J,\ s\ge0}
\left\{\frac{\|B_{\tilde\theta}w-p_0\|^2}{\|p_0\|^2}
+\lambda^2\|D(Cw+s e_R-b)\|^2\right\}.
\end{equation}
Because the objective is a sum of squares of affine functions of $(w,s)$ and the only remaining constraints are nonnegativity, problem~\eqref{eq_sim_penalty} is a nonnegative least-squares problem. We solve it by the active-set method of~\citet{lawson1995solving}, with $\lambda=10^{15}$ for $d=1$ and $\lambda=10^{10}$ for $d=3$. All fits use the fixed grids described in Section~\ref{sec_simulation} and are computed in $1500$-bit arithmetic for $d=1$ and $600$-bit arithmetic for $d=3$.
Let $(\widehat w,\widehat s)$ denote the masses and slack returned by the solver, and define
\begin{equation}\label{eq_sim_fitted_delta}
\widehat\Delta_p(\tilde\theta)
=\frac{\|B_{\tilde\theta}\widehat w-p_0\|}{\|p_0\|}.
\end{equation}
Since a finite penalty does not enforce the constraints exactly, $\widehat\Delta_p$ is a numerical approximation to $\Delta_p$.


For the penalized least-squares calculations, we classify an evaluated value $t$ as numerically equivalent if $\widehat\Delta_p(\tilde\theta_t)^2\le\tau$ and as nonequivalent otherwise, with $\tau\approx10^{-301}$ for $d=1$ and $\tau\approx10^{-120}$ for $d=3$. For $d=1$, $p=7,8,9$ and $d=3$, $p=9$, we instead classify $t$ as numerically equivalent when we verify that nonnegative grid masses satisfy all the constraints, and as nonequivalent when we verify that no such masses exist. These checks account for rounding error.

To estimate $t_p$, we search over $a=\log_{10}t$ in two stages. In the coarse stage, the search starts at $a=0$ and decreases $a$ until it finds a numerically equivalent point or reaches a prescribed lower search limit. Once two consecutive points receive different classifications, write their log values as $a_{\rm L}<a_{\rm U}$, with the lower point numerically equivalent and the upper point nonequivalent. In the refinement stage, each bisection step evaluates the midpoint $a_{\rm M}=(a_{\rm L}+a_{\rm U})/2$, corresponding to $t=10^{a_{\rm M}}$, and replaces $a_{\rm L}$ by $a_{\rm M}$ if this point is numerically equivalent and $a_{\rm U}$ by $a_{\rm M}$ otherwise. Bisection continues until $a_{\rm U}-a_{\rm L}\le0.1$.

\end{document}